\pdfoutput=1
\documentclass{siamart1116}

\usepackage{lipsum}
\usepackage{amsfonts}
\usepackage{graphicx}
\usepackage{epstopdf}
\usepackage{algorithmic}
\ifpdf
  \DeclareGraphicsExtensions{-eps-converted-to,.pdf,.png,.jpg}
\else
  \DeclareGraphicsExtensions{-eps-converted-to}
\fi

\numberwithin{theorem}{section}

\newcommand{\U}{\mathbf{U}}

\newcommand{\Q}{\mathbf{Q}} 
\newcommand{\q}{\mathbf{q}} 
\newcommand{\f}{\mathbf{f}} 
\newcommand{\g}{\mathbf{g}}
\newcommand{\h}{\mathbf{h}} 
 
\newcommand{\B}{\mathbf{B}} 
\newcommand{\F}{\mathbf{F}} 
\newcommand{\G}{\mathbf{G}} 
\newcommand{\x}{\mathbf{x}} 
\newcommand{\xxi}{\boldsymbol{\xi}}

\newcommand{\w}{\mathbf{w}}

\renewcommand{\v}{\mathbf{v}}

\newcommand{\quotew}[1]{``#1''}

\newcommand{\Popt}{\mathbf{P}_{\text{\sf opt}}}
\DeclareMathOperator{\argmin}{argmin}

\newcommand{\TheTitle}{CWENO schemes for conservation laws on unstructured meshes} 
\newcommand{\TheAuthors}{M. Dumbser, W. Boscheri, M. Semplice, and G. Russo}

\headers{\TheTitle}{\TheAuthors}

\title{{Central WENO schemes for hyperbolic conservation laws on fixed and moving unstructured meshes}
\thanks{Submitted to the editors DATE.} 
}

\author{
  Michael Dumbser\thanks{University of Trento, Laboratory of Applied Mathematics, Via Mesiano 77, I-38123 Trento, Italy 
    (\email{michael.dumbser@unitn.it},\email{walter.boscheri@unitn.it}).} 
  \and
  Walter Boscheri\footnotemark[2] 
  \and
  Matteo Semplice\thanks{Dipartmento di Matematica, Universit\`a di Torino, Via C. Alberto 10, I-10123 Torino, Italy 
	  (\email{matteo.semplice@unito.it}).} 
	\and 
	Giovanni Russo\thanks{Department of Mathematics and Informatics, University of Catania, Viale A. Doria 6, I-95125 Catania, Italy 
	(\email{russo@dmi.unict.it}).} 
}

\usepackage{amsopn}

\ifpdf
\hypersetup{
  pdftitle={\TheTitle},
  pdfauthor={\TheAuthors}
}
\fi

\begin{document}

\maketitle

\begin{abstract}
We present a novel arbitrary high order accurate central WENO spatial reconstruction procedure (CWENO) for the solution of nonlinear systems of hyperbolic conservation  
laws on fixed and moving unstructured simplex meshes in two and three space dimensions. Starting from the given cell averages of a function on a triangular or tetrahedral 
control volume and its neighbors, the nonlinear CWENO reconstruction yields a high order accurate and essentially non-oscillatory polynomial that is defined everywhere 
in the cell. 
Compared to other WENO schemes on unstructured meshes, the total stencil size is the minimum possible one, as in classical point-wise WENO schemes of Jiang and Shu. However,
the linear weights can be chosen arbitrarily, which makes the practical implementation on general unstructured meshes particularly simple. 
We make use of the piecewise polynomials generated by the CWENO reconstruction operator inside the framework of fully discrete and high order accurate one-step ADER finite 
volume schemes on fixed Eulerian grids as well as on moving Arbitrary-Lagrangian-Eulerian (ALE) meshes. The computational efficiency of the high order finite volume schemes 
based on the new CWENO reconstruction is tested on several two- and three-dimensional benchmark problems for the compressible Euler and MHD equations and is found to be more 
efficient in terms of memory consumption and computational efficiency with respect to classical WENO reconstruction schemes on unstructured meshes.  
We also provide evidence that the new algorithm is suitable for implementation on massively parallel distributed memory supercomputers, showing two numerical examples 
run with more than one billion degrees of freedom in space.  
\end{abstract}

\begin{keywords}
central WENO reconstruction (CWENO); 
fixed and moving unstructured meshes; 
fully-discrete one-step ADER approach; 
high order in space and time;  
Arbitrary-Lagrangian-Eulerian (ALE) finite volume schemes; 
hyperbolic conservation laws in multiple space dimensions; 
large scale parallel HPC computations;  
\end{keywords}

\begin{AMS}
  65M08, 76N15, 35L60
\end{AMS}

\section{Introduction}
\label{sec.introduction}

Finite volume schemes are widely employed for the numerical solution of nonlinear systems of hyperbolic conservation laws and are thus very important in many application fields in science and 
engineering. In a finite volume scheme one stores and evolves only the cell averages of the conserved variables on a computational grid. The advantages of using higher order schemes on a coarser 
mesh rather than a first order one on a very fine mesh in terms of accuracy and computational efficiency are unquestionable. However, in order to advance the solution from time $t^n$ to $t^{n+1}$ 
with high order of accuracy, more spatial information are needed than just the mere cell averages. Therefore, it is necessary to employ a {\em reconstruction} or \textit{recovery} procedure, 
that can produce high order piecewise polynomials from the known cell averages in an appropriate neighborhood of the control volume under consideration, the so-called reconstruction stencil. 

Many reconstruction strategies of this kind have been developed and one of the most successful ones is the Weighted Essentially Non Oscillatory reconstruction (WENO), which was introduced in the seminal  papers \cite{shu_efficient_weno,balsarashu} and is described together with its many improvements and variants in the reviews \cite{Shu97,Shu:2009:WENOreview}. The WENO procedure is targeted at providing 
high-order non-oscillatory point values of the conserved variables at cell boundaries. Its implementation is rather simple and very efficient in the one-dimensional case or on multi-dimensional 
uniform and adaptive Cartesian meshes, see \cite{TitarevToroWENO3D,SubcellWENO,AMR3DCL,Buchmuller2014,Buchmuller2015}. On unstructured meshes, however, reconstruction becomes much more cumbersome, 
and early pioneering works are those of Barth et al. \cite{BarthJespersen,barthlsq}, Abgrall \cite{abgrall_eno} and others \cite{goochaltena}. Several extensions of the popular WENO reconstruction to unstructured meshes 
have been introduced in two and three space dimensions, see for example \cite{friedrich,HuShuTri,kaeserjcp,Dumbser2007204,Dumbser2007693,MixedWENO2D,MixedWENO3D,ZhangShu3D}. 

The essential idea of the WENO reconstruction is to reproduce the point values of a high order central interpolating polynomial by means of convex combinations of the point values of lower degree polynomials having smaller and directionally biased stencils. The coefficients of the convex combination have some optimal values (called {\em linear} or {\em optimal weights}) that are determined by satisfying some accuracy requirements. The coefficients actually employed in the reconstruction (called {\em nonlinear weights}) are typically derived from the optimal values by a nonlinear procedure whose task is to discard any information that might lead to oscillatory polynomials.

The main problem in applying the WENO idea to unstructured meshes or to reconstructions at points inside the computational cells comes from the definition of the linear weights. In fact these must satisfy accuracy requirements that depend on the location of the reconstruction point and on the size and relative location of the neighboring cells. The first difficulty is particularly evident in the fact that even in one space dimension the linear weights for the reconstruction at the center of the cell either do not exist (for accuracy $3+4k$) or involve non-positive numbers (for accuracy $5+4k$), see \cite{QiuShu:02}. The second difficulty arises from the very complicated formulae which define the linear weights in one-dimensional non-uniform meshes in terms of the ratio of the cell sizes (see e.g. \cite{WangFengSpiteri,PS:shentropy}) and even more in the higher-dimensional setting, where they additionally also depend on the disposition of the neighboring cells (see e.g. \cite{HuShuTri,Dumbser2007693}). 

The central WENO reconstruction (CWENO) was originally introduced in the one-dimensional context by Levy, Puppo and Russo \cite{LPR:99} in order to obtain a third-order accurate reconstruction at 
the cell center that was needed for the construction of a central finite volume scheme and that could not be provided by the classical WENO3 scheme of Jiang and Shu. The technique was later 
extended to fifth order in \cite{Capdeville:08} and subsequently the properties of the third order versions were studied in detail on uniform and non-uniform meshes in \cite{Kolb:2014} and 
\cite{CS:epsweno}. 

The CWENO procedure is again based on the idea of reproducing a high-order cell centered (central) polynomial with a convex combination of other polynomials and on the use of oscillation indicators to prevent the onset of spurious oscillations, but it differs from WENO in a seemingly minor but substantial aspect. Indeed, the linear weights for CWENO do not need to satisfy any accuracy requirements and can thus be chosen independently from the local mesh topology (which dramatically simplifies the multi-dimensional extension to unstructured meshes) and also independently from the location of the reconstruction point inside the cell. In turn, this means that CWENO employs one single set of linear weights and thus one single set of nonlinear weights that are valid for any point in the cell. This is of course equivalent to saying that the output of the CWENO reconstruction is \textit{not} given by some point values, but instead CWENO provides an entire \textit{reconstruction polynomial} that is defined everywhere inside a cell. Such polynomials can then be evaluated where needed within the cell, which is particularly important for hyperbolic PDE with source terms or non-conservative products, see
\cite{Pares2006,Castro2006}.
A very general framework for CWENO reconstructions that highlights these features has been recently presented in \cite{CWENOframework} and the reconstruction is exploited there to obtain high order accurate data representation also at inner quadrature points for the discretization of source terms, also in the context of well-balanced schemes.

To the best of our knowledge, the only realizations of the CWENO procedure in more than one space dimension have been presented in \cite{LPR:00,LPR:02} for Cartesian meshes and in \cite{SCR:CWENOquadtree} on two-dimensional quadrangular meshes, locally refined in a non-conforming fashion in a quad-tree type grid. In particular, \cite{SCR:CWENOquadtree} is the first example where the the independence of the linear weights on the mesh topology was exploited to obtain a very simple and robust reconstruction and the global definition of the piecewise polynomial reconstruction was used for the computation of a quadrature appearing in the error indicator (the numerical entropy production, \cite{PS:entropy,PS:shentropy}) and for the quadratures that compute sub-cell averages during mesh refinement.

Although the need for a reconstruction procedure is present in all time-marching schemes, in particular for the ADER approach 
\cite{toro3,toro10,DumbserEnauxToro,Dumbser2008,CastroToro,GoetzIske} that relies on the approximate solution of generalized Riemann problems (GRP) 
one needs a reconstruction procedure that is able to provide a globally defined piecewise polynomial reconstruction, 
with the entire polynomial defined everywhere in each computational cell. This requirement rules out the classical pointwise WENO schemes \cite{shu_efficient_weno,balsarashu,HuShuTri,ZhangShu3D} and leaves one with the choice of the ENO techniques \cite{eno,abgrall_eno} or the WENO procedure of \cite{friedrich,Dumbser2007204,Dumbser2007693}. Both are characterized by having very large stencils since all polynomials involved, even the directionally biased ones, should have the correct order of accuracy. This in turn creates difficulties for example when more than one discontinuity is located close to the current computational cell and suggests the use of many directionally biased polynomials in order to prevent the generation of spurious oscillations.

In this paper we introduce a novel CWENO reconstruction procedure that is, at least in principle, arbitrary high order accurate in space and which reconstructs piecewise polynomials from given cell  averages on conforming two- and three-dimensional unstructured simplex meshes. The reconstruction is tested relying on the ADER approach for hyperbolic systems of conservation laws, both in the Eulerian setting \cite{Dumbser2008} and in the Arbitrary-Lagrangian-Eulerian (ALE) framework \cite{Lagrange3D} that includes also the mesh motion from the current time step to the next one. The efficiency of the new reconstruction is studied by comparing the new ADER-CWENO schemes with the classical ADER-WENO methods based on the WENO reconstruction detailed in 
\cite{Dumbser2007204,Dumbser2007693}. In order to demonstrate that the new method is also well suited for the implementation on \textbf{massively parallel} distributed memory \textbf{supercomputers}, 
we also show two computational examples that involve hundreds of millions of elements and a total of more than \textbf{one billion} degrees of freedom to represent the piecewise polynomial 
reconstruction in each time step. To our knowledge, these are the largest simulations ever run so far with WENO schemes on unstructured meshes. 

The rest of the paper is organized as follows. In Section \ref{sec.numethod} the reconstruction and the fully discrete numerical scheme are described. In Section \ref{sec.validation} numerous two- and three-dimensional tests for conservation laws are considered in both the Eulerian and the ALE setting, i.e. on fixed and moving meshes, respectively. Finally, in Section \ref{sec.concl} we summarize the main findings of this work and give some outlook to future developments.

\section{Numerical method}
\label{sec.numethod} 

In this paper we consider time-dependent nonlinear hyperbolic systems of conservation laws that can be cast into the form
\begin{equation}
\label{eqn.PDE}
 \frac{\partial \Q}{\partial t} + \nabla \cdot \F(\Q) = \mathbf{0}, \qquad \x \in \Omega(t) \subset \mathbb{R}^d, \quad t \in \mathbb{R}_0^+, \quad \Q \in \Omega_{\Q} \subset \mathbb{R}^\nu,     
\end{equation} 
where $d\in[2,3]$ is the space dimension. We present the entire algorithm considering $d=3$, since the two-dimensional version can be easily derived from the former. In the governing PDE \eqref{eqn.PDE} $\x=(x,y,z)$ denotes the spatial coordinate vector, $\Q$ is the vector of conserved variables defined in the space of the admissible states $\Omega_{\Q}\subset \mathbb{R}^\nu $ and $\F(\Q)=\left( \f(\Q),\g(\Q),\h(\Q) \right)$ represents the nonlinear flux tensor. 

The computational domain $\Omega(t) \in \mathbb{R}^d$ can in general be time-dependent and is covered by a total number $N_E$ of non-overlapping control volumes $T_i^n$, whose union
\begin{equation}
\mathcal{T}_{\Omega}^n = \bigcup \limits_{i=1}^{N_E}{T^n_i}
\label{trian}
\end{equation}
is referred to as the \textit{current mesh configuration}.

\begin{figure}[!htbp]
\centering
\includegraphics[width=0.9\textwidth]{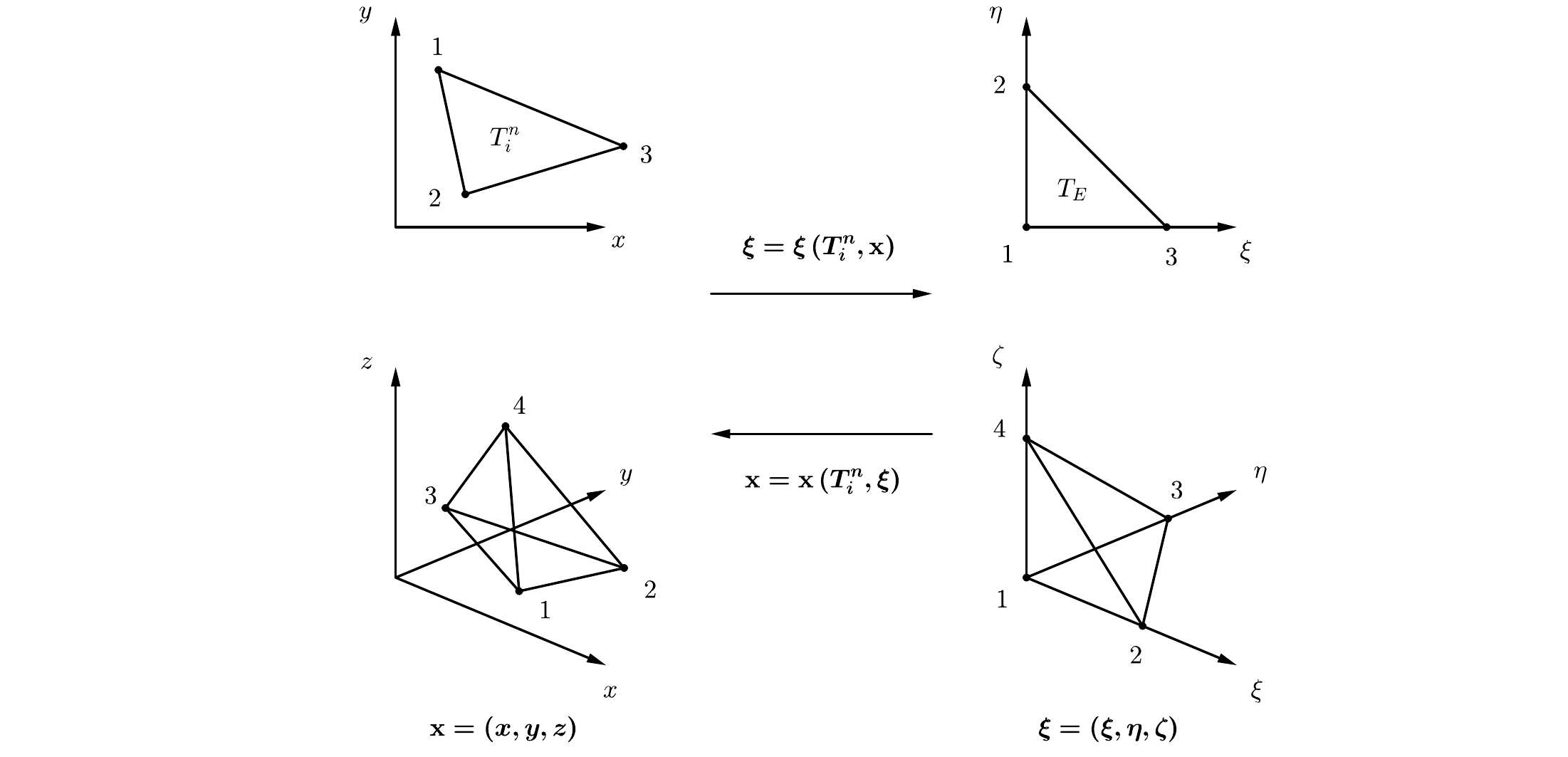}  
\caption{Spatial mapping from the physical element $T^n_i$ given in physical coordinates $\mathbf{x}$ to the unit reference element $T_E$ in $\boldsymbol{\xi}$ for triangles (top) 
and tetrahedra (bottom). Vertices are numbered according to a local connectivity.} 
\label{fig.refSystem}
\end{figure}

In the ALE context, due to the mesh motion, the domain is time-dependent, i.e. $\Omega=\Omega(t)$, but the mesh topology remains always fixed. Only the vertices of the computational grid are moved, as fully explained in Section \ref{sec.meshMot}, so that $T_i^n$ will denote the element configuration at time $t^n$. For this reason, it is convenient to adopt the following linear mapping from the physical element $T^n_i$ to a reference element $T_E$:
\begin{equation} 
 \mathbf{x} = \mathbf{X}^n_{1,i} + 
\left( \mathbf{X}^n_{2,i} - \mathbf{X}^n_{1,i} \right) \xi + 
\left( \mathbf{X}^n_{3,i} - \mathbf{X}^n_{1,i} \right) \eta + 
\left( \mathbf{X}^n_{4,i} - \mathbf{X}^n_{1,i} \right) \zeta,
 \label{xietaTransf} 
\end{equation}
where $\boldsymbol{\xi}=(\xi,\eta,\zeta)$ is the reference position vector and $\mathbf{X}^n_{k,i}=(X^n_{k,i},Y^n_{k,i},Z^n_{k,i})$ represents the vector of physical spatial coordinates of the $k$-th vertex of element $T^n_i$. Figure \ref{fig.refSystem} depicts the reference element $T_E$ which is the unit triangle in 2D, defined by nodes $\boldsymbol{\xi}_{e,1}=(\xi_{e,1},\eta_{e,1})=(0,0)$, $\boldsymbol{\xi}_{e,2}=(\xi_{e,2},\eta_{e,2})=(1,0)$ and $\boldsymbol{\xi}_{e,3}=(\xi_{e,3},\eta_{e,3})=(0,1)$, or the unit tetrahedron in 3D with vertices $\boldsymbol{\xi}_{e,1}=(\xi_{e,1},\eta_{e,1},\zeta_{e,1})=(0,0,0)$, $\boldsymbol{\xi}_{e,2}=(\xi_{e,2},\eta_{e,2},\zeta_{e,2})=(1,0,0)$, $\boldsymbol{\xi}_{e,3}=(\xi_{e,3},\eta_{e,3},\zeta_{e,3})=(0,1,0)$ and $\boldsymbol{\xi}_{e,4}=(\xi_{e,4},\eta_{e,4},\zeta_{e,4})=(0,0,1)$.

The cell averages of the conserved quantities $\Q(\x,t)$ are denoted by
\begin{equation}
  \Q_i^n = \frac{1}{|T_i^n|} \int_{T^n_i} \Q(\mathbf{x},t^n) d\x,     
 \label{eqn.cellaverage}
\end{equation}  
with $|T_i^n|$ representing the volume of cell $T_i^n$.

\subsection{Polynomial CWENO reconstruction} 
\label{ssec.cweno} 
The task of the reconstruction procedure is to compute a high order non-oscillatory polynomial representation $\w_i(\x,t^n)$ of the solution $\Q_i(\x,t^n)$ for each computational cell $T_i^n$. The degree $M$ of the polynomial can be chosen arbitrarily and it provides a nominal spatial order of accuracy of $M+1$. The total number of unknown degrees of freedom $\mathcal{M}$ of $\w_i(\x,t^n)$ is given by
\begin{equation}
  \mathcal{M}(M,d)=\frac{1}{d!} \prod \limits_{l=1}^{d} (M+l).
\label{eqn.cwenoM}
\end{equation}

The requirements of high order of accuracy and non-oscillatory behavior are contradictory due to the Godunov theorem \cite{godunov:linear}, since one needs a large stencil centered in $T_i^n$ in order to achieve high accuracy,  but this choice produces oscillations close to discontinuities, the well-known Gibbs' phenomenon. In order to fulfill both requirements, here we propose a CWENO reconstruction strategy that can be cast in the framework introduced in \cite{CWENOframework}. The procedure adopted in this work is inspired by the one already employed in \cite{SCR:CWENOquadtree} for two-dimensional unstructured meshes of quad-tree type.

The reconstruction starts from the computation of a polynomial $\Popt$ of degree $M$. In order to define $\Popt$ in a robust manner, following \cite{barthlsq,goochaltena,kaeserjcp,SCR:CWENOquadtree}, we consider a stencil $\mathcal{S}_i^0$ with $n_e=d \cdot \mathcal{M}(M,d)$ cells
\begin{equation}
\mathcal{S}_i^0 = \bigcup \limits_{k=1}^{n_e} T^n_{j(k)}, 
\label{stencil}
\end{equation}
where $j=j(k)$ denotes a mapping from the set of integers $k\in[1,n_e]$ to the global indices $j$ of the cells in the mesh given by \eqref{trian}. Stencil $S_i^0$ includes the current cell $T_i^n$ and is filled by recursively adding neighbors of elements that have been already selected, until the desired number $n_e$ is reached. For convenience, we assume that $j(1)=i$ so that the first cell in the 
stencil is always the element for which we are computing the reconstruction. The polynomial $\Popt$ is then defined by imposing that its cell averages on each $T^n_{j}$ match the computed averages 
$\Q^n_j$ in a weak form. Since $n_e>\mathcal{M}$, this of course leads to an overdetermined linear system which is solved using a constrained least-squares technique (LSQ) \cite{Dumbser2007693} as
\begin{equation}
\label{CWENO:Popt}
\begin{aligned}
\Popt &= \underset{{\mathbf{p}\in\mathcal{P}_i}}{\argmin}  
			\sum_{T_j^n \in \mathcal{S}_i^0} 
			\left( \mathbf{Q}^n_{j}-\frac{1}{|T_{j}^n|} \int_{T^n_{j}} \mathbf{p}(\x) d\x \right)^2
       \text{, with}
\\
\mathcal{P}_i &= \left\{\mathbf{p} \in \mathbb{P}_M: \Q^n_{i}=\frac{1}{|T_{i}^n|} \int_{T^n_{i}} \mathbf{p}(\x) d\x\right\}
			\subset \mathbb{P}_M,
\end{aligned}
\end{equation}
where $\mathbb{P}_M$ is the set of all polynomials of degree at most $M$. In other words, the polynomial $\Popt$ has exactly the cell average $\Q^n_{i}$ on the cell $T_i^n$ and matches all the other cell averages in the stencil in the least-square sense. The polynomial $\Popt$ is expressed in terms of the orthogonal Dubiner-type basis functions $\psi_l(\boldsymbol{\xi})$ \cite{Dubiner,orth-basis,CBS-book} on the reference element $T_E$, i.e.
\begin{equation}
\label{eqn.recpolydef} 
\Popt(\x,t^n) = \sum \limits_{l=1}^\mathcal{M} \psi_l(\boldsymbol{\xi}) \hat{\mathbf{p}}^{n}_{l,i} =: \psi_l(\boldsymbol{\xi}) \hat{\mathbf{p}}^{n}_{l,i},   
\end{equation}
where the mapping to the reference coordinate system is given by \eqref{xietaTransf} and $\hat{\mathbf{p}}_{l,i}$ denote the unknown expansion coefficients. In the rest of the paper we will use classical tensor index notation based on the Einstein summation convention, which implies summation over repeated indices. The integrals appearing in \eqref{CWENO:Popt} are computed in the reference system $\boldsymbol{\xi}$ using Gaussian quadrature rules of suitable order, see \cite{stroud}. The transformation from the physical to the reference coordinates is provided by \eqref{xietaTransf} and it prevents the reconstruction matrix of system \eqref{CWENO:Popt} to be ill-conditioned due to scaling effects. 

\begin{figure}[!htbp]
\begin{center}
\begin{tabular}{cc} 
\includegraphics[width=0.44\textwidth]{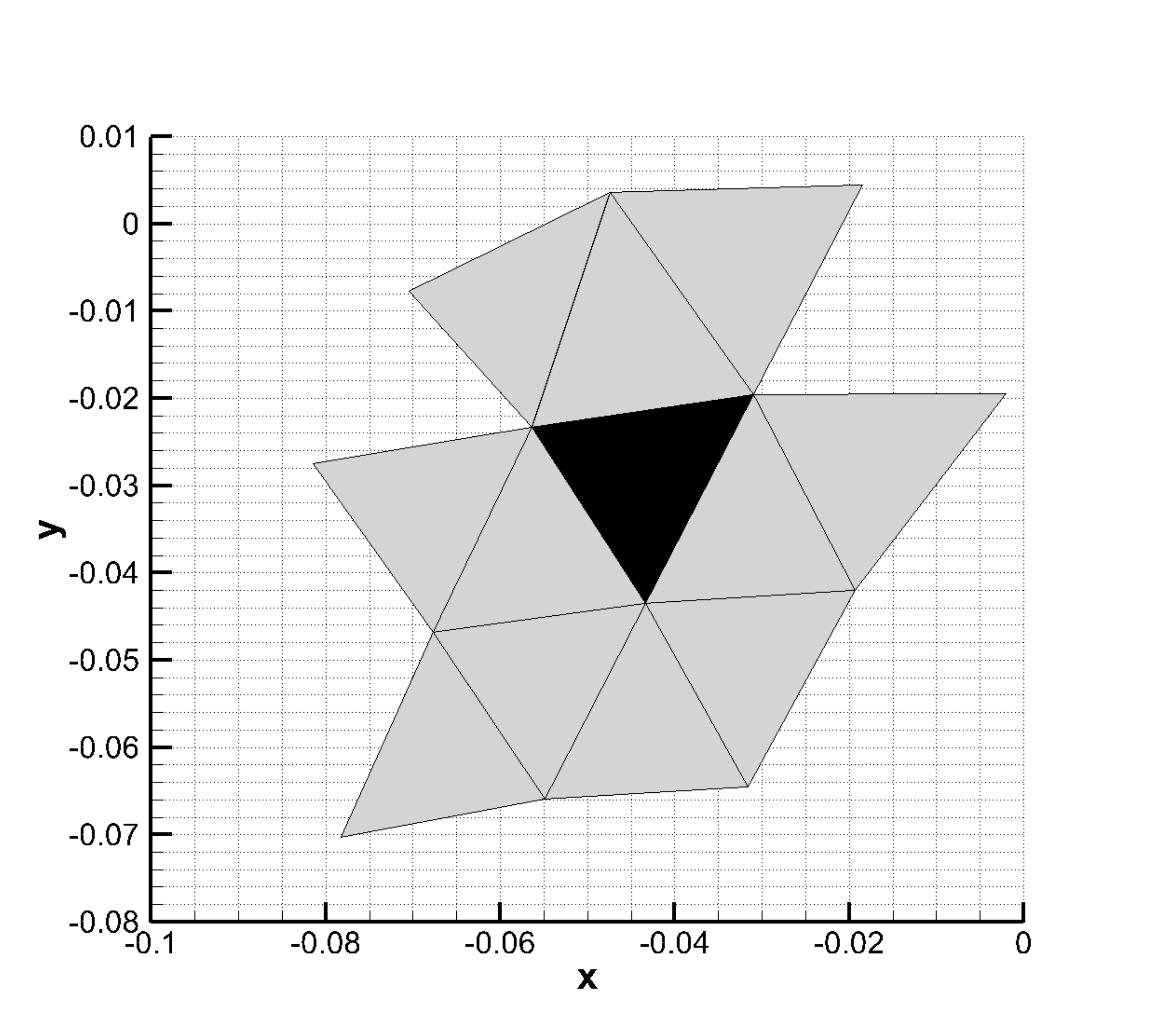} & 
\includegraphics[width=0.44\textwidth]{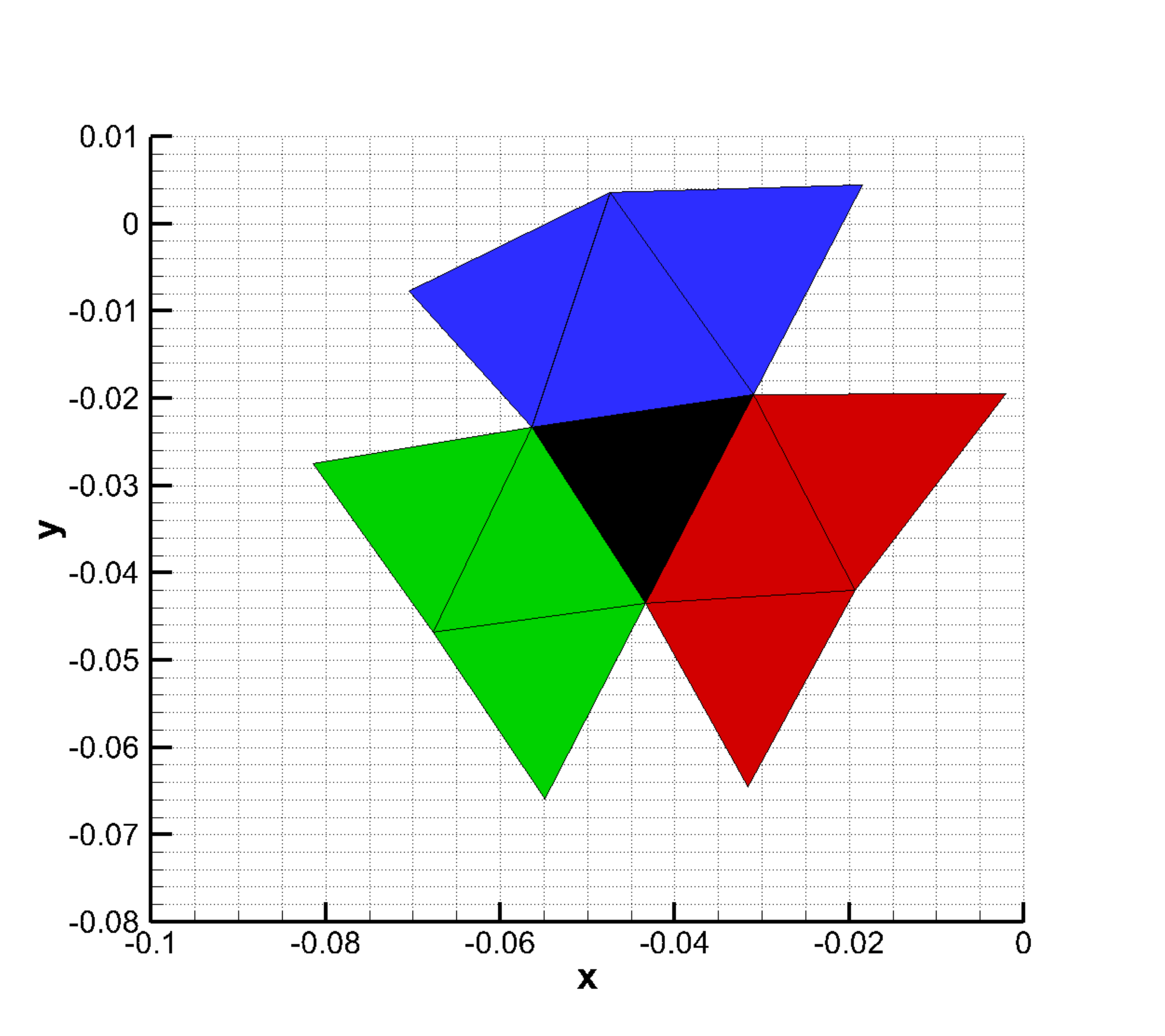}   \\
\includegraphics[width=0.44\textwidth]{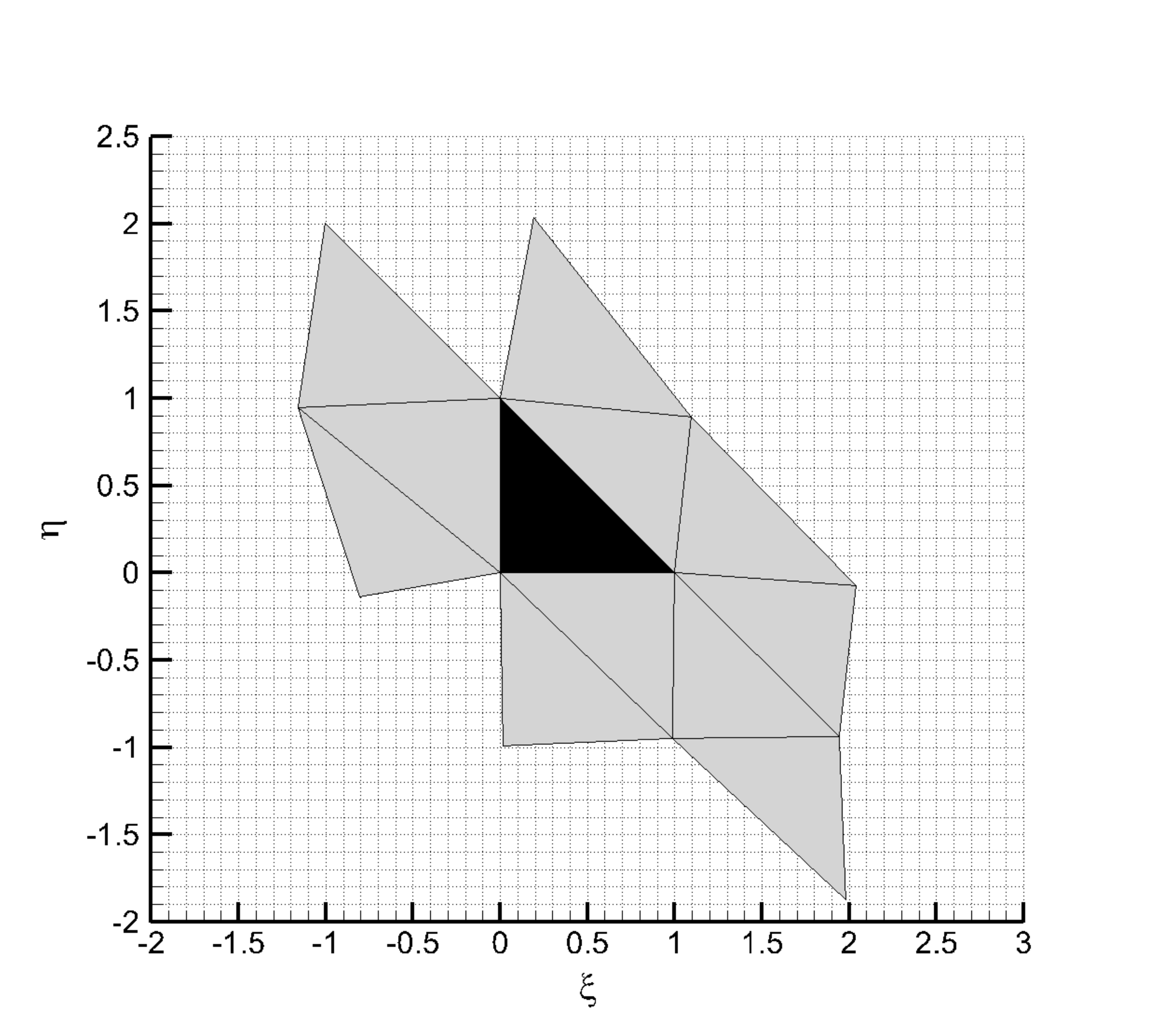} & 
\includegraphics[width=0.44\textwidth]{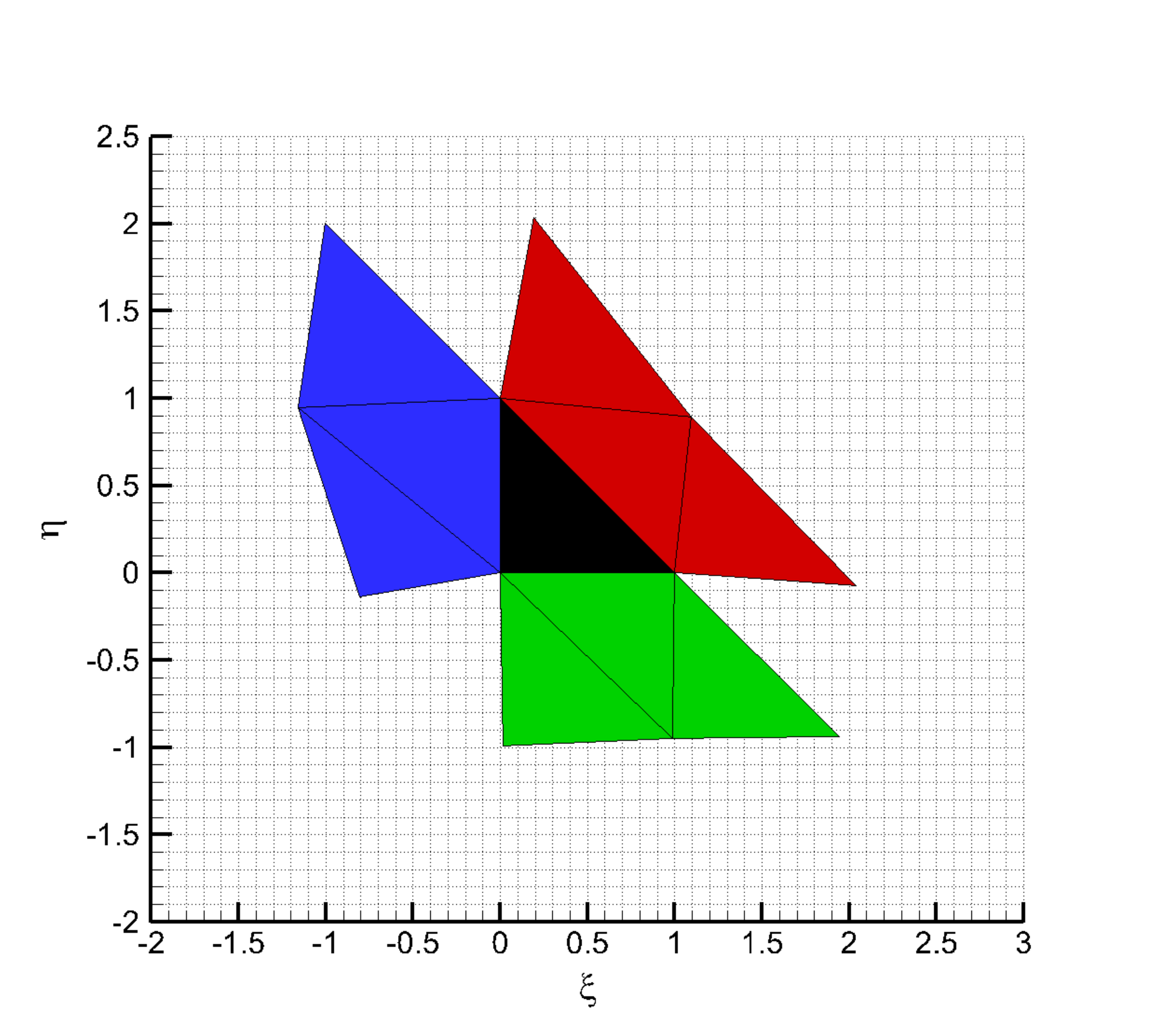}   \\
\end{tabular}
\caption{Two-dimensional reconstruction stencils for $M=2$, hence $n_e=12$, in the physical (top row) and in the reference (bottom row) coordinate system. One central stencil (left) and three sectorial stencils (right).}
\label{fig.stencil2D}
\end{center}
\end{figure}

\begin{figure}[!htbp]
\begin{center}
\begin{tabular}{cc} 
\includegraphics[width=0.44\textwidth]{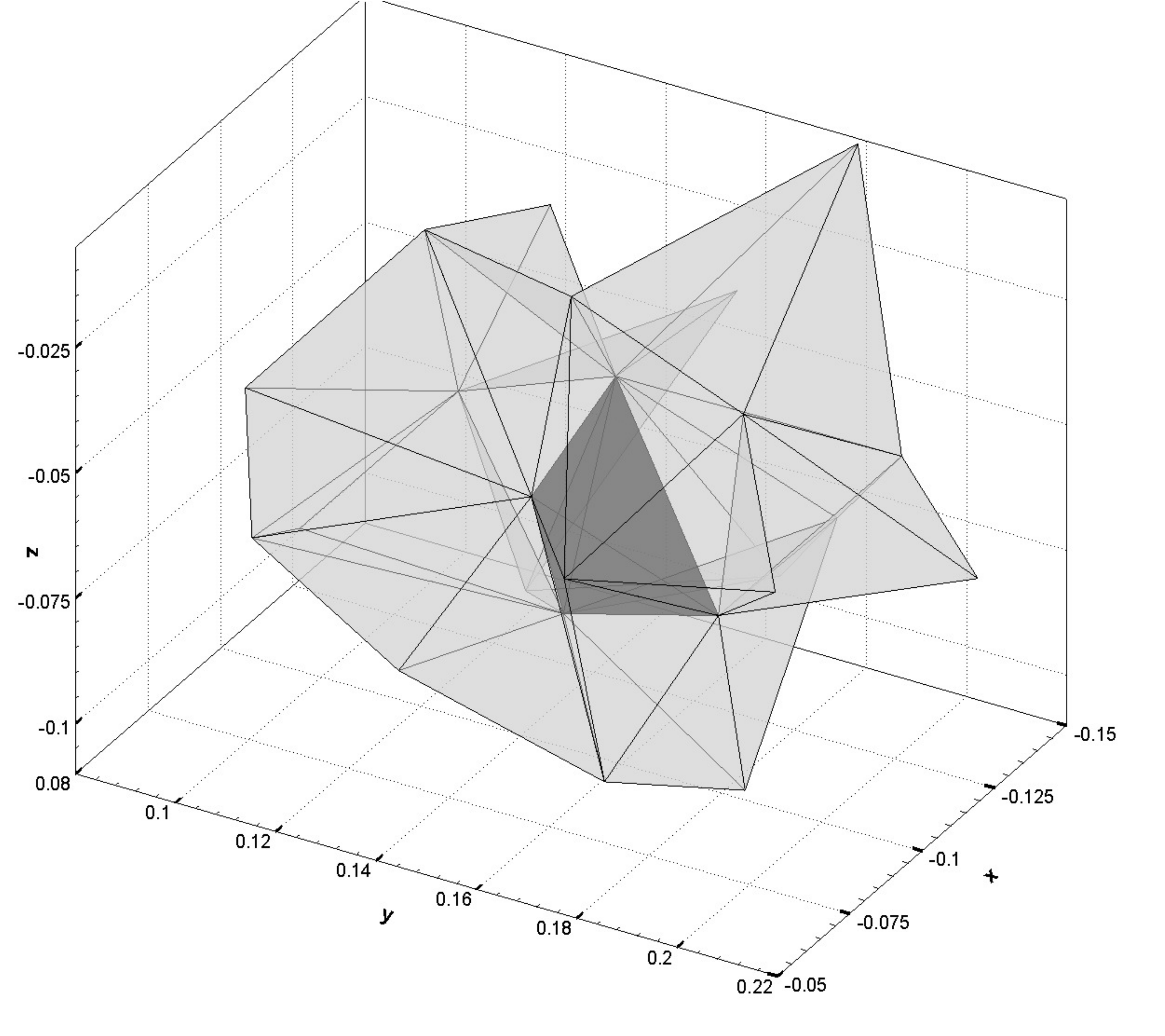} & 
\includegraphics[width=0.44\textwidth]{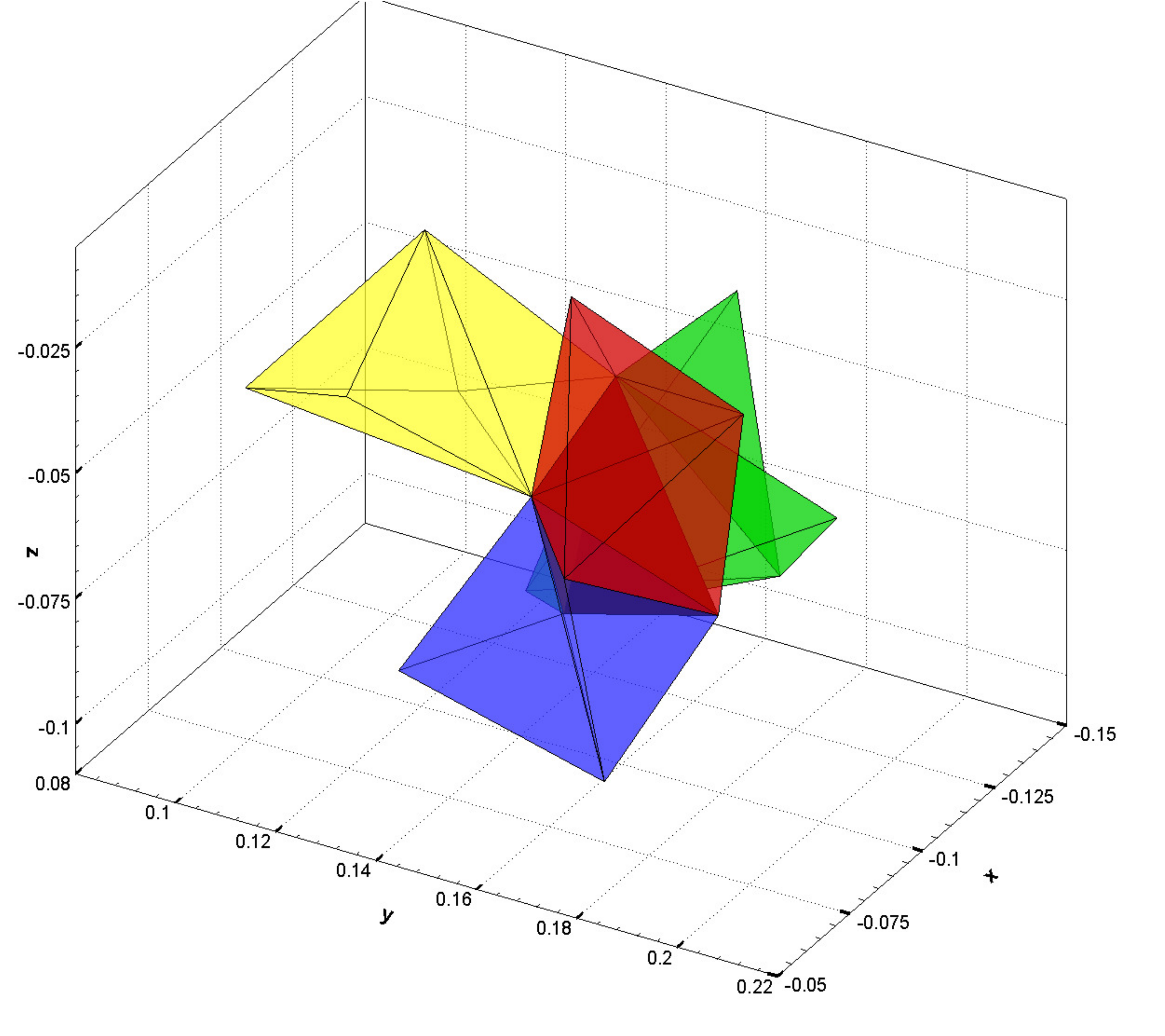}   \\
\includegraphics[width=0.44\textwidth]{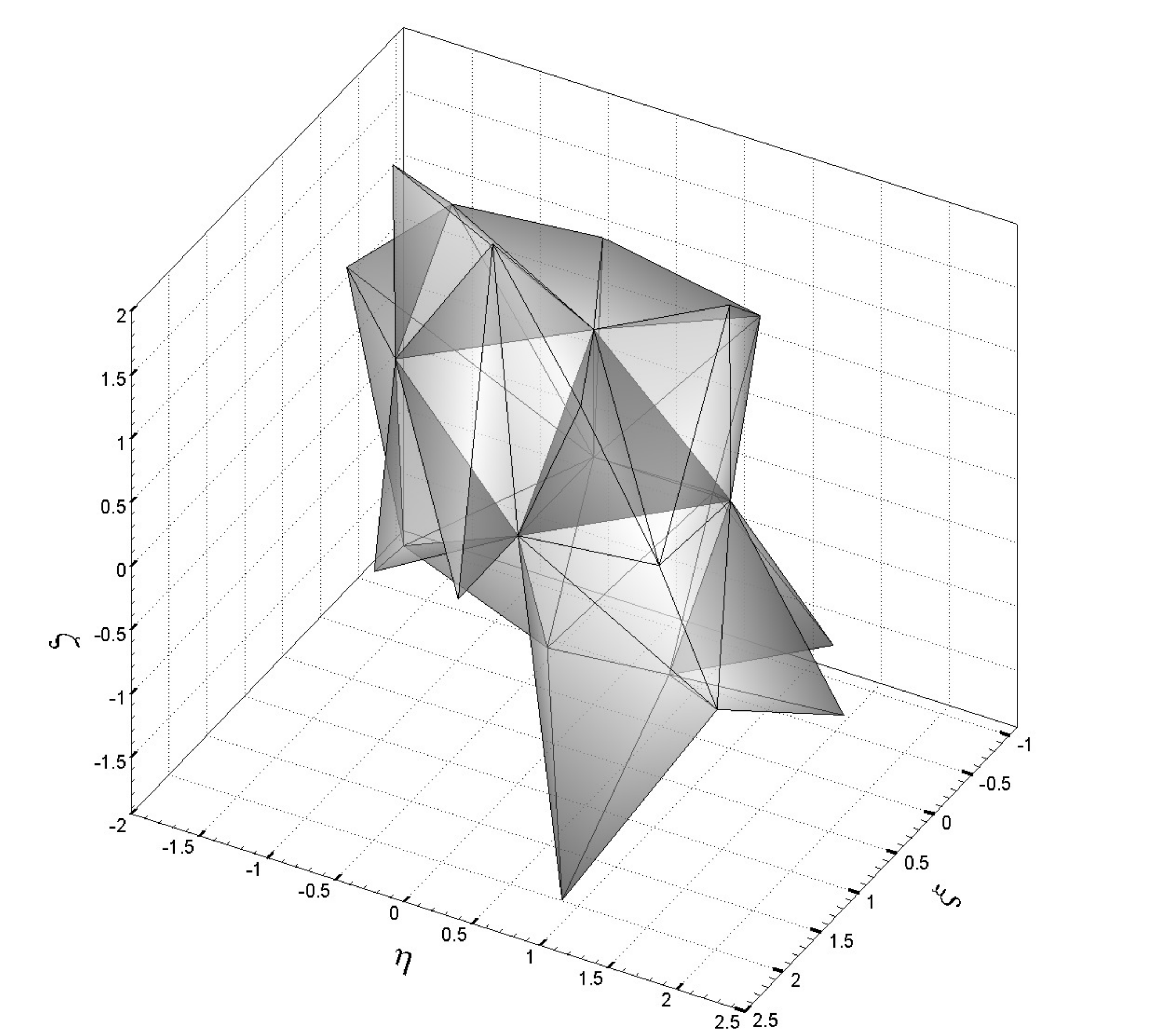} & 
\includegraphics[width=0.44\textwidth]{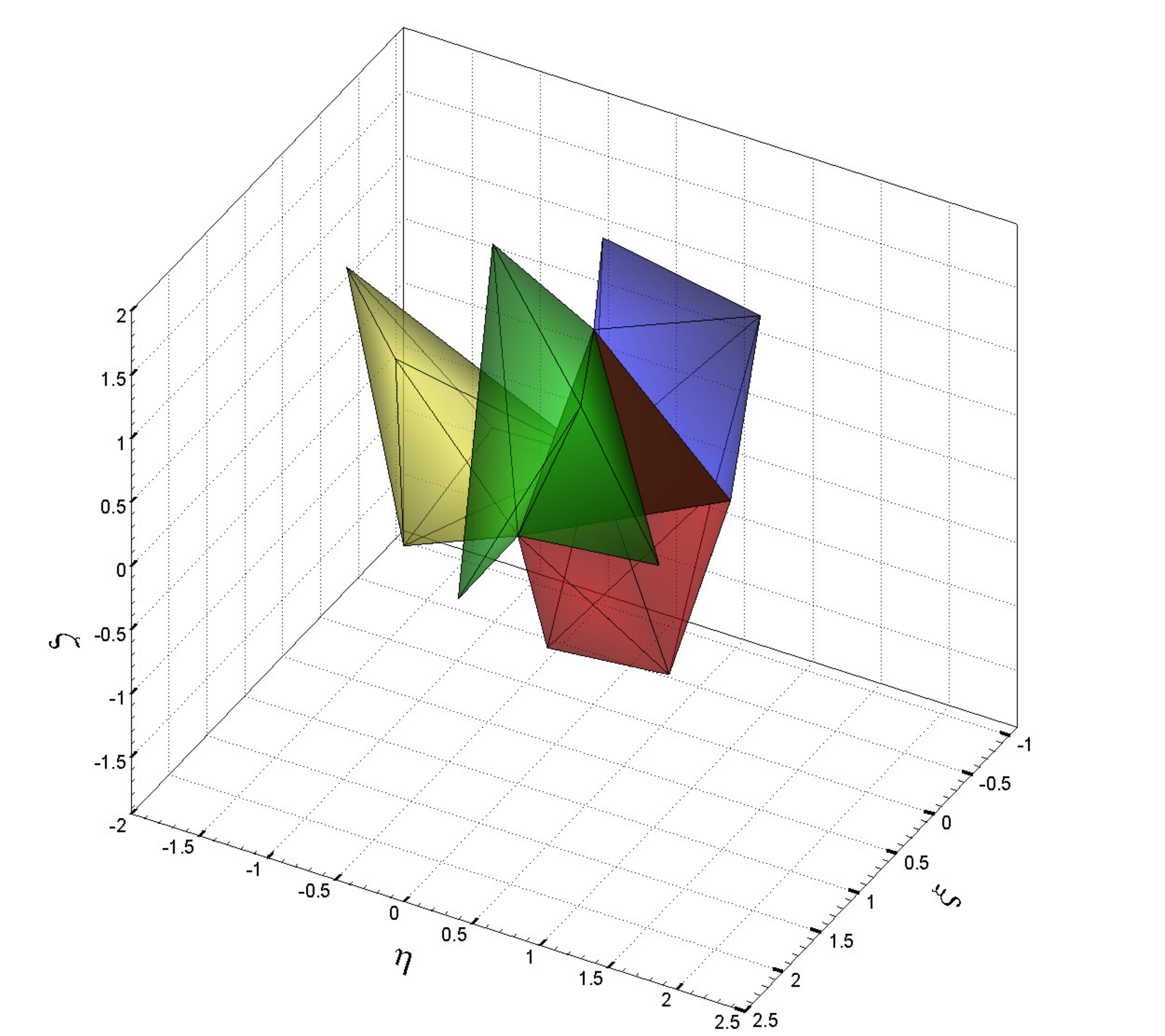}   \\
\end{tabular}
\caption{Three-dimensional reconstruction stencils for $M=2$, hence $n_e=30$, in the physical (top row) and in the reference (bottom row) coordinate system. One central stencil (left) and three sectorial stencils (right).}
\label{fig.stencil3D}
\end{center}
\end{figure}

The CWENO reconstruction makes use also of other polynomials of lower degree and in this paper we choose a total number of $N_p=(d+1)$ interpolating polynomials of degree $M^s=1$ referred to as \textit{sectorial polynomials}. More precisely, we consider $N_p$ stencils $S_i^s$ with $s\in[1,N_p]$, each of them containing exactly $\hat{n}_e=\mathcal{M}(M^s,d) = (d+1)$ cells. $S_i^s$ includes again the central cell $T^n_i$ and is filled by the same recursive algorithm used for $S_i^0$, but selecting only those elements whose barycenter lies in the open cone defined by one vertex and the opposite face of the $T_i^n$, as proposed in \cite{kaeserjcp,Dumbser2007693} and as shown in Figures \ref{fig.stencil2D}-\ref{fig.stencil3D}. 

Let again $j=j^s(k)$ denote the mapping from $k \in [1,\hat{n}_e]$ to the global indices associated with $\mathcal{T}_{\Omega}$. For each stencil $S_i^s$ we compute a linear polynomial by solving the reconstruction systems 
\begin{equation}
\label{CWENO:Ps}
\mathbf{P}_s\in\mathbb{P}_1 
\text{  s.t. } \forall T^n_j \in S_i^s: 
\Q^n_{j}=\frac{1}{|T_{j}^n|} \int_{T^n_j} \mathbf{P}_s(\x) d\x,
\end{equation}
which are not overdetermined and have a unique solution for non-degenerate locations of the cell barycenters. 

Following the general framework introduced in \cite{CWENOframework}, we select a set of positive coefficients $\lambda_0,\ldots,\lambda_{N_p}$ such that $\sum_{l=0}^{N_p}\lambda_l=1$ and we define
\begin{equation}
\label{CWENO:P0}
\mathbf{P}_0 = \frac{1}{\lambda_0}\left(\Popt - \sum_{l=1}^{N_p} \lambda_l \mathbf{P}_l \right) \in\mathbb{P}_M.
\end{equation}
In this way, the linear combination of the polynomials $\mathbf{P}_0,\ldots,\mathbf{P}_{N_p}$ with the coefficients $\lambda_0,\ldots,\lambda_{N_p}$ is equal to $\Popt$. For this reason, these are called {\em optimal coefficients}. Note that $\mathbf{P}_0$ is not directly reconstructed in the CWENO approach, but it is computed by subtracting the weighted sectorial polynomials $\mathbf{P}_l$ with $l \geq 1$ from $\Popt$.  We point out that, in contrast to standard WENO reconstructions, the accuracy of the CWENO reconstruction does not depend on the choice of the optimal coefficients. They can thus be chosen arbitrarily, satisfying only the normalization constraint and the positivity assumption. Furthermore this avoids the appearance of negative weights which could be source of instabilities that should otherwise be cured, for example as described in \cite{ShiHuShu:2002}.

Finally, the sectorial polynomials $\mathbf{P}_s$ with $s\in[1,N_p]$ are nonlinearly combined with the $\mathbf{P}_0$, obtaining $\w_i(\x,t^n)$ as
\begin{equation}
\label{CWENO:Prec}
\w_i(\x,t^n) = \sum_{l=0}^{N_p} \omega_s \mathbf{P}_s(\x),
\end{equation}
where the normalized {\em nonlinear weights} $\boldsymbol{\omega}_s$ are given by
\begin{equation}
\label{eqn.weights}
\boldsymbol{\omega}_s = \frac{\tilde{\boldsymbol{\omega}}_s}{\sum \limits_{m=0}^{N_p} \tilde{\boldsymbol{\omega}}_m}.
\end{equation} 
In the above expression the non-normalized weights $\tilde{\boldsymbol{\omega}}_s$ depend on the linear weights $\lambda_s$ and the oscillation indicators $\boldsymbol{\sigma}_s$, hence
\begin{equation}
\tilde{\boldsymbol{\omega}}_s = \frac{\lambda_s}{\left(\boldsymbol{\sigma}_s + \epsilon \right)^r},
\label{eqn.weights2}
\end{equation}
with the parameters $\epsilon=10^{-14}$ and $r=4$ chosen according to \cite{Dumbser2007693}.
Note that in smooth areas, $\omega_s\simeq\lambda_s$ and then $w_i\simeq\Popt$, so that we recover optimal accuracy. On the other hand, close to a discontinuity, $\mathbf{P}_0$ and some of the low degree polynomials $\mathbf{P}_s$ would be oscillatory and have high oscillation indicators, leading to $\omega_s\simeq0$ and in these cases only lower order non-oscillatory data are employed in $w_i$, guaranteeing the non-oscillatory property of the reconstruction.

 As linear weights we take $\lambda_0=10^5$ for $\mathcal{S}_i^0$ and $\lambda_s=1$ for the sectorial stencils, as suggested in \cite{Dumbser2007693}. The oscillation indicators $\boldsymbol{\sigma}_s$ appearing in \eqref{eqn.weights} are given by 
\begin{equation}
\boldsymbol{\sigma}_s = \Sigma_{lm} \hat{\mathbf{p}}^{n,s}_{l,i} \hat{\mathbf{p}}^{n,s}_{m,i},
\end{equation}
according to \cite{Dumbser2007693}, using the universal oscillation indicator matrix $\Sigma_{lm}$ that reads
\begin{equation}
\Sigma_{lm} = \sum \limits_{ 1 \leq \alpha + \beta + \gamma \leq M}  \, \, \int \limits_{T_E} \frac{\partial^{\alpha+\beta+\gamma} \psi_l(\xi,\eta,\zeta)}{\partial \xi^\alpha \partial \eta^\beta \partial \zeta^\gamma} \cdot 
                                                                         \frac{\partial^{\alpha+\beta+\gamma} \psi_m(\xi,\eta,\zeta)}{\partial \xi^\alpha \partial \eta^\beta \partial \zeta^\gamma} d\xi d\eta d\zeta.   
\end{equation}  
Since the entire reconstruction procedure is carried out on the reference system $\boldsymbol{\xi}$, matrix $\Sigma_{lm}$ only depends on the reconstruction basis functions $\psi_m(\boldsymbol{\xi})$ and \textit{not} on the mesh, therefore it can be conveniently precomputed once and stored.

Note that also the picking of the stencil elements can be performed and stored at the beginning of the computation because the mesh topology remains fixed. In the Eulerian case, i.e. when the mesh velocity is zero, the matrices involved in the local least-squares problem \eqref{CWENO:Popt} can be precomputed, inverted and saved in the pre-processing stage. Unfortunately in the ALE setting the element matrices change with the motion of the mesh and thus the local problem must be assembled and solved at each time step, which is of course computationally more demanding but also has a smaller memory footprint.

Finally we point out that the use of the new central WENO reconstruction instead of the classical unstructured WENO formulation \cite{Dumbser2007693,Lagrange3D} in the ADER context improves the overall algorithm efficiency, since the sectorial stencils are filled with a smaller number of elements because they are only of degree $M^s=1$. Furthermore, we only need a total number of $N_p=(d+2)$ 
stencils, which is not the case in the aforementioned WENO formulations, where it was necessary to consider $N_p=7$ stencils in 2D or even $N_p=9$ stencils in 3D. This fact has been properly demonstrated by monitoring the computational time for the simulations reported in Section \ref{sec.validation} and the results are given in Table \ref{tab.efficiency}.

\subsection{Local space-time Galerkin predictor} 
\label{sec.localCG}
In order to achieve high order of accuracy in time we rely on the ADER approach, which makes use of the high order spatial reconstruction polynomials $\w_i(\x,t^n)$ to construct an element-local solution 
\textit{in the small} \cite{eno} of the Cauchy problem with given initial data $\mathbf{w}_i(\mathbf{x},t^n)$. This solution in the small is sought inside each space-time control volume 
$T_i(t) \times [t^n; t^{n+1}]$ under the form of a space-time polynomial of degree $M$, which is employed later to evaluate the numerical fluxes at element interfaces. The original formulation was proposed by Toro et al. \cite{toro3,toro10,CastroToro} and was based on the so-called Cauchy-Kovalewski procedure in which time derivatives are replaced by space derivatives using repeatedly the governing conservation laws \eqref{eqn.PDE} in their differential form. This formulation is based on a Taylor expansion in time, hence problems arise in case of PDE with stiff source terms. Furthermore, the Cauchy-Kovalewski procedure becomes very cumbersome or even unfeasible for general complex nonlinear systems of conservation laws. For that reason, this is why we employ here a more recent version of the ADER  method that is based on a local space-time finite element formulation introduced in \cite{DumbserEnauxToro,Dumbser2008}. The reconstruction polynomials $\w_i(\x,t^n)$ computed 
at the current time $t^n$ are locally evolved in time by applying a weak form of the governing PDE in space and time, but considering only the element itself. No neighbor information are required during  this \textit{predictor} stage because the coupling with the neighbor elements occurs only later in the final one-step finite volume scheme (see Section \ref{sec.SolAlg}). Therefore, the proposed
approach can also be seen as a \textit{predictor-corrector} method. This procedure has already been successfully applied in the context of moving mesh finite volume schemes, see e.g. 
\cite{Lagrange3D,LagrangeMHD,ALEMOOD1}. 

First the governing PDE \eqref{eqn.PDE} is written in the space-time reference system defined by the coordinate vector $\boldsymbol{\tilde{\xi}}=(\xi,\eta,\zeta,\tau)$ as
\begin{equation}
\frac{\partial \Q}{\partial \tau} \frac{\partial \tau}{\partial t} + \frac{\partial \Q}{\partial \xxi} \frac{\partial \xxi}{\partial t} + \left(\frac{\partial \xxi}{\partial \x}\right)^T \nabla_{\xxi} \cdot \F(\Q) = \mathbf{0}
\label{eqn.PDExi}
\end{equation}
with
\begin{equation}
\nabla_{\xxi} = \left( \begin{array}{c} \frac{\partial}{\partial \xi} \\ \frac{\partial}{\partial \eta} \\ \frac{\partial}{\partial \zeta}  \end{array} \right), \qquad
\left(\frac{\partial \xxi}{\partial \x}\right) = \left( \begin{array}{ccc} \xi_x & \eta _x & \zeta _x \\ \xi_y & \eta_y & \zeta _y \\ \xi_z & \eta_z & \zeta _z \end{array} \right).
\label{eqn.nabla}
\end{equation}
The linear transformation \eqref{xietaTransf} provides the spatial coordinate vector $\xxi=(\xi,\eta,\zeta)$, while for the reference time $\tau$ we adopt the simple mapping 
\begin{equation}
\tau = \frac{(t-t^n) }{\Delta t} \quad \textnormal{ with } \quad \Delta t = \textnormal{CFL} \, \min \limits_{T_i^n} \frac{h_i}{|\lambda_{\max,i}|} \qquad \forall T_i^n \in \Omega^n,
\label{eqn.timestep}
\end{equation}
where the time step $\Delta t$ is computed under a classical (global) Courant-Friedrichs-Levy number (CFL) stability condition with $\textnormal{CFL} \leq \frac{1}{d}$. Furthermore $h_i$ represents a characteristic element size, either the incircle or the insphere diameter for triangles or tetrahedra, respectively, while $|\lambda_{\max,i}|$ is taken to be the maximum absolute value of the eigenvalues computed from the current solution $\Q_i^n$ in $T_i^n$. According to \eqref{eqn.timestep} one obtains $\frac{\partial \tau}{\partial t}=\frac{1}{\Delta t}$, therefore Eqn. \eqref{eqn.PDExi} can be written more compactly as
\begin{equation}
\frac{\partial \Q}{\partial \tau} + \Delta t \, \mathbf{H}(\Q) = \mathbf{0}
\label{eqn.PDExi2}
\end{equation}
by means of the abbreviation
\begin{equation}
\mathbf{H}(\Q) :=  \frac{\partial \Q}{\partial \xxi} \frac{\partial \xxi}{\partial t} + \left(\frac{\partial \xxi}{\partial \x}\right)^T \nabla_{\xxi} \cdot \F(\Q).
\label{eqn.Hterm}
\end{equation}
Note that the ALE framework yields a \textit{moving} control volume which is correctly taken into account by the term $\frac{\partial \xxi}{\partial t}$ in \eqref{eqn.PDExi} that is not present for the Eulerian case introduced in \cite{Dumbser2008}.

\begin{figure}[!htbp]
\centering
\begin{tabular}{cc} 
\includegraphics[width=0.3\textwidth]{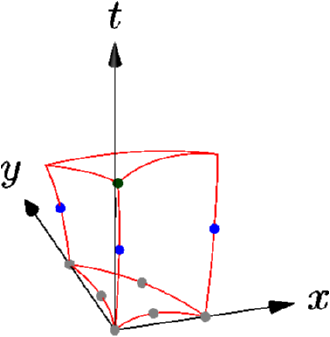}    & 
\includegraphics[width=0.3\textwidth]{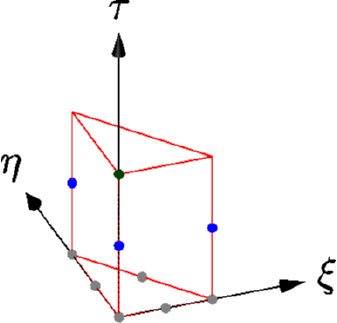}     
\end{tabular}
\caption{Isoparametric mapping of the physical space-time element $T_i(t) \times \Delta t$ (left) to the space-time reference element $T_E \times [0,1]$ (right) used within the local space-time Galerkin predictor for a triangular control volume. Space-time nodes, i.e. the degrees of freedom of the transformation, are highlighted for $M=2$.}
\label{fig.stcg.map}
\end{figure}

To obtain a high order predictor solution, the conservation law \eqref{eqn.PDExi2} is multiplied by a space-time test function $\theta_k=\theta_k(\boldsymbol{\tilde{\xi}})$ and then integrated over the space-time reference element $T_E \times [0,1]$ leading to
\begin{equation}
\int \limits_{0}^{1} \int \limits_{T_E} \theta_k \cdot \frac{\partial \Q}{\partial \tau} \, d\xxi d\tau= -\int \limits_{0}^{1} \int \limits_{T_E} \Delta t \, \theta_k \cdot \mathbf{H}(\Q) \, d\xxi d\tau
\label{eqn.PDEweak}
\end{equation}
which is a \textit{weak formulation} of the governing equations \eqref{eqn.PDE}. The predictor solution is given by a high order space-time polynomial $\q_h$ that is approximated by means of a set of space-time basis functions $\theta_{l}(\boldsymbol{\tilde{\xi}})$ as
\begin{equation}
\q_h = \q_h(\boldsymbol{\tilde{\xi}}) = \sum \limits_{l=1}^{\mathcal{L}} \theta_l(\boldsymbol{\tilde{\xi}}) \widehat{\q}_{l}:= \theta_l(\boldsymbol{\tilde{\xi}}) \widehat{\q}_{l},
\label{eqn.qh}
\end{equation}
where $\mathcal{L}=\mathcal{M}(M,d+1)$ represents the total number of degrees of freedom $\widehat{\q}_{l}$ needed to reach the formal order of accuracy $M+1$ in space \textit{and} time, 
as given by \eqref{eqn.cwenoM} in $d+1$ dimensions. The basis functions $\theta_l$ are assumed to be the same as the test functions $\theta_k$ and they are defined by the Lagrange interpolation polynomials passing through a set of space-time nodes $\boldsymbol{\tilde{\xi}}_m$ with $m\in[1,\mathcal{L}]$ explicitly specified in \cite{Dumbser2008}, see Figure \ref{fig.stcg.map} for $M=2$.

Since the functions $\theta_l$ yield a \textit{nodal basis}, the numerical approximation for the fluxes $\F(\Q)$ as well as for the abbreviation $\mathbf{H}(\Q)$ simply reads
\begin{equation}
\begin{array}{l}
\F_h = \theta_l(\boldsymbol{\tilde{\xi}}) \, \F(\widehat{\q}_{l}) := \theta_l(\boldsymbol{\tilde{\xi}}) \, \widehat{\F}_{l}, \\
\mathbf{H}_h = \theta_l(\boldsymbol{\tilde{\xi}}) \, \mathbf{H}(\widehat{\q}_{l}) := \theta_l(\boldsymbol{\tilde{\xi}}) \, \widehat{\mathbf{H}}_{l}.
\end{array}
\label{eqn.PDEterm}
\end{equation}
As a consequence the weak formulation \eqref{eqn.PDEweak} becomes
\begin{equation}
\int \limits_{0}^{1} \int \limits_{T_E} \theta_k \cdot \frac{\partial \theta_l}{\partial \tau} \widehat{\q}_{l} \, d\xxi d\tau = -\Delta t \int \limits_{0}^{1} \int \limits_{T_E} \theta_k \theta_l \widehat{\mathbf{H}}_{l} \, d\xxi d\tau
\label{eqn.PDEweak2}
\end{equation}
that simplifies to
\begin{equation}
\mathbf{K_\tau} \widehat{\q}_{l} = -\Delta t \, \mathbf{M} \widehat{\mathbf{H}}_{l} \quad \textnormal{ with } \quad
\mathbf{K_\tau}:= \int \limits_{0}^{1} \int \limits_{T_E} \theta_k \cdot \frac{\partial \theta_l}{\partial \tau} \, d\xxi d\tau, \quad
\mathbf{M} := \int \limits_{0}^{1} \int \limits_{T_E} \theta_k \theta_l \, d\xxi d\tau.
\label{eqn.cgPDE}
\end{equation}
The above expression constitutes a \textit{nonlinear} algebraic system of equations that is solved by an iterative procedure. According to \cite{Dumbser2008} the vector of the degrees of freedom $\widehat{\q}_{l}$ is split into two parts, namely the \textit{known} degrees of freedom $\widehat{\q}_{l}^0$ that come from the reconstruction polynomial $\w_i$ at the current reference time $\tau=0$ and the \textit{unknown} degrees of freedom $\widehat{\q}_{l}^1$ given for $\tau>0$. Therefore $\widehat{\q}_{l}=(\widehat{\q}_{l}^0,\widehat{\q}_{l}^1)$ and subsequently matrix $\mathbf{K_\tau}$ can be also divided into $\mathbf{K_\tau}=(\mathbf{K_\tau}^0,\mathbf{K_\tau}^1)$, so that the iterative procedure for the solution of system \eqref{eqn.cgPDE} reads
\begin{equation}
\mathbf{K_\tau}^1 \widehat{\q}_{l}^{1, r+1} = -\Delta t \, \mathbf{M} \widehat{\mathbf{H}}_{l}^r - \mathbf{K_\tau}^0 \widehat{\q}_{l}^{0, r},  
\label{eqn.CGfinal}
\end{equation}
where $r$ denotes the iteration number.

If the mesh moves, the space-time volume changes its configuration in time, hence implying the following ODE to be considered:
\begin{equation}
\frac{d\x}{dt} = \mathbf{V}(\x,t),
\label{eqn.ODEmesh}
\end{equation}
which is typically addressed as \textit{trajectory equation}. The element geometry defined by $\x$ as well as the local mesh velocity $\mathbf{V}(\x,t)=(U,V,W)$ are approximated using again the same basis functions $\theta_l$ as
\begin{equation}
\x=\x_h=\theta_l(\boldsymbol{\tilde{\xi}}) \, \widehat{\x}_{l}, \qquad \mathbf{V}=\mathbf{V}_h=\theta_l(\boldsymbol{\tilde{\xi}}) \, \widehat{\mathbf{V}}_{l},
\label{eqn.xh_vh}
\end{equation}
leading to an isoparametric approach. The ODE \eqref{eqn.ODEmesh} is solved employing the same strategy adopted for system \eqref{eqn.cgPDE}, that is
\begin{equation}
\mathbf{K_\tau}^1 \widehat{\x}_{l}^{1, r+1} = \Delta t \, \mathbf{M} \widehat{\mathbf{V}}_{l}^r - \mathbf{K_\tau}^0 \widehat{\x}_{l}^{0, r}
\label{eqn.ODEfinal}
\end{equation}
and the iteration procedure is carried out \textit{together} with the solution of the nonlinear system \eqref{eqn.CGfinal} in a coupled manner until the residuals of both equations are less than a prescribed tolerance. 

The local Galerkin procedure described above has to be performed for all elements $T_i(t)$ of the computational domain, hence producing the space-time predictor for the solution $\q_h$, for the fluxes 
$\F_h=(\f_h,\g_h,\h_h)$ and also for the mesh velocity $\mathbf{V}_h$. We stress again that this predictor is computed \textit{without} considering any neighbor information, but only solving the 
evolution equations \eqref{eqn.PDE} locally in the space-time control volume $T_i(t) \times [t^n,t^{n+1}]$.

\subsection{Mesh motion}
\label{sec.meshMot}
In the ALE approach the mesh velocity $\mathbf{V}$ can be chosen independently from the local fluid velocity $\mathbf{v}$, therefore containing both, Eulerian and Lagrangian-like algorithms as special cases
for $\mathbf{V}=\mathbf{0}$ and $\mathbf{V}=\mathbf{v}$, respectively. If the mesh velocity is not assumed to be zero, then it must be properly evaluated and in this section we are going to describe how we derive it in our direct ALE scheme.

The mesh motion is already considered within the Galerkin predictor procedure by the trajectory equation \eqref{eqn.ODEmesh} and the high order space-time representation of the local mesh velocity $\mathbf{V}_h$ is available. Since the predictor step is carried out for each element $T_i(t)$ without involving any coupling with the neighbor elements, at the next time level $t^{n+1}$ the mesh might be discontinuous. In other words, each vertex $k$ of the mesh may be assigned with a different velocity vector $\mathbf{V}_{k,j}$ that is computed from the predictor solution $\mathbf{V}_h$ of the neighbor element $T_j$, hence yielding either holes or overlapping between control volumes in the new mesh configuration $\mathcal{T}^{n+1}_{\Omega}$. The so-called \textit{nodal solver} algorithms aim at fixing a \textit{unique} velocity vector $\mathbf{V}_k$ for each node $k$ of the grid. A lot of effort has been put in the past into the design and the implementation of nodal solvers \cite{chengshu1,chengshu3,chengshu4,Maire2009,Maire2010,Maire2011,Despres2009} and for a comparison among different nodal solvers applied to ADER-WENO finite volume schemes the reader is referred to \cite{LagrangeMHD}. 

Here, we rely on a simple and very robust technique that evaluates the node velocity as a mass weighted sum among the local velocity contributions $\mathbf{V}_{k,j}$ of the Voronoi neighborhood $\mathcal{V}_k$ of vertex $k$, which is composed by all those elements $T_j$ that share the common node $k$. The computation is simply given by
\begin{equation}
\mathbf{V}_k = \frac{1}{\mu_k}\sum \limits_{T_j \in \mathcal{V}_k}{\mu_{k,j}\mathbf{V}_{k,j}}, 
\label{eqnNScs}
\end{equation}
with
\begin{equation}
\mu_k = \sum \limits_{T_j \in \mathcal{V}_k}{\mu_{k,j}}, \qquad \mu_{k,j}=\rho^n_j |T_j^n|.
\label{eqn.NScs.weights}
\end{equation}
In the above expression $\mu_{k,j}$ are the local weights, i.e. the masses of the Voronoi neighbors $T_j$ which are defined multiplying the cell averaged value of density $\rho^n_j$ with the cell volume $|T_j^n|$ at the current time $t^n$. The velocity contributions $\mathbf{V}_{k,j}$ are taken to be the time integral of the high order extrapolated velocity at node $k$ from element $T_j$, i.e.
\begin{equation}
\mathbf{V}_{k,j} = \left( \int \limits_{0}^{1} \theta_l(\boldsymbol{\tilde{\xi}}_{m(k)}) d \tau \right) \widehat{\mathbf{V}}_{l,j}, 
\label{NodesVel}
\end{equation} 
where $m(k)$ is a mapping from the global node number $k$ to the local node number in element $T_j^n$ defined by the position vector $\boldsymbol{\tilde{\xi}}_{m(k)}$ in the reference system.

The new vertex coordinates $\mathbf{X}_k^{n+1}$ are then obtained starting from the old ones $\mathbf{X}_k^{n}$ as
\begin{equation}
\mathbf{X}_k^{n+1} = \mathbf{X}_k^{n} + \Delta t \cdot \mathbf{V}_k,
\label{eqn.newX}
\end{equation}
that is the solution of the trajectory equation \eqref{eqn.ODEmesh}. In the Eulerian case the mesh velocity is set to zero, i.e. $\mathbf{V}=\mathbf{0}$, thus obtaining $\mathbf{X}_k^{n+1} = \mathbf{X}_k^{n}$. 

If the node velocities $\mathbf{V}_k$ lead to bad quality control volumes, e.g. highly stretched or compressed or even tangled elements that might arise from the local fluid velocity, 
in particular in the presence of strong shear flows, an \textit{optimization} is typically performed to improve the geometrical quality of the cells. This is done through appropriate rezoning algorithms \cite{KnuppRezoning,MaireRezoning,ShashkovFSM} by minimizing a local goal functional that does not depend on physical quantities but is only based on the node positions. As a result such strategies provide for each vertex a new coordinate that ensures a better mesh quality, hence avoiding the occurrence of tangled elements or very small time steps according to \eqref{eqn.timestep}. A detailed description of the rezoning algorithm employed in our approach can be found in \cite{Lagrange3D}.

\subsection{Finite volume scheme}
\label{sec.SolAlg}
Our numerical method belongs to the category of direct ALE finite volume schemes, therefore we first have to define the space-time control volume $\mathcal{C}^n_i = T_i(t) \times \left[t^{n}; t^{n+1}\right]$ used to perform the time evolution of the governing PDE \eqref{eqn.PDE}. For each element $T_i(t)$ the old vertex coordinates $\mathbf{X}_{m(k)}^{n}$ are connected to the new vertex positions $\mathbf{X}_{m(k)}^{n+1}$ with \textit{straight} line segments and the resulting space-time volume is bounded by five sub-surfaces or six sub-volumes in 2D or in 3D, respectively: specifically, the bottom and the top of $\mathcal{C}^n_i$ are given by the element configurations $T_i^{n}$ at the current time level and $T_i^{n+1}$ at the new time level, respectively, while the volume is closed by a total number of $\mathcal{N}_i=(d+1)$ lateral sub-volumes $\partial C^n_{ij} = \partial T_{ij}(t) \times [t^n;t^{n+1}]$ that are generated by the evolution of each face $\partial T_{ij}(t)$ shared between element $T_i$ and its direct neighbor $T_j$, hence
\begin{equation}
\partial C^n_i = \left( \bigcup \limits_{T_j \in \mathcal{N}_i} \partial C^n_{ij} \right) 
\,\, \cup \,\, T_i^{n} \,\, \cup \,\, T_i^{n+1}.  
\label{eqn.dCi}
\end{equation}     

\begin{figure}[!htbp]
\centering
\begin{tabular}{c} 
\includegraphics[width=0.85\textwidth]{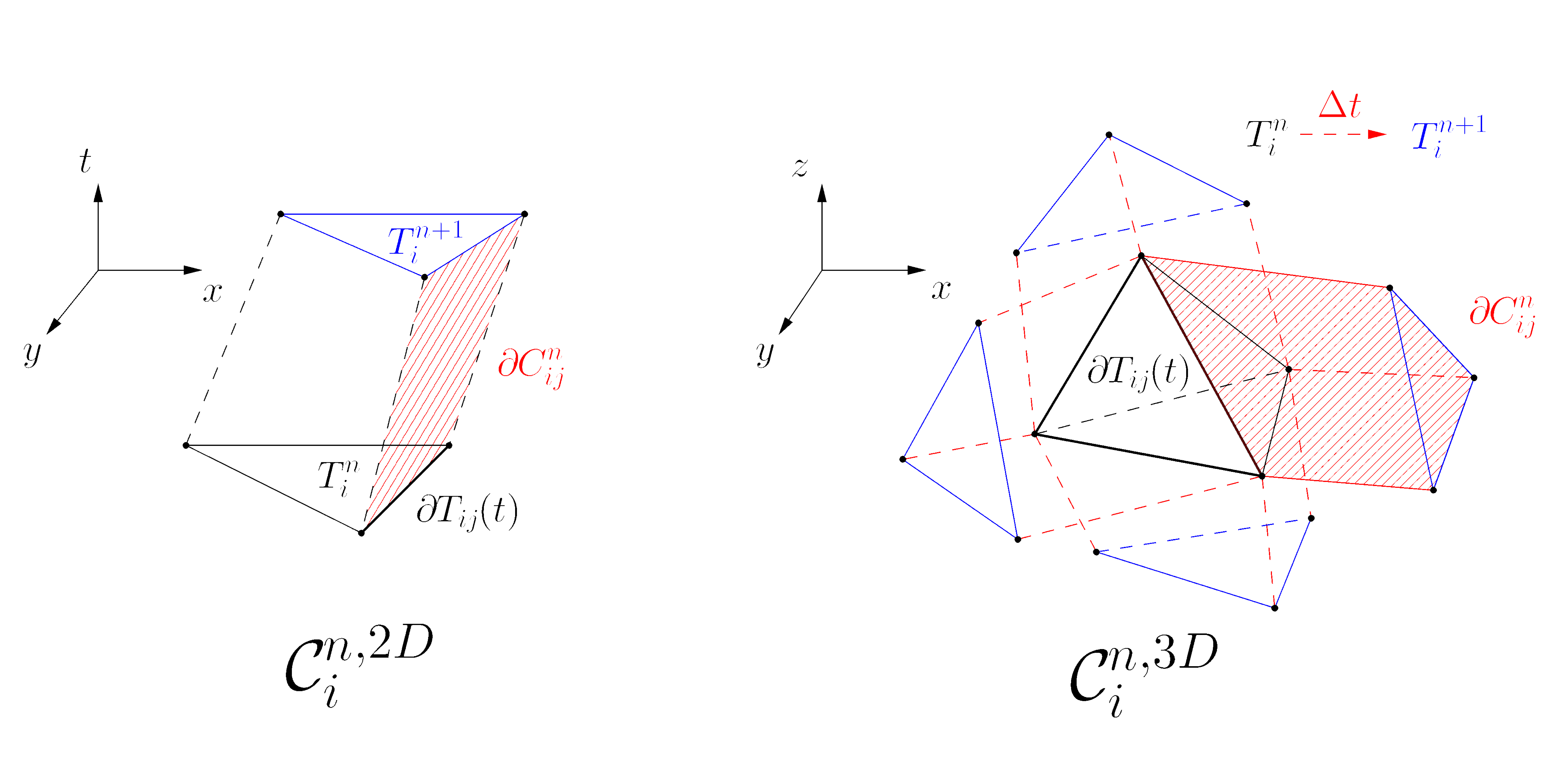}
\end{tabular}
\caption{Space-time control volume $\mathcal{C}_i^n$ given by the evolution of element $T_i$ within one timestep $\Delta t$ in two (left) and three (right) space dimensions. The dashed red lines denote the evolution in time of the faces $\partial C^n_{ij}$ of the control volume $T_i$, whose configuration at the current time level $t^n$ and at the new time level $t^{n+1}$ is depicted in black and blue, respectively.}
\label{fig.STelem}
\end{figure}

To simplify the integral computation each sub-volume is mapped to a reference element: for the bottom and the top of $\mathcal{C}^n_i$ we use the transformation \eqref{xietaTransf} with $(\xi,\eta,\zeta)\in[0,1]$, while for the lateral sub-volumes $\partial T_{ij}(t)$ we define a local reference system $\boldsymbol{\chi}=(\chi_1,\chi_2,\tau)$ that is parametrizes the face $\partial T_{ij}(t)$ of element 
$T_i(t)$. The spatial coordinates $(\chi_1,\chi_2)$ lie on the face $\partial T_{ij}(t)$ and are orthogonal to the time coordinate $\tau$ and the space-time surface is parametrized using a set of bilinear basis  functions $\beta_l(\chi_1,\chi_2,\tau)$ which read 
\begin{eqnarray}
\beta_1(\chi_1,\chi_2,\tau) = (1-\chi_1-\chi_2) \tau , \quad & \beta_4(\chi_1,\chi_2,\tau) = (1-\chi_1-\chi_2)(1-\tau)  \nonumber \\
\beta_2(\chi_1,\chi_2,\tau) = \chi_1\tau,              \quad & \beta_5(\chi_1,\chi_2,\tau) = \chi_1(1-\tau)    \nonumber \\
\beta_3(\chi_1,\chi_2,\tau) = \chi_2\tau,              \quad & \beta_6(\chi_1,\chi_2,\tau) = \chi_2(1-\tau)    
\label{eqn.BetaBaseFunc}
\end{eqnarray}
so that the analytical description of the space-time sub-volume $\partial C^n_{ij}$ results in 
\begin{equation}
 \mathbf{\tilde{x}} \left( \chi_1,\chi_2,\tau \right) = 
 \sum\limits_{l=1}^{N_l} \beta_l(\chi_1,\chi_2,\tau) \, \mathbf{\tilde{X}}_{ij,l} ,	
\label{eqn.SurfPar}
\end{equation}
with $0 \leq \chi_1 \leq 1$, $0 \leq \chi_2 \leq 1 - \chi_1$, $0 \leq \tau \leq 1$. 
The total number of degrees of freedom $\mathbf{\tilde{X}}_{ij,l}$ is equal to $N_l=2 \cdot d$ and they are given by the space-time node coordinates which define the common face $\partial T_{ij}$ at time $t^n$ and $t^{n+1}$, as depicted in Figure \ref{fig.STelem}.

The Jacobian matrix of the transformation from the physical to the reference system can be written formally as
\begin{equation}
\mathcal{T}= \left( \begin{array}{c} \hat{\mathbf{e}} \\ \frac{\partial \mathbf{\tilde{x}}}{\partial \chi1} \\ \frac{\partial \mathbf{\tilde{x}}}{\partial \chi2} \\ \frac{\partial \mathbf{\tilde{x}}}{\partial \tau} \end{array} \right)^T = \left( \begin{array}{cccc} 
\hat{e}_1 & \hat{e}_2 & \hat{e}_3 & \hat{e}_4 \\ 
\frac{\partial x}{\partial \chi_1} & \frac{\partial y}{\partial \chi_1} & \frac{\partial z}{\partial \chi_1} & 0 \\
\frac{\partial x}{\partial \chi_2} & \frac{\partial y}{\partial \chi_2} & \frac{\partial z}{\partial \chi_2} & 0 \\
\frac{\partial x}{\partial \tau} & \frac{\partial y}{\partial \tau} & \frac{\partial z}{\partial \tau} & \Delta t 
 \end{array} \right)
\label{eqn.ChiTransf}
\end{equation}
where $\hat{e}_p$ represents the unit vector aligned with the $p$-th axis of the physical coordinate system $(x,y,z,t)$. The space-time volume $|\partial C_{ij}^n|$ as well as the space-time normal vector $\mathbf{\tilde n}_{ij}$ of $\partial C_{ij}^n$ can be conveniently obtained by computing the determinant of $\mathcal{T}$, that is
\begin{equation}
|\partial C_{ij}^n| = |\mathbf{\tilde n}_{ij}|, \quad \mathbf{\tilde n}_{ij} = \sum \limits_{p=1}^{d+1} \hat{e}_p \, (-1)^{1+p}|\mathcal{T}_{[1,p]}|
\label{eqn.detChi}
\end{equation}
with $|\mathcal{T}_{[1,p]}|$ representing the determinant of the cofactor matrix of $\hat{e}_p$.  
 
Once the space-time control volume $C_i^n$ has been determined, we write the governing balance law \eqref{eqn.PDE} in a more compact space-time divergence formulation as 
\begin{equation}
\tilde \nabla \cdot \tilde{\F} = \mathbf{0}  \qquad \textnormal{ with } \quad 
\tilde \nabla  = \left( \frac{\partial}{\partial x}, \, \frac{\partial}{\partial y}, \, \frac{\partial}{\partial z}, \, \frac{\partial}{\partial t} \right)^T, \quad \tilde{\F}  = \left( \mathbf{f}, \, \mathbf{g}, \, \mathbf{h}, \, \Q \right),
\label{eqn.st.pde}
\end{equation}
that is integrated over the space-time control volume $C^n_i$ yielding
\begin{equation}
 \int \limits_{\mathcal{C}^n_i} \tilde \nabla \cdot \tilde{\F} \, d\mathbf{x} dt = \mathbf{0}.
\label{STPDE}
\end{equation}
The integral in the above expression can be reformulated by means of Gauss theorem as an integral over the space-time surface $\partial \mathcal{C}^{n}_i$ of dimension $d$ defined by the outward pointing space-time unit normal vector $\mathbf{\tilde n}$, thus becoming
\begin{equation}
\int \limits_{\partial \mathcal{C}^{n}_i} \tilde{\F} \cdot \ \mathbf{\tilde n} \, dS = \mathbf{0}.
\label{STPDEgauss}
\end{equation}
Using the surface decomposition \eqref{eqn.dCi} and the mapping to the reference system $\boldsymbol{\chi}$ with \eqref{eqn.ChiTransf}-\eqref{eqn.detChi}, the final one-step ALE ADER finite volume scheme derives directly from \eqref{STPDEgauss} and writes
\begin{equation}
|T_i^{n+1}| \, \Q_i^{n+1} = |T_i^n| \, \Q_i^n - \sum \limits_{T_j \in \mathcal{N}_i} \,\, {\int \limits_0^1 \int \limits_0^1 \int \limits_{0}^{1-\chi_1} 
| \partial C_{ij}^n| \, \tilde{\G}_{ij} \cdot \mathbf{\tilde{n}}_{ij} \, d\chi_2 d\chi_1 d\tau}, 
\label{PDEfinal}
\end{equation}
with the numerical flux function $\tilde{\G}_{ij} \cdot \mathbf{\tilde{n}}_{ij}$ that properly takes into account the discontinuities of the predictor solution $\q_h$ on the element boundaries $\partial C_{ij}$ where a left state $\q_h^-$ in $T_i$ and a right state $\q_h^+$ in $T_j$ interact with each other. In \eqref{PDEfinal} the symbol $\mathcal{N}_i$ denotes the Neumann neighbors of element $T_i$, i.e. those elements $T_j^n$ which share a common face with $T_i^n$ across which the numerical fluxes are computed. The fluxes can be evaluated using a simple Rusanov-type scheme 
\begin{equation}
  \tilde{\G}_{ij} \cdot \mathbf{\tilde n}_{ij} =  
  \frac{1}{2} \left( \tilde{\F}(\q_h^+) + \tilde{\F}(\q_h^-)  \right) \cdot \mathbf{\tilde n}_{ij}  - 
  \frac{1}{2} s_{\max} \left( \q_h^+ - \q_h^- \right),  
  \label{eqn.rusanov} 
\end{equation}
where $s_{\max}$ denotes the maximum eigenvalue of the ALE Jacobian matrix in spatial normal direction which is given by
\begin{equation} 
\mathbf{A}^{\!\! \mathbf{V}}_{\mathbf{n}}(\Q):=\left(\sqrt{\tilde n_x^2 + \tilde n_y^2 + \tilde n_z^2 }\right)\left[\frac{\partial \mathbf{F}}{\partial \Q} \cdot \mathbf{n}  - 
(\mathbf{V} \cdot \mathbf{n}) \,  \mathbf{I}\right], \qquad    
\mathbf{n} = \frac{(\tilde n_x, \tilde n_y, \tilde n_z)^T}{\sqrt{\tilde n_x^2 + \tilde n_y^2 + \tilde n_z^2 }},  
\end{equation} 
with the local normal mesh velocity $\mathbf{V} \cdot \mathbf{n}$ and the identity matrix $\mathbf{I}$. The Rusanov flux is very robust and can handle quite challenging and strong discontinuities by supplying numerical dissipation to the scheme. A less dissipative flux function, namely a generalization of the Osher-Solomon scheme \cite{osherandsolomon}, has been successfully proposed in 
\cite{OsherUniversal,OsherNC,ApproxOsher} and subsequently applied also in the moving mesh context \cite{Lagrange3D}. In this case the numerical flux takes the form  
\begin{equation}
  \tilde{\G}_{ij} \cdot \mathbf{\tilde n}_{ij} =  
  \frac{1}{2} \left( \tilde{\F}(\q_h^+) + \tilde{\F}(\q_h^-)  \right) \cdot \mathbf{\tilde n}_{ij}  - 
  \frac{1}{2} \left( \int \limits_0^1 \left| \mathbf{A}^{\!\! \mathbf{V}}_{\mathbf{n}}(\boldsymbol{\Psi}(s)) \right| ds \right) \left( \q_h^+ - \q_h^- \right),  
  \label{eqn.osher} 
\end{equation} 
where $\boldsymbol{\Psi}(s)$ is the straight-line segment path $\boldsymbol{\Psi}(s) = \q_h^- + s \left( \q_h^+ - \q_h^- \right)$ with $s\in[0,1]$ that connects the left and the right state across the discontinuity. The absolute value of the dissipation matrix in \eqref{eqn.osher} is evaluated as usual as 
\begin{equation}
 |\mathbf{A}| = \mathbf{R} |\boldsymbol{\Lambda}| \mathbf{R}^{-1},  \qquad |\boldsymbol{\Lambda}| = \textnormal{diag}\left( |\lambda_1|, |\lambda_2|, ..., |\lambda_\nu| \right),  
\end{equation}
with $\mathbf{R}$ and $\mathbf{R}^{-1}$ representing the right eigenvector matrix and its inverse, respectively.  

For computing the integrals we adopt Gaussian quadrature formulae of suitable order, see \cite{stroud}. However, in the Eulerian case where the mesh does not move in time, the Jacobians $|\partial C_{ij}^n|$ as well as the space-time normal vectors $\mathbf{\tilde{n}}_{ij}$ in \eqref{PDEfinal} are \textit{constant}, therefore a very efficient \textit{quadrature-free} formulation can be derived as fully detailed in \cite{Dumbser2007204}: the Galerkin predictor stage provides the constant degrees of freedom $\widehat{\q}_{l}$ and $\widehat{\F}_{l}$, hence allowing the space-time integrals appearing in \eqref{PDEfinal} to be precomputed and stored once and for all on the space-time reference element $T_E \times [0,1]$ using the isoparametric approximations \eqref{eqn.PDEterm}. Thus, the integrals are simply and efficiently obtained by one single multiplication between the degrees of freedom and the precomputed integrals of the basis functions $\theta_l$ on the space-time reference element that does never change. For recent quadrature-free high order Lagrangian algorithms, see \cite{LagrangeQF}. 

Finally, we underline that the numerical scheme \eqref{PDEfinal} satisfies the so-called \textit{geometrical conservation law} (GCL) \textit{by construction}, since an integration over a closed 
space-time control volume is carried out and application of Gauss theorem directly leads to 
\begin{equation}
 \int_{\partial \mathcal{C}_i^n} \mathbf{\tilde n} \, dS = 0.
 \label{eqn.gcl} 
\end{equation}
For more details, see also the appendix in \cite{Lagrange3D}. 

\subsection{Some comments on the implementation}

ADER finite volume schemes are high order accurate and fully-discrete one-step schemes. The predictor step for the computation of the space-time polynomials $\mathbf{q}_h(\mathbf{x},t)$ is done in a completely \textit{local} manner, without the need of any MPI communication. Further to that, the step is arithmetically very intensive, but without requiring access to large and disjoint memory 
patterns. Instead, the iterative procedure of the predictor step works only on the same element-local space-time degrees of freedom of the state and the fluxes inside each element. 
Hence, it is very well suited for the current cache-based CPU architectures with high arithmetic velocity but fairly slow connection to the main memory (RAM). On fixed meshes, the entire scheme 
can be written in a 
quadrature-free manner \cite{Dumbser2007204}, requiring essentially only matrix-matrix multiplications, for the reconstruction step \eqref{CWENO:Popt}, for the predictor (local time-evolution) stage 
\eqref{eqn.CGfinal} as well as for the integration of the numerical fluxes on the element boundaries in the final finite volume scheme \eqref{PDEfinal}. This task can be easily optimized at the aid
of standard linear algebra packages (BLAS), e.g. using the Intel Math Kernel Library (MKL) or the particularly optimized Intel library \texttt{libxsmm} for small matrix-matrix multiplications 
\cite{libxsmm}. 
Thanks to the fully-discrete one-step nature of ADER schemes, the CWENO reconstruction is carried out only once per time step, for any order of accuracy in time, while classical 
Runge-Kutta time-stepping requires the reconstruction (and the associated MPI communication) to be carried out in each Runge-Kutta substage again. ADER schemes can thus be called 
genuinely \textit{communication avoiding} methods. Concerning the generation and partition of very big unstructured meshes, we proceed as follows. First, a coarse unstructured mesh is 
generated with only a few million elements. This coarse mesh is then partitioned and distributed among the MPI ranks using the free graph partitioning 
software Metis/ParMetis \cite{metis}. Next, our algorithm proceeds with the generation of a refined mesh, producing for each coarse element a fine subgrid with nodes 
and connectivity given in \cite{DGLimiter3}. This mesh refinement step is done in parallel, where each MPI rank refines its own elements 
and their Voronoi neighbors, i.e. those elements sharing a common node. This is required in order to produce a sufficient overlap with the neighbor CPUs, necessary 
for building the CWENO reconstruction stencils. In this way our implementation of the unstructured CWENO scheme generates the final unstructured mesh in parallel and 
is thus able to handle hundreds of millions of elements, leading to billions of degrees of freedom of the resulting CWENO polynomials, as shown on two examples in the 
following section.

\section{Test problems}
\label{sec.validation} 

In the following we present a set of test problems in order to validate the ADER-CWENO finite volume schemes presented in this paper. We run our simulations using either the Eulerian or the Arbitrary-Lagrangian-Eulerian (ALE) version of the scheme and this section is split accordingly. In the ALE framework we set the mesh velocity to be equal to the local fluid velocity, hence $\mathbf{V}=\mathbf{v}$, and for all the test cases both the order of the simulation $M+1$ as well as the numerical flux function that has been used are explicitly written.

\paragraph{Euler equations} The first hyperbolic system of the form \eqref{eqn.PDE} is given by the Euler equations of compressible gas dynamics with 
\begin{equation}
\label{eulerTerms}
\Q = \left( \begin{array}{c} \rho   \\ \rho u  \\ \rho v \\ \rho w \\ \rho E \end{array} \right), \quad
\mathbf{F} = \left( \begin{array}{ccc}  \rho u       & \rho v        & \rho w       \\ 
                                        \rho u^2 + p & \rho u v      & \rho w u     \\
																			  \rho u v     & \rho v^2 + p  & \rho w v     \\ 
																				\rho u w     & \rho v w      & \rho w^2 + p \\ 
																			 u(\rho E + p) & v(\rho E + p) & w(\rho E + p)  
										\end{array} \right).  
\end{equation}
The vector of conserved variables $\Q$ involves the fluid density $\rho$, the momentum density vector $\rho \v=(\rho u, \rho v,\rho w)$ and the total energy density $\rho E$. 
The fluid pressure $p$ is related to the fluid velocity vector $\mathbf{v}=(u,v,w)$ using the equation of state for an ideal gas 
\begin{equation}
\label{eqn.eos} 
p = (\gamma-1)\left(\rho E - \frac{1}{2} \rho \mathbf{v}^2 \right), 
\end{equation}
where $\gamma$ is the ratio of specific heats so that the speed of sound takes the form $c=\sqrt{\frac{\gamma p}{\rho}}$. To assign the initial condition for the test problems discussed in this article we may also use the vector of primitive variables $\U=(\rho,u,v,w,p)$.

\paragraph{Ideal MHD equations} We also consider the equations of ideal classical magnetohydrodynamics (MHD) that result in a more complicated hyperbolic conservation law. Here, the state vector $\Q$ and the flux tensor $\F=(\f,\g,\h)$ read
\begin{equation}
\label{MHDTerms} 
  \Q = \left( \begin{array}{c} \rho \\   \rho \v \\   \rho E \\ \B \\ \psi \end{array} \right), \quad
  \F(\Q) = \left( \begin{array}{c} 
  \rho \v  \\ 
  \rho \v \v + p + p_m \mathbf{I} - \frac{1}{4 \pi} \B \B \\ 
  \v (\rho E + p + p_m ) - \frac{1}{4 \pi} \B ( \v \cdot \B ) \\ 
  \v \B - \B \v + \psi \mathbf{I} \\
  c_h^2 \B 
  \end{array} \right),
\end{equation} 
where the total pressure  $p+p_m$ is computed as the sum of the hydrodynamic pressure $p$ and the magnetic pressure $p_m = \frac{1}{8 \pi} \B^2$; $\mathbf{I}$ is the identity matrix 
and the rest of the notation is the same employed for the Euler equations. System \eqref{eqn.PDE} with \eqref{MHDTerms} is again closed by the ideal gas equation of state 
\begin{equation}
p = \left(\gamma - 1 \right) \left(\rho E - \frac{1}{2}\mathbf{v}^2 - \frac{\mathbf{B}^2}{8\pi}\right).
\label{MHDeos}
\end{equation}
In \eqref{MHDTerms} one more linear scalar PDE has been added to the system with the additional variable $\psi$ in order to transport the divergence errors out of the computational domain with an artificial divergence cleaning speed $c_h$, according to the hyperbolic version of the generalized Lagrangian multiplier (GLM) divergence cleaning approach \cite{Dedneretal}. This strategy is needed because the MHD equations require a constraint on the divergence of the magnetic field $\B$ to be respected, that is 
\begin{equation}
\nabla \cdot \mathbf{B} = 0.
\label{eqn.divB}
\end{equation}
Although at the continuous level it is enough to initialize the magnetic field with divergence-free data, at the discrete level condition \eqref{eqn.divB} might be in principle violated. 
Another possible technique that ensures such constraint to be fulfilled is presented in \cite{MHDdivFree2015}. Similar to the Euler equations, the vector of primitive variables for the MHD 
system reads $\U=(\rho,u,v,w,p,B_x,B_y,B_z,\psi)$.

\subsection{Numerical convergence studies}
\label{sec.conv}
The numerical convergence of the new ADER-CWENO schemes is studied by considering a test problem proposed in \cite{HuShuTri} for the Euler equations of compressible gas dynamics. A smooth isentropic vortex is convected on the horizontal plane $x-y$ with velocity $\v_c=(1,1,0)$ and the initial condition is given in primitive variables as a superposition of a homogeneous background field and some perturbations $\delta$, that is
\begin{equation}
\U = (1+\delta \rho, 1+\delta u, 1+\delta v, 1+\delta w, 1+\delta p),
\label{eqn.ShuVortIC}
\end{equation}
with
\begin{eqnarray}
\label{ShuVortDelta}
\left(\begin{array}{c} \delta u \\ \delta v \\ \delta w \end{array}\right) = \frac{\epsilon}{2\pi}e^{\frac{1-r^2}{2}} \left(\begin{array}{c} -(y-5) \\ \phantom{-}(x-5) \\ 0 \end{array}\right), & \quad \delta T = -\frac{(\gamma-1)\epsilon^2}{8\gamma\pi^2}e^{1-r^2}, \nonumber \\
\quad \delta \rho = (1+\delta T)^{\frac{1}{\gamma-1}}-1, & \quad \delta p = (1+\delta T)^{\frac{\gamma}{\gamma-1}}-1.
\end{eqnarray}
$\delta T$ is the perturbation for temperature, the vortex strength is set to $\epsilon=5$ and the ratio of specific heats is $\gamma=1.4$. According to \cite{HuShuTri} we assume that the entropy fluctuation is zero, hence $\delta S=0$ with $S=\log \left(\frac{p}{\rho^\gamma(\gamma-1)}\right)$. The generic radial position on the $x-y$ plane is denoted by $r^2=(x-5)^2+(y-5)^2$ and the initial computational domain is the square $\Omega(0)=[0;10]\times[0;10]$ in 2D and the box $\Omega(0)=[0;10]\times[0;10]\times[0;5]$ in 3D, where periodic boundaries have been set everywhere. The final time of the simulation is chosen to be $t_f=1$ and the exact solution is given by the time-shifted initial condition as $\Q_e(\x,t_f)=\Q(\x-\v_c t_f,0)$.
The high order reconstructed solution $\w_i(\x,t_f)$ is used to measure the error at time $t_f$ adopting the continuous $L_2$ norm as
\begin{equation}
  \epsilon_{L_2} = \sqrt{ \int \limits_{\Omega(t_f)} \left( \Q_e(\x,t_f) - \w_i(\x,t_f) \right)^2 d\x },  
	\label{eqnL2error}
\end{equation}
where $h(\Omega(t_f))$ denotes the characteristic mesh size that is assumed to be the maximum diameter of the circumcircles or the circumspheres of the control volumes at time $t_f$. A sequence of successively refined grids is employed to run this test problem and Tables \ref{tab.convEUL} and \ref{tab.convALE} report the numerical convergence studies for the Eulerian and the ALE setting
on fixed and moving grids, respectively. The Osher-type flux \eqref{eqn.osher} has been used in all computations and the results show that the desired order of accuracy 
is obtained in both space and time up to $\mathcal{O}(5)$. Note that the ADER approach is also uniformly high order accurate in time, hence all tests can be run with the maximum admissible Courant 
number of CFL$=1/d$.

\begin{table}[!htbp]  
\caption{Numerical convergence results for the multidimensional compressible Euler equations using the ADER-CWENO finite volume schemes on \textit{fixed meshes} from third up to fifth order of accuracy. The error norms refer to the variable $\rho$ (density) at time $t=10$.}  
\begin{center} 
\begin{small}
\renewcommand{\arraystretch}{1.0}
\begin{tabular}{c|cccccc} 
\multicolumn{7}{c}{} \\
\hline
         \multicolumn{7}{c}{Eulerian schemes on fixed meshes} \\
\hline
\textbf{2D}  & \multicolumn{2}{c}{$\mathcal{O}3$} & \multicolumn{2}{c}{$\mathcal{O}4$}  & \multicolumn{2}{c}{$\mathcal{O}5$}  \\
\hline
  $h(\Omega(t_f))$ & $\epsilon_{L_2}$ & $\mathcal{O}(L_2)$ & $\epsilon_{L_2}$ & $\mathcal{O}(L_2)$ &  $\epsilon_{L_2}$ & $\mathcal{O}(L_2)$ \\ 
\hline
2.48E-01 & 9.7789E-02 & -   & 4.9493E-02 & -   & 6.1176E-02 & -    \\ 
1.28E-01 & 1.9593E-02 & 2.4 & 1.7276E-03 & 5.1 & 1.5472E-03 & 5.6  \\ 
6.31E-02 & 2.8801E-03 & 2.7 & 1.1273E-04 & 3.9 & 6.1297E-05 & 4.6  \\ 
3.21E-02 & 3.7076E-04 & 3.0 & 6.8197E-06 & 4.1 & 1.9332E-06 & 5.1  \\ 
\hline 
				 \multicolumn{7}{c}{} \\
\hline
\textbf{3D}  & \multicolumn{2}{c}{$\mathcal{O}3$} & \multicolumn{2}{c}{$\mathcal{O}4$}  & \multicolumn{2}{c}{$\mathcal{O}5$}  \\
\hline
  $h(\Omega(t_f))$ & $\epsilon_{L_2}$ & $\mathcal{O}(L_2)$ & $\epsilon_{L_2}$ & $\mathcal{O}(L_2)$ &  $\epsilon_{L_2}$ & $\mathcal{O}(L_2)$ \\ 
\hline
5.92E-01 & 1.2824E-01 & -   & 1.1112E-01 & -   & 1.8716E-03 & -    \\ 
3.61E-01 & 3.9004E-02 & 2.4 & 1.9492E-02 & 3.5 & 6.7054E-02 & 2.1  \\ 
2.31E-01 & 1.2255E-02 & 2.6 & 3.4648E-03 & 3.8 & 8.0960E-03 & 4.7  \\ 
1.80E-01 & 5.9360E-03 & 2.9 & 1.3851E-03 & 3.7 & 1.9927E-03 & 5.7  \\ 
\hline
\end{tabular}
\end{small}
\end{center}
\label{tab.convEUL}
\end{table}

\begin{table}[!htbp]  
\caption{Numerical convergence results for the multidimensional compressible Euler equations using the ALE ADER-CWENO finite volume schemes on \textit{moving meshes} from third up to fifth order of accuracy. The error norms refer to the variable $\rho$ (density) at time $t=1$.}  
\begin{center} 
\begin{small}
\renewcommand{\arraystretch}{1.0}
\begin{tabular}{c|cccccc} 
\multicolumn{7}{c}{} \\
\hline
         \multicolumn{7}{c}{Direct Arbitrary-Lagrangian-Eulerian (ALE) schemes on moving meshes} \\
\hline
\textbf{2D}  & \multicolumn{2}{c}{$\mathcal{O}3$} & \multicolumn{2}{c}{$\mathcal{O}4$}  & \multicolumn{2}{c}{$\mathcal{O}5$}  \\
\hline
  $h(\Omega(t_f))$ & $\epsilon_{L_2}$ & $\mathcal{O}(L_2)$ & $\epsilon_{L_2}$ & $\mathcal{O}(L_2)$ &  $\epsilon_{L_2}$ & $\mathcal{O}(L_2)$ \\ 
\hline
3.30E-01 & 1.6270E-02 & -   & 4.4800E-03 & -   & 4.5492E-03 & -    \\ 
2.51E-01 & 7.0051E-03 & 3.1 & 1.7353E-03 & 3.5 & 1.2842E-03 & 4.7  \\ 
1.68E-01 & 2.3028E-03 & 2.7 & 4.3117E-04 & 3.4 & 2.2611E-04 & 4.3  \\ 
1.28E-01 & 9.3371E-04 & 3.3 & 1.3580E-04 & 4.3 & 5.8177E-05 & 5.0  \\ 
\hline 
				 \multicolumn{7}{c}{} \\
\hline
\textbf{3D}  & \multicolumn{2}{c}{$\mathcal{O}3$} & \multicolumn{2}{c}{$\mathcal{O}4$}  & \multicolumn{2}{c}{$\mathcal{O}5$}  \\
\hline
  $h(\Omega(t_f))$ & $\epsilon_{L_2}$ & $\mathcal{O}(L_2)$ & $\epsilon_{L_2}$ & $\mathcal{O}(L_2)$ &  $\epsilon_{L_2}$ & $\mathcal{O}(L_2)$ \\ 
\hline
5.92E-01 & 9.9325E-02 & -   & 4.3194E-02 & -   & 3.3929E-03 & -    \\ 
3.61E-01 & 3.0677E-02 & 2.4 & 8.3789E-03 & 3.3 & 5.4086E-03 & 3.7  \\ 
2.31E-01 & 9.0675E-03 & 2.7 & 1.7496E-03 & 3.5 & 8.6444E-04 & 4.1  \\ 
1.80E-01 & 4.1534E-03 & 3.1 & 6.0517E-04 & 4.3 & 2.5793E-04 & 4.9  \\ 
\hline
\end{tabular}
\end{small}
\end{center}
\label{tab.convALE}
\end{table}

\clearpage
\subsection{Numerical results on fixed meshes (Eulerian schemes)}
\label{ssec.EUL}

\subsubsection{Riemann problems}
\label{ssec.RP}
\begin{figure}[!htbp]
\begin{center}
\begin{tabular}{cc} 
\includegraphics[width=0.44\textwidth]{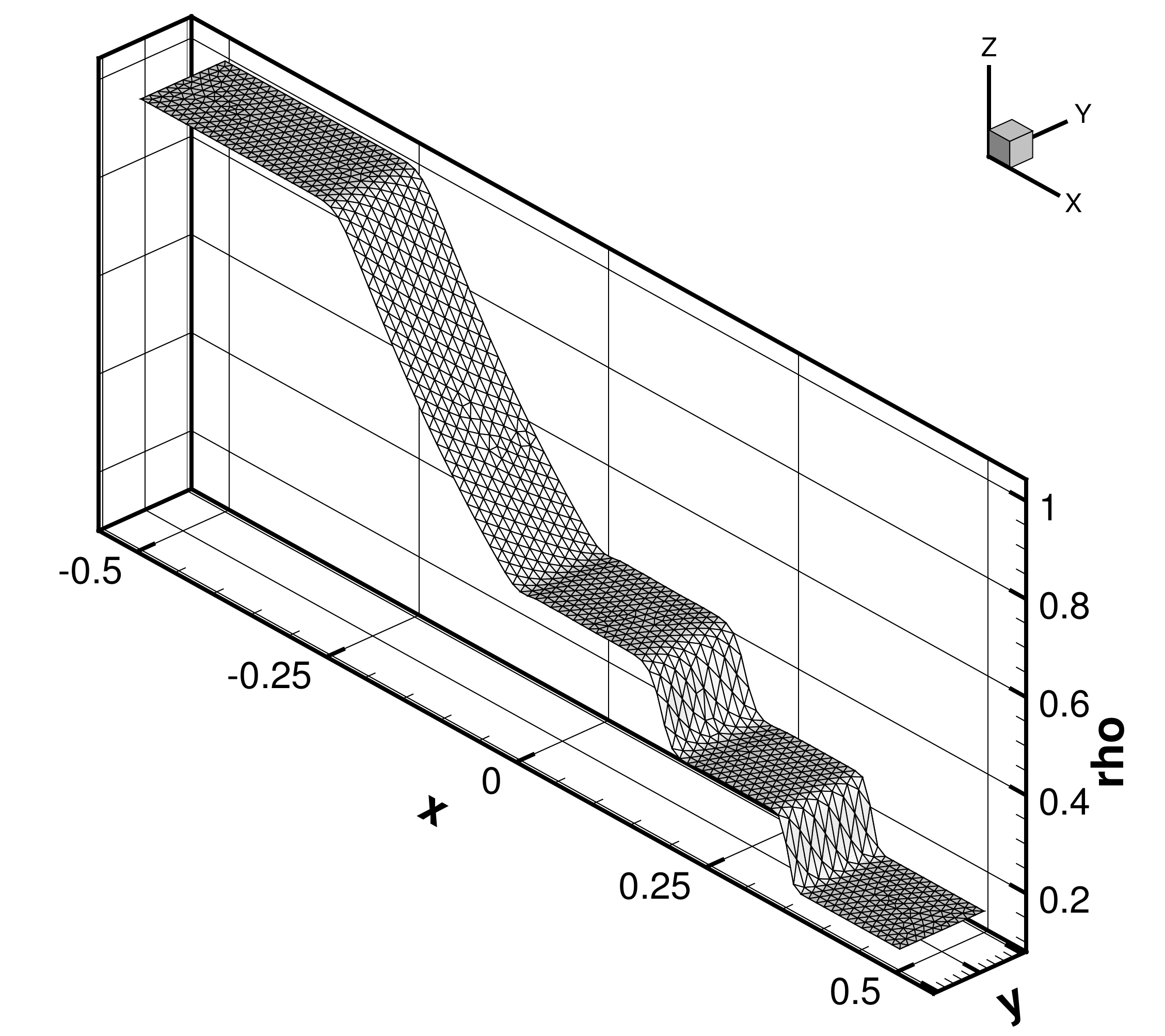}  &           
\includegraphics[width=0.44\textwidth]{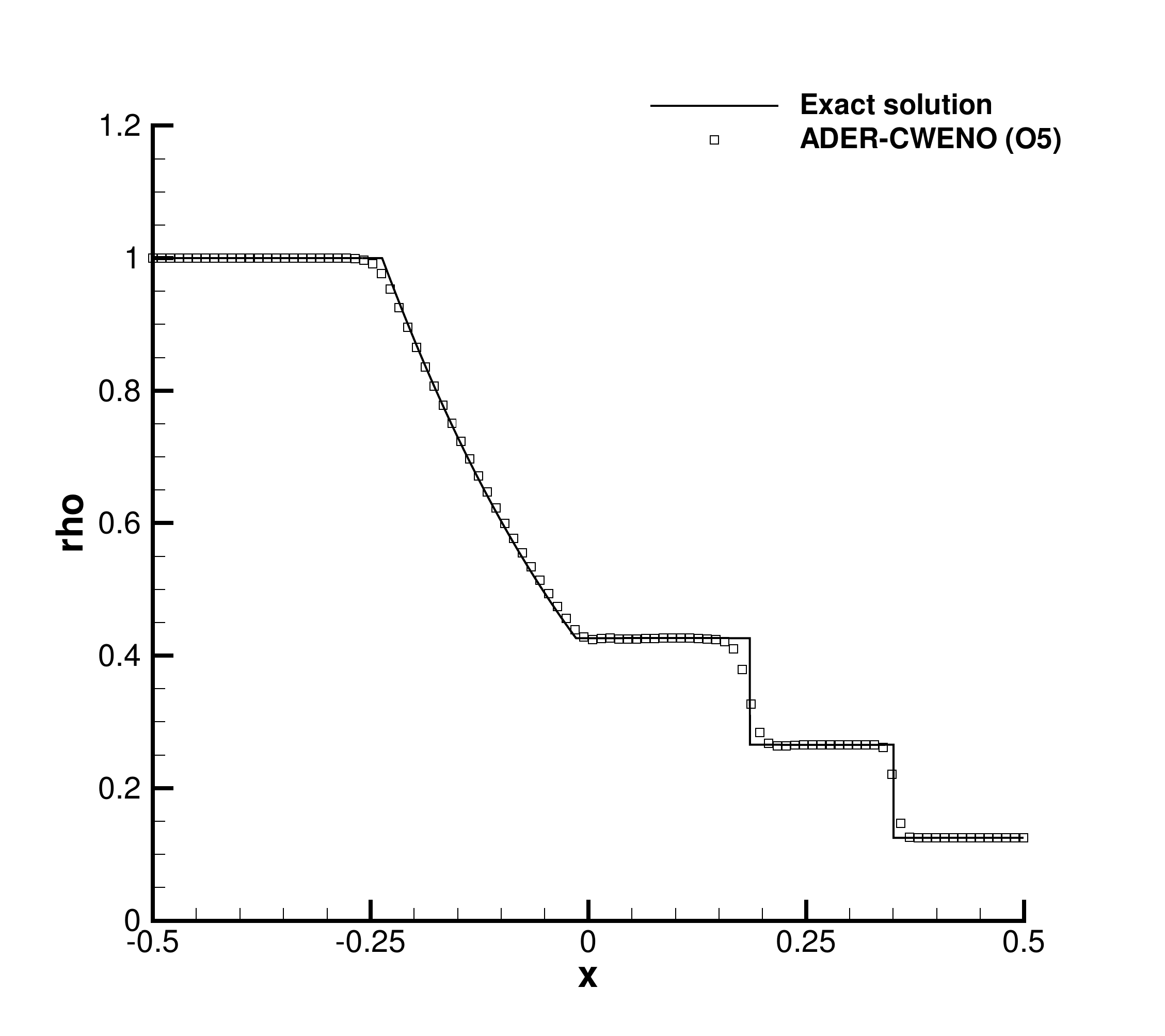} \\
\includegraphics[width=0.44\textwidth]{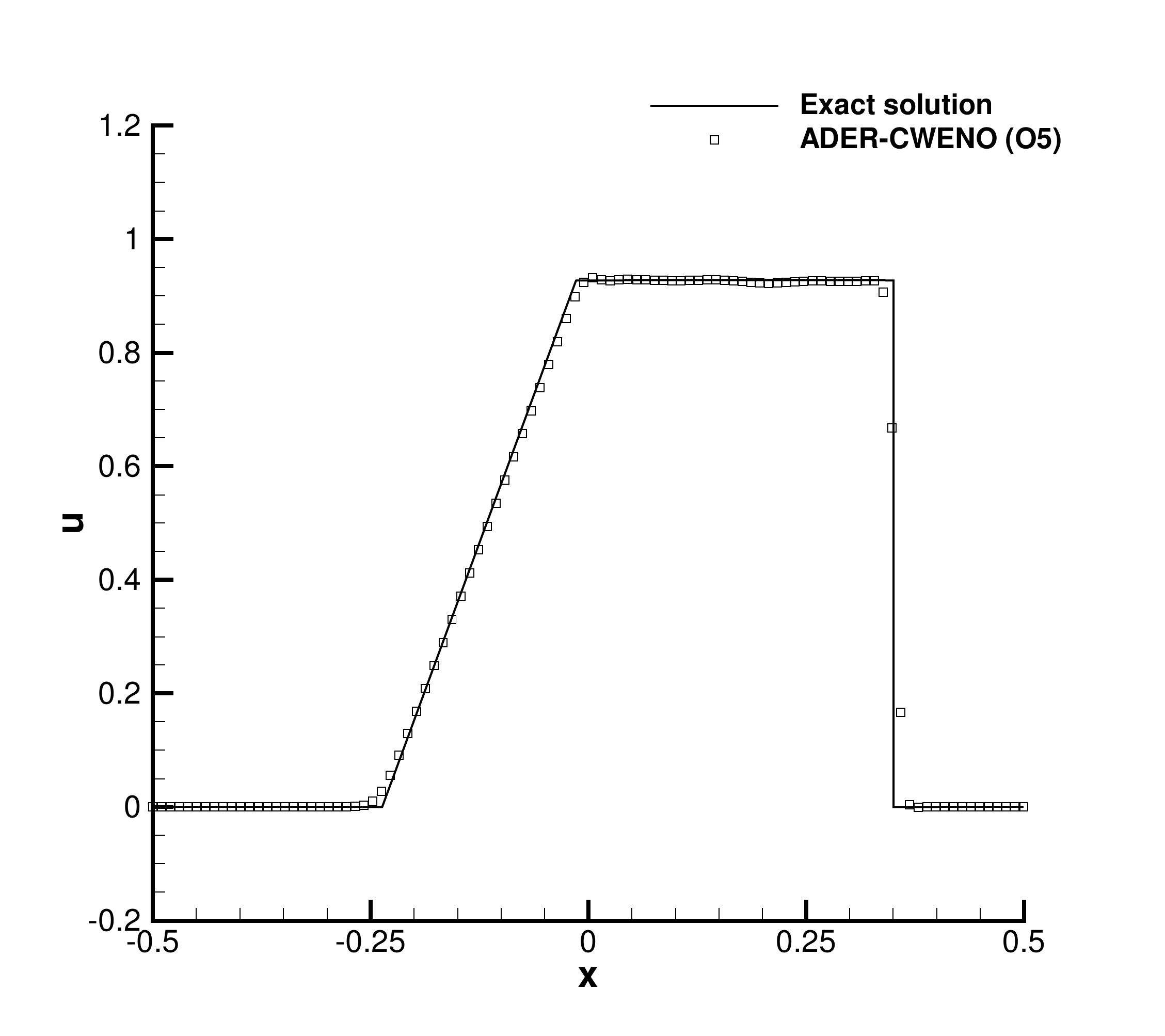}  &           
\includegraphics[width=0.44\textwidth]{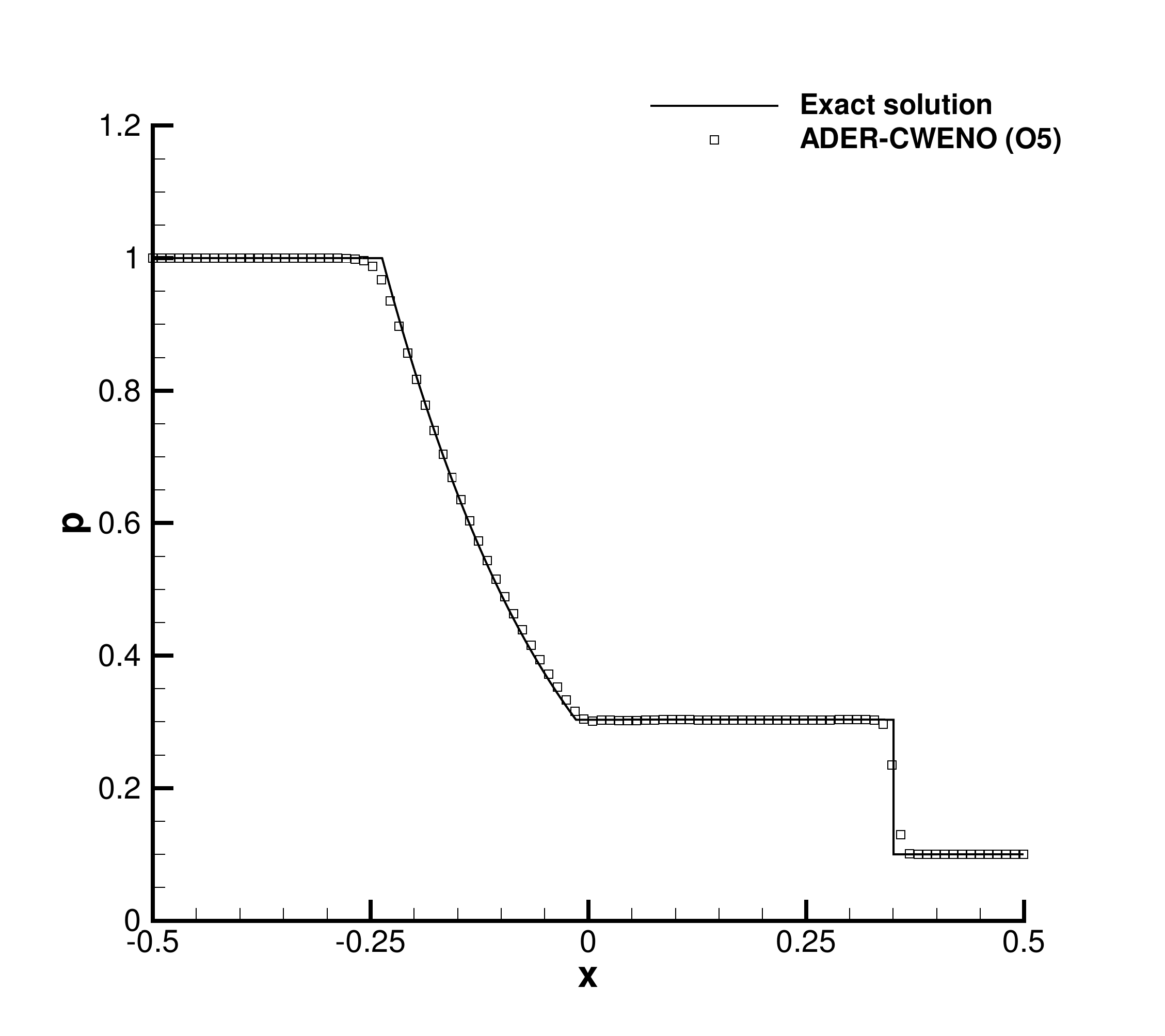} \\
\end{tabular} 
\caption{Numerical results for the Sod problem ($M=4$) at time $t=0.2$: 3D view of the density distribution and comparison against the exact solution for density, velocity $u$ and pressure.}
\label{fig.Sod}
\end{center}
\end{figure}

\begin{figure}[!htbp]
\begin{center}
\begin{tabular}{cc} 
\includegraphics[width=0.44\textwidth]{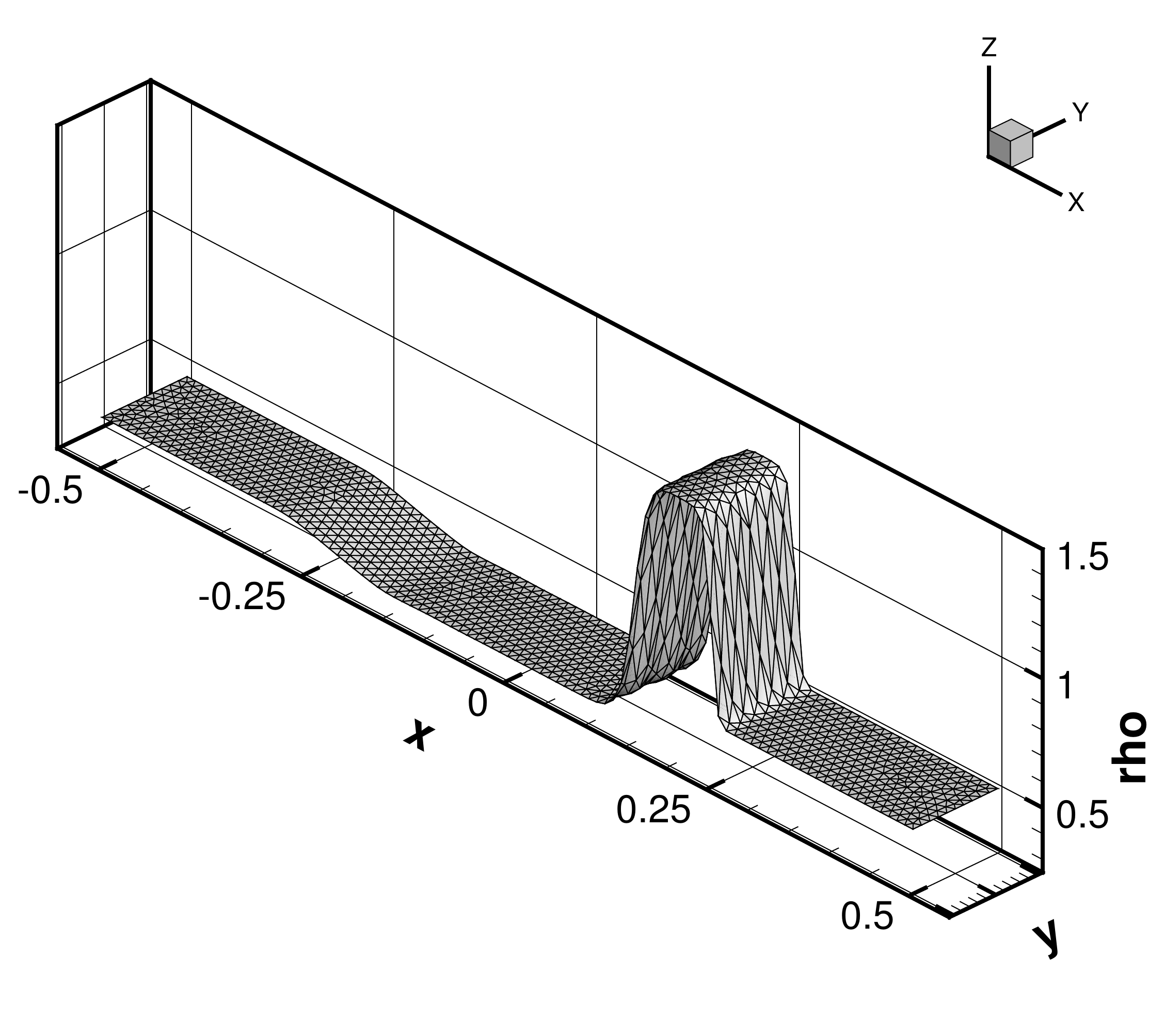}  &           
\includegraphics[width=0.44\textwidth]{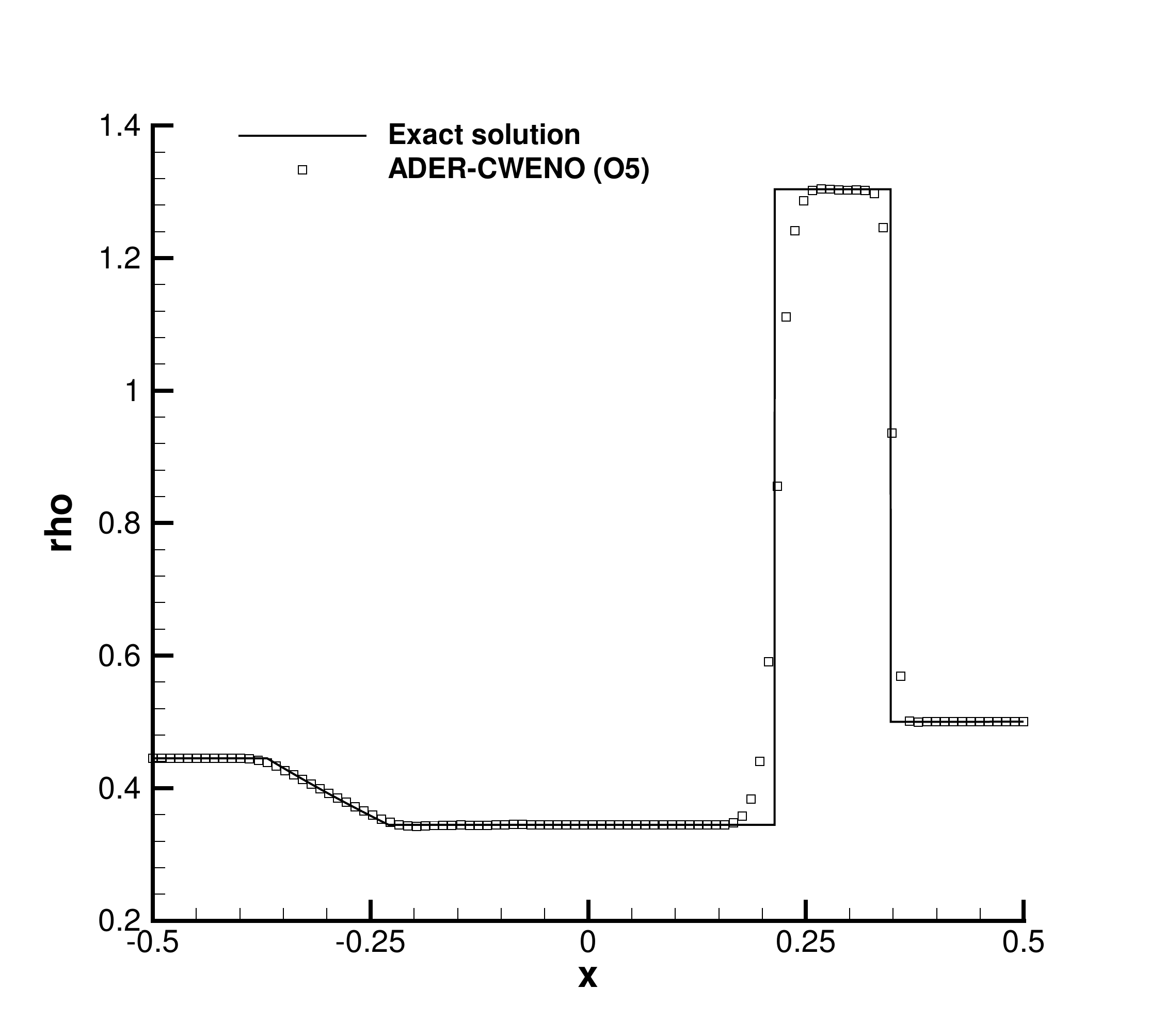} \\
\includegraphics[width=0.44\textwidth]{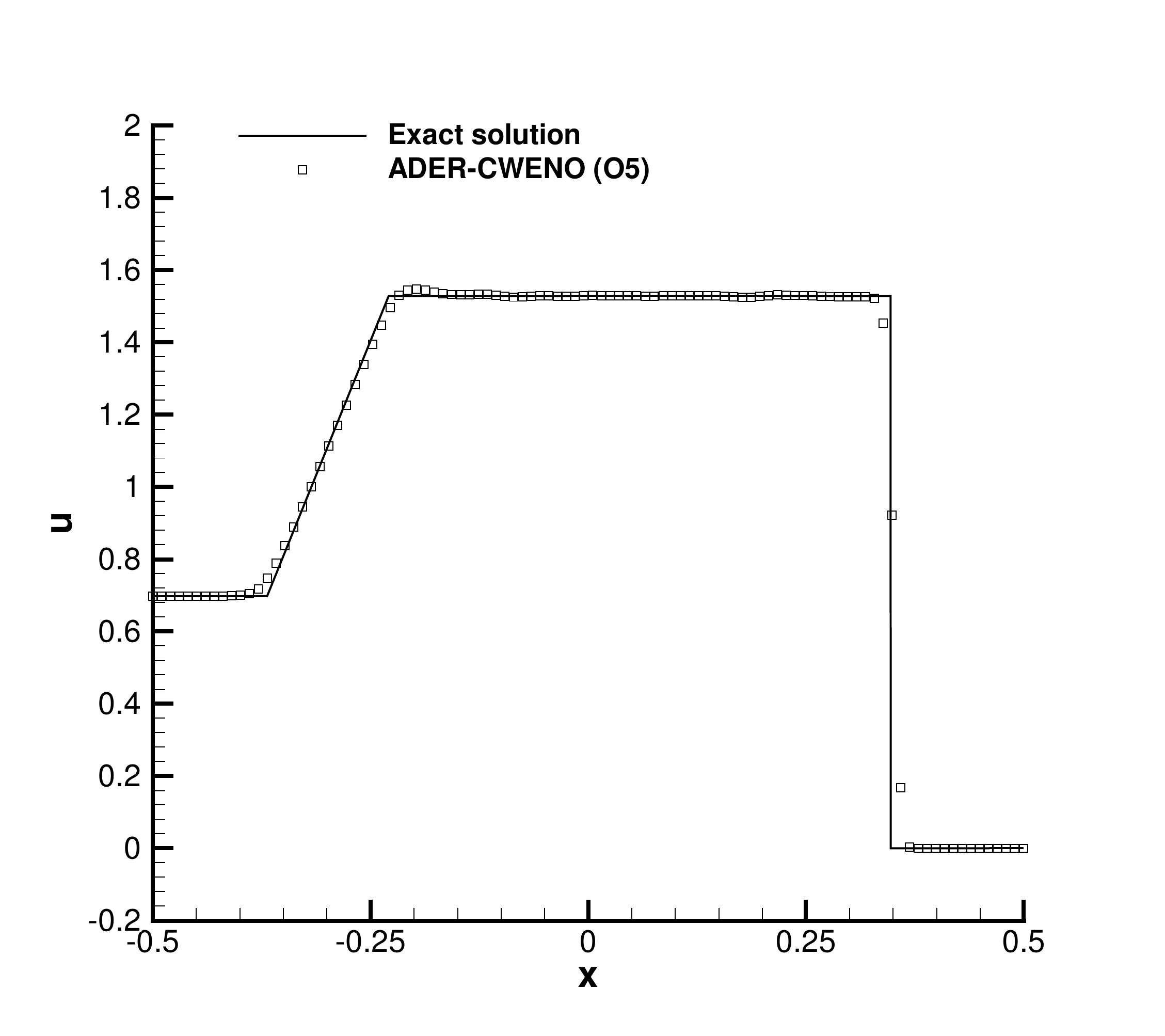}  &           
\includegraphics[width=0.44\textwidth]{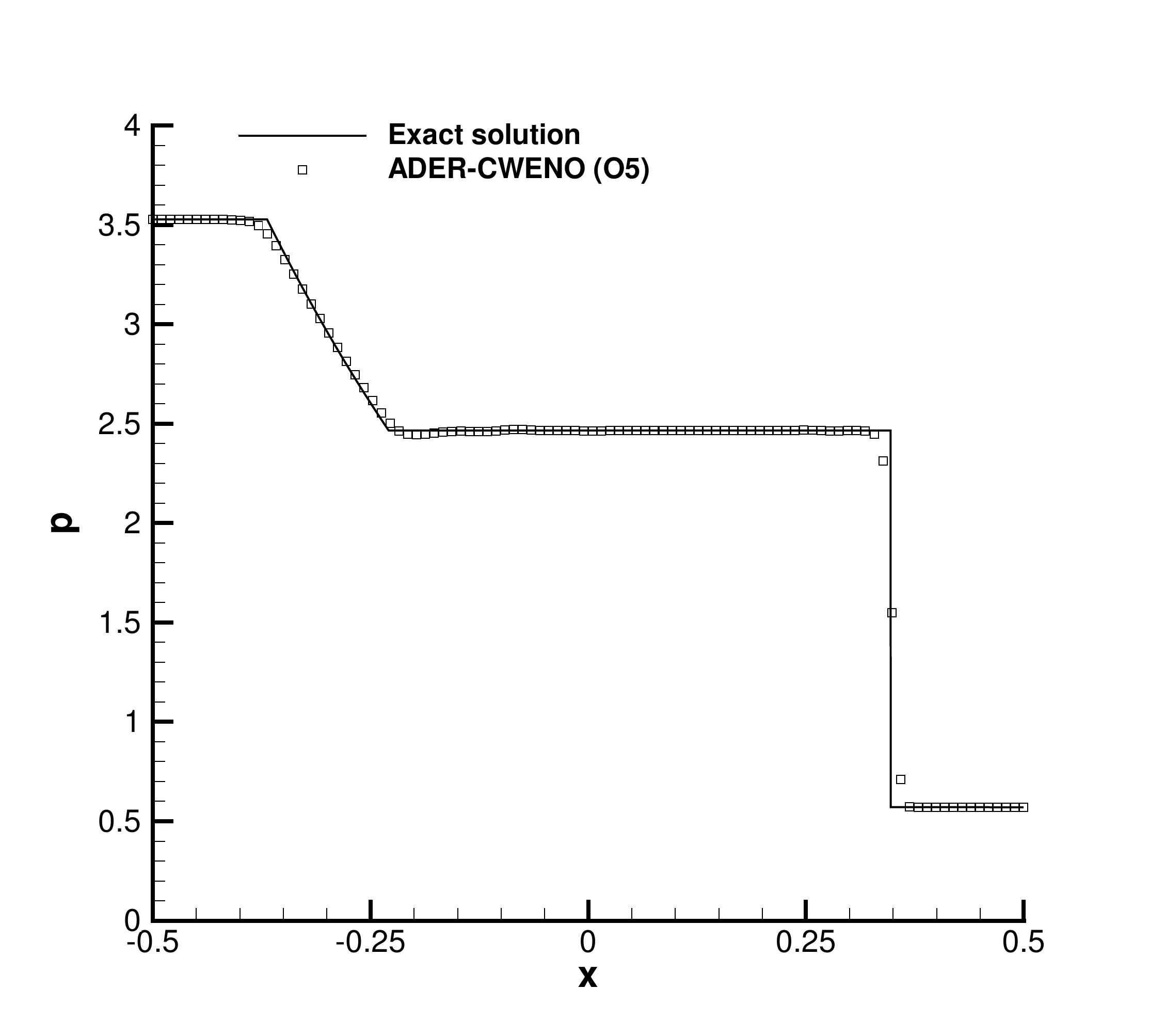} \\
\end{tabular} 
\caption{Numerical results for the Lax problem ($M=4$) at time $t=0.14$: 3D view of the density and comparison against the exact solution for density, velocity and pressure.}
\label{fig.Lax}
\end{center}
\end{figure}

\begin{figure}[!htbp]
\begin{center}
\begin{tabular}{cc} 
\includegraphics[width=0.44\textwidth]{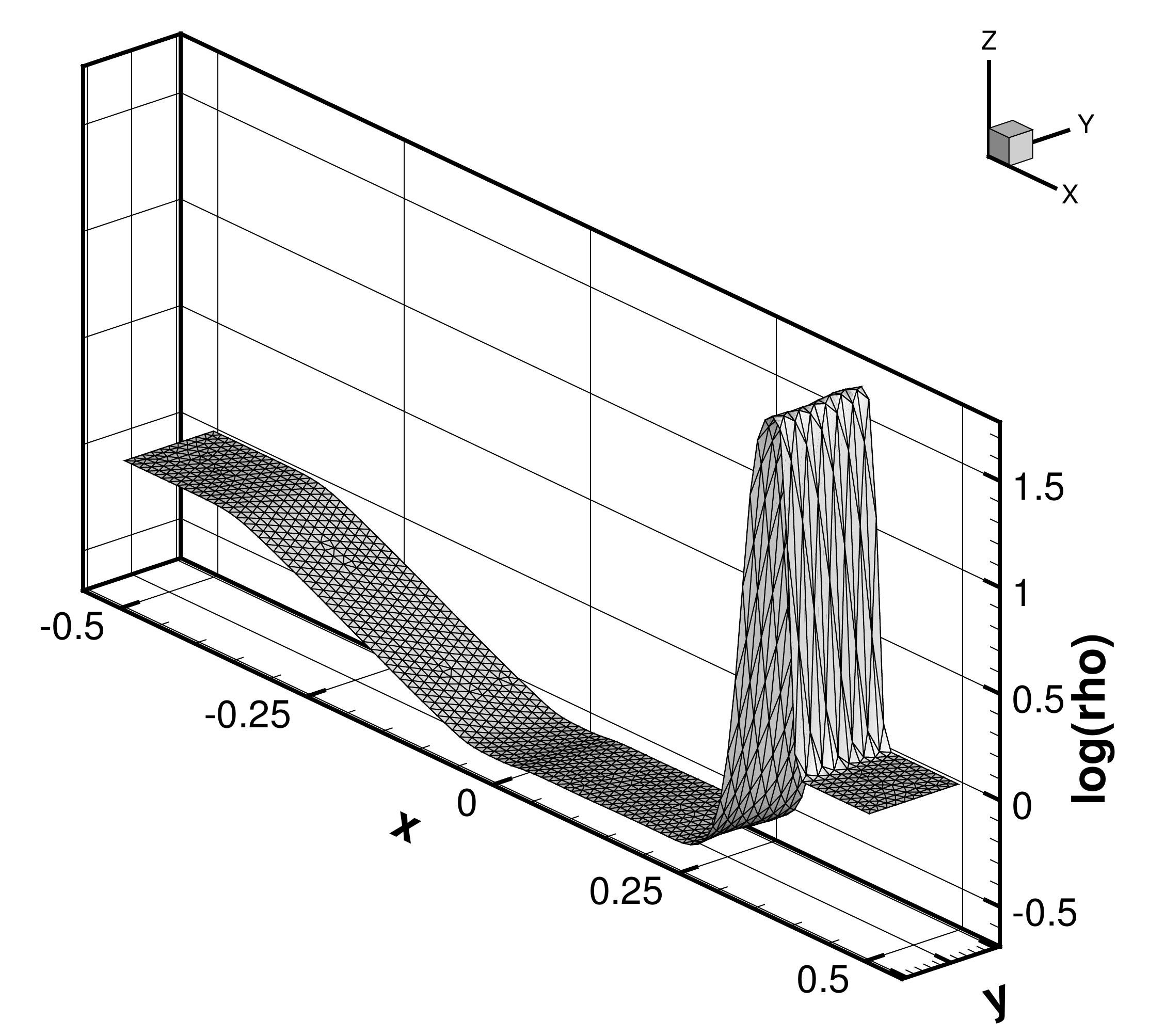}  &           
\includegraphics[width=0.44\textwidth]{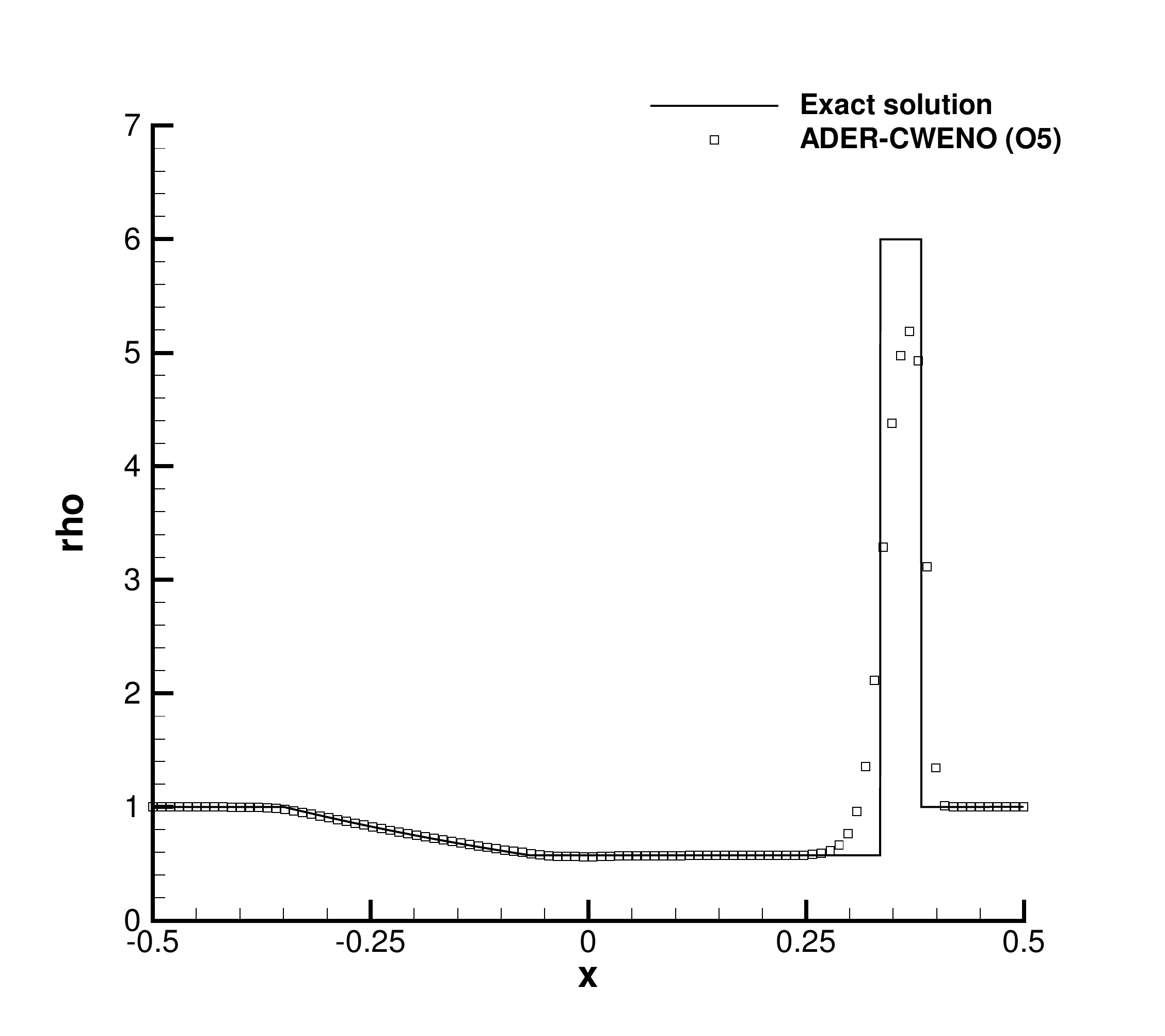} \\
\includegraphics[width=0.44\textwidth]{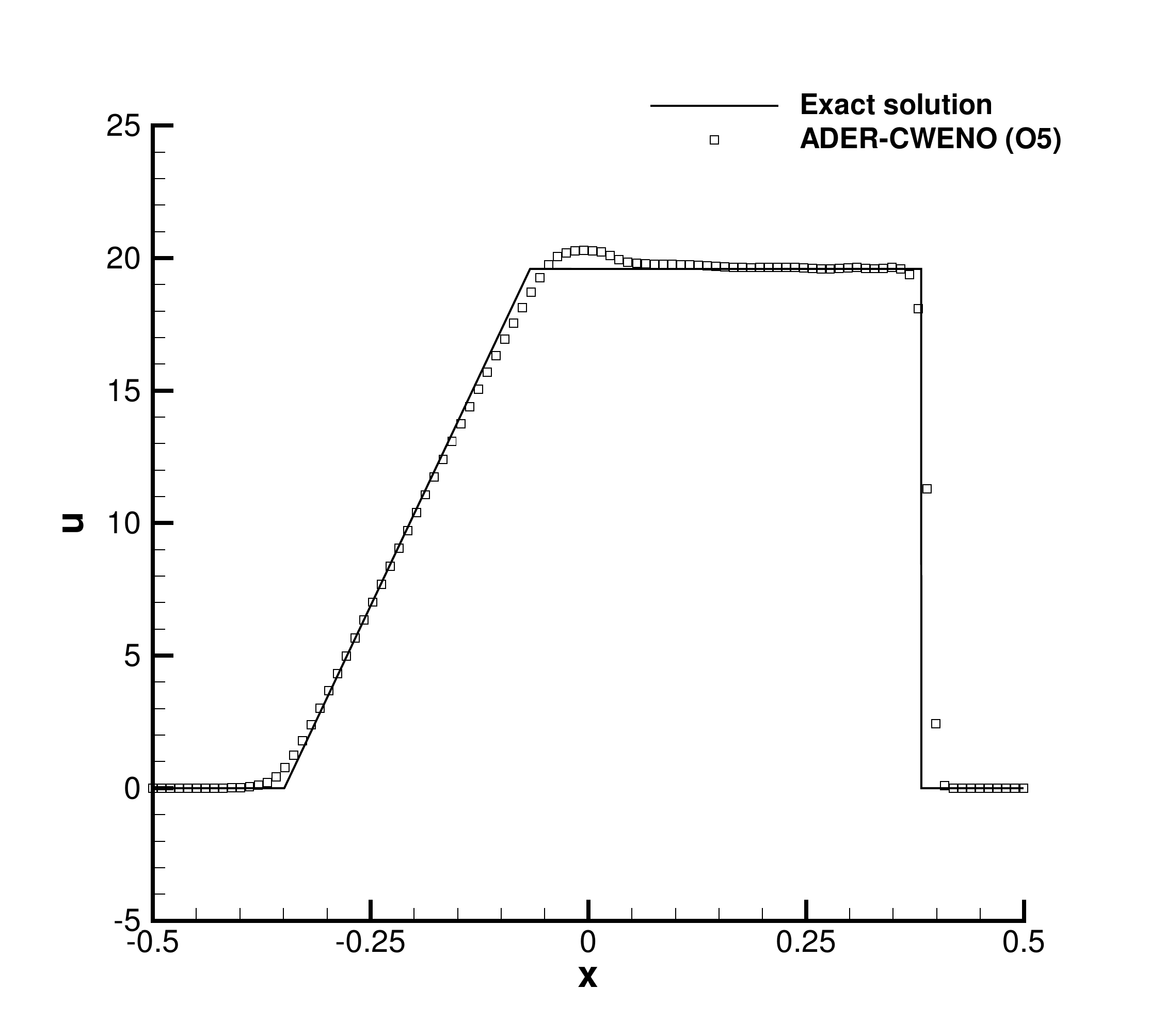}  &           
\includegraphics[width=0.44\textwidth]{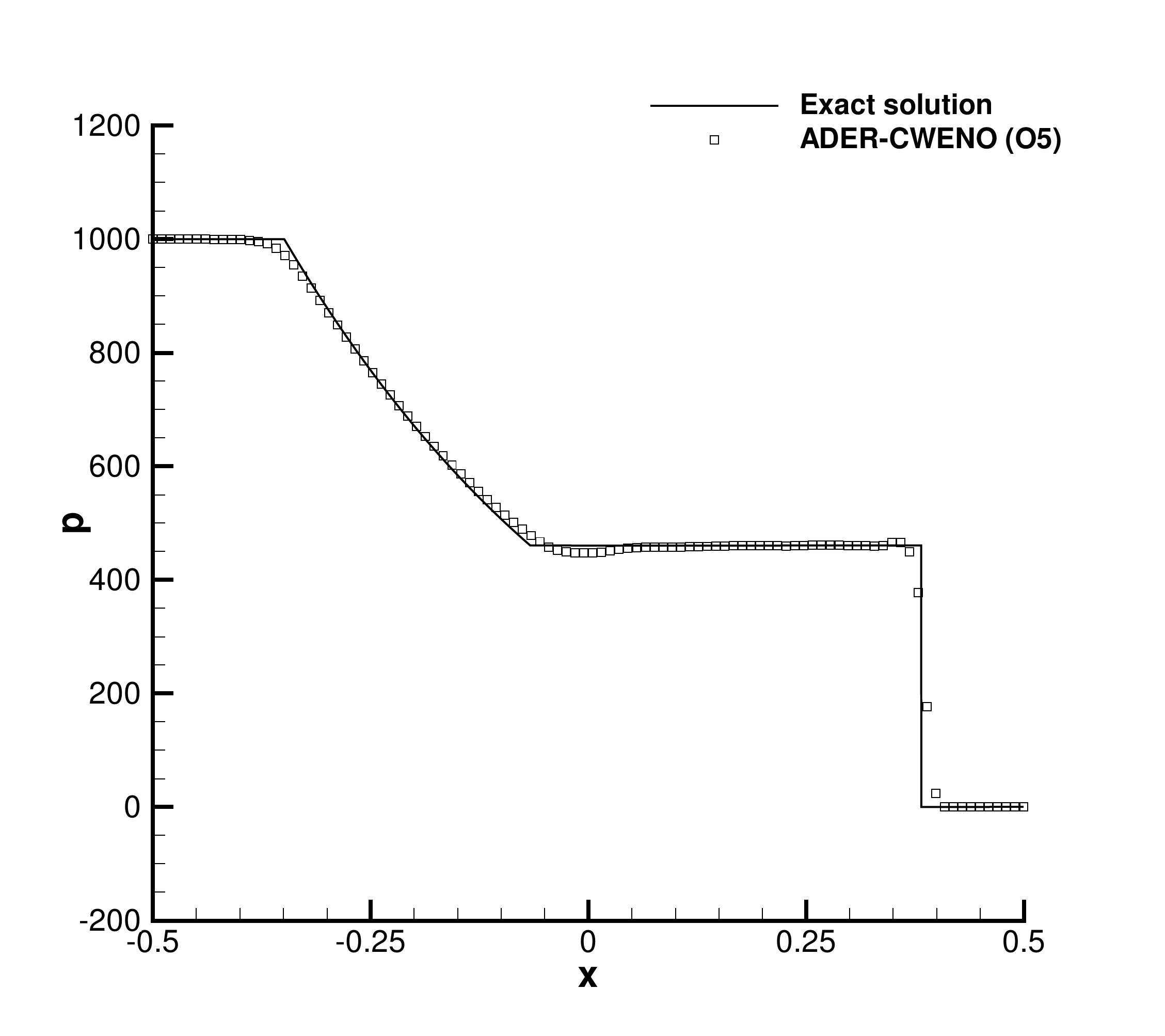} \\
\end{tabular} 
\caption{Numerical results for Riemann problem RP3 ($M=4$) at time $t=0.012$: 3D view of $\log(\rho)$ and comparison against the exact solution for density, velocity and pressure.}
\label{fig.RP3}
\end{center}
\end{figure}

\begin{figure}[!htbp]
\begin{center}
\begin{tabular}{cc} 
\includegraphics[width=0.44\textwidth]{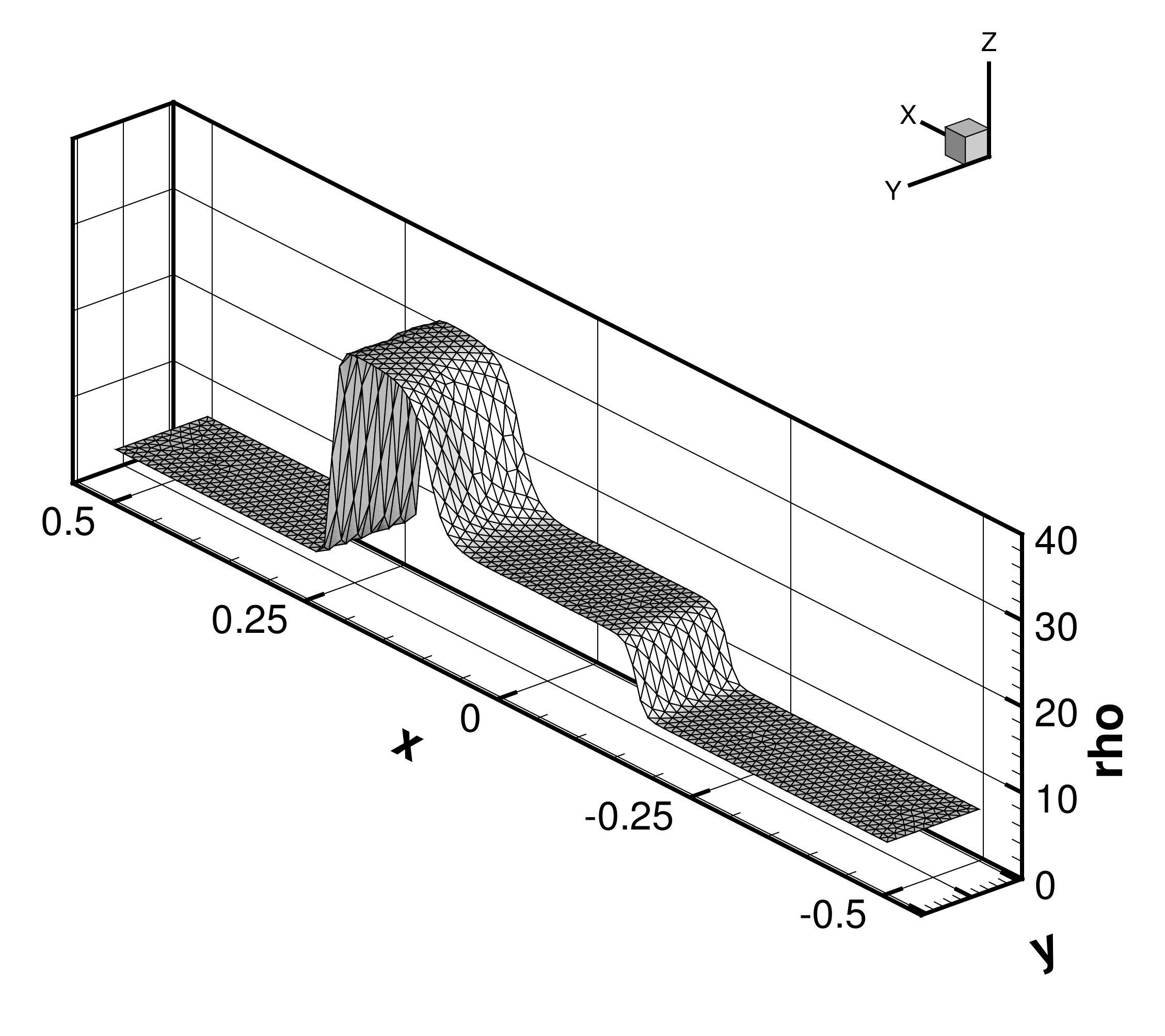}  &           
\includegraphics[width=0.44\textwidth]{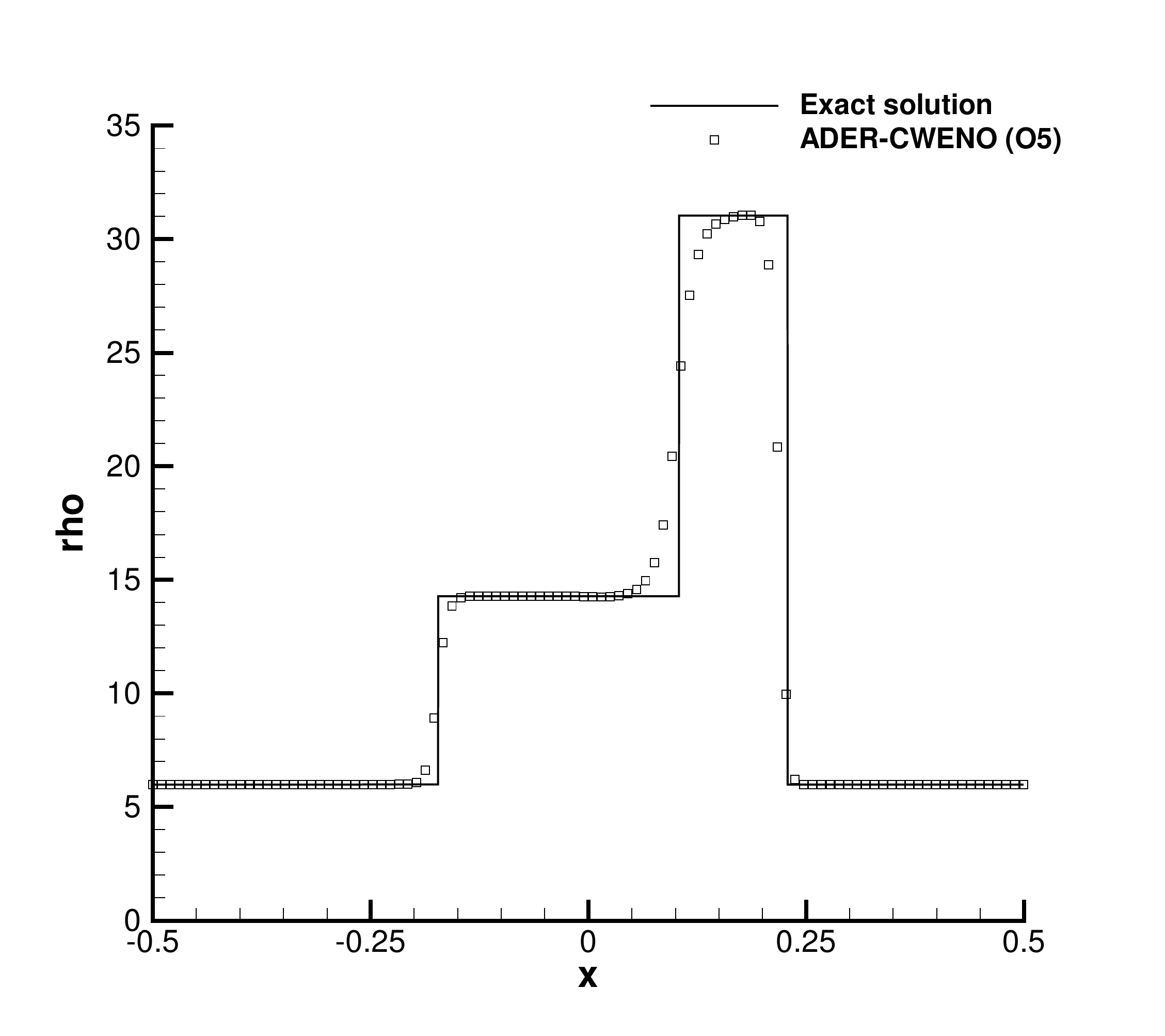} \\
\includegraphics[width=0.44\textwidth]{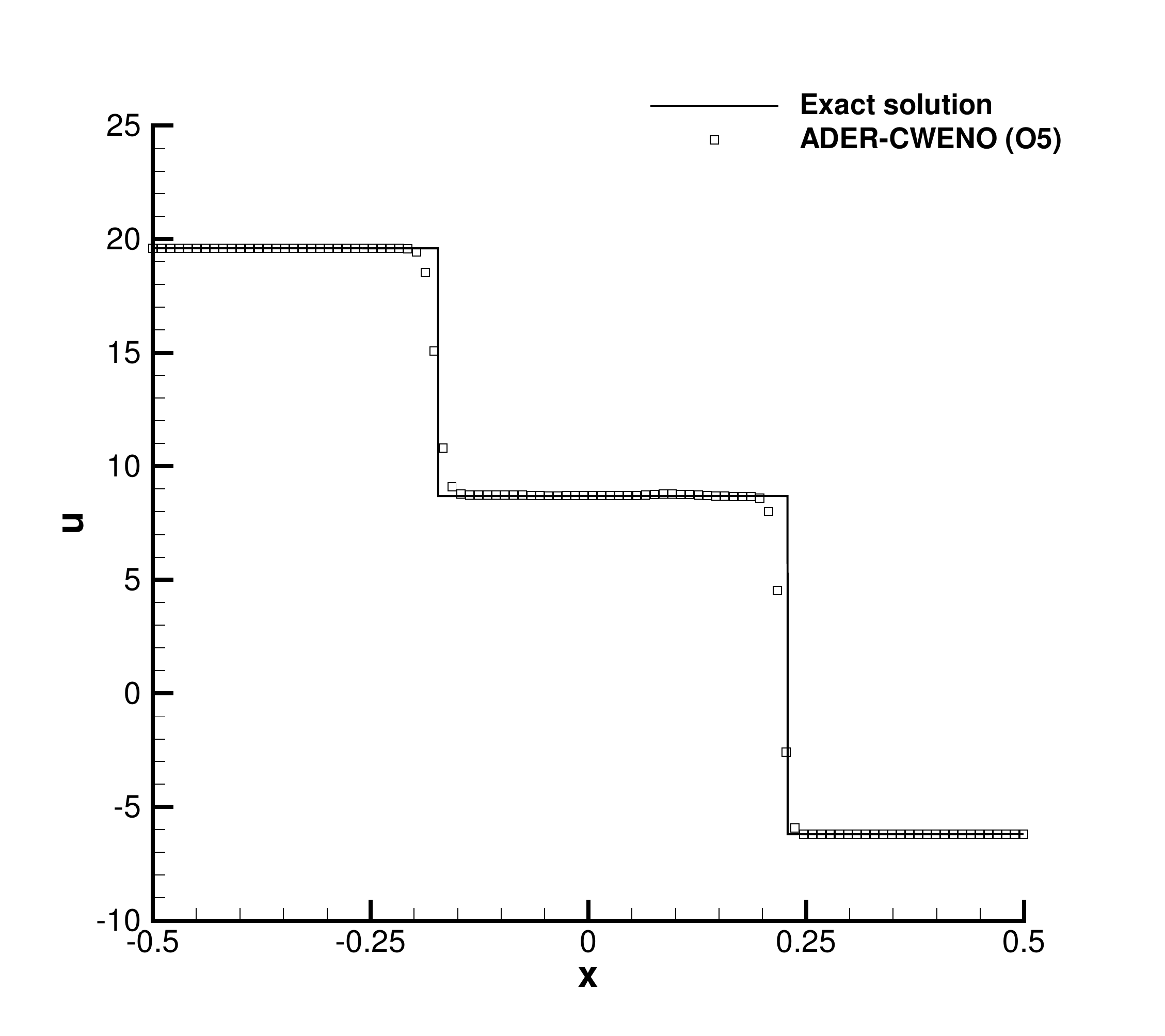}  &           
\includegraphics[width=0.44\textwidth]{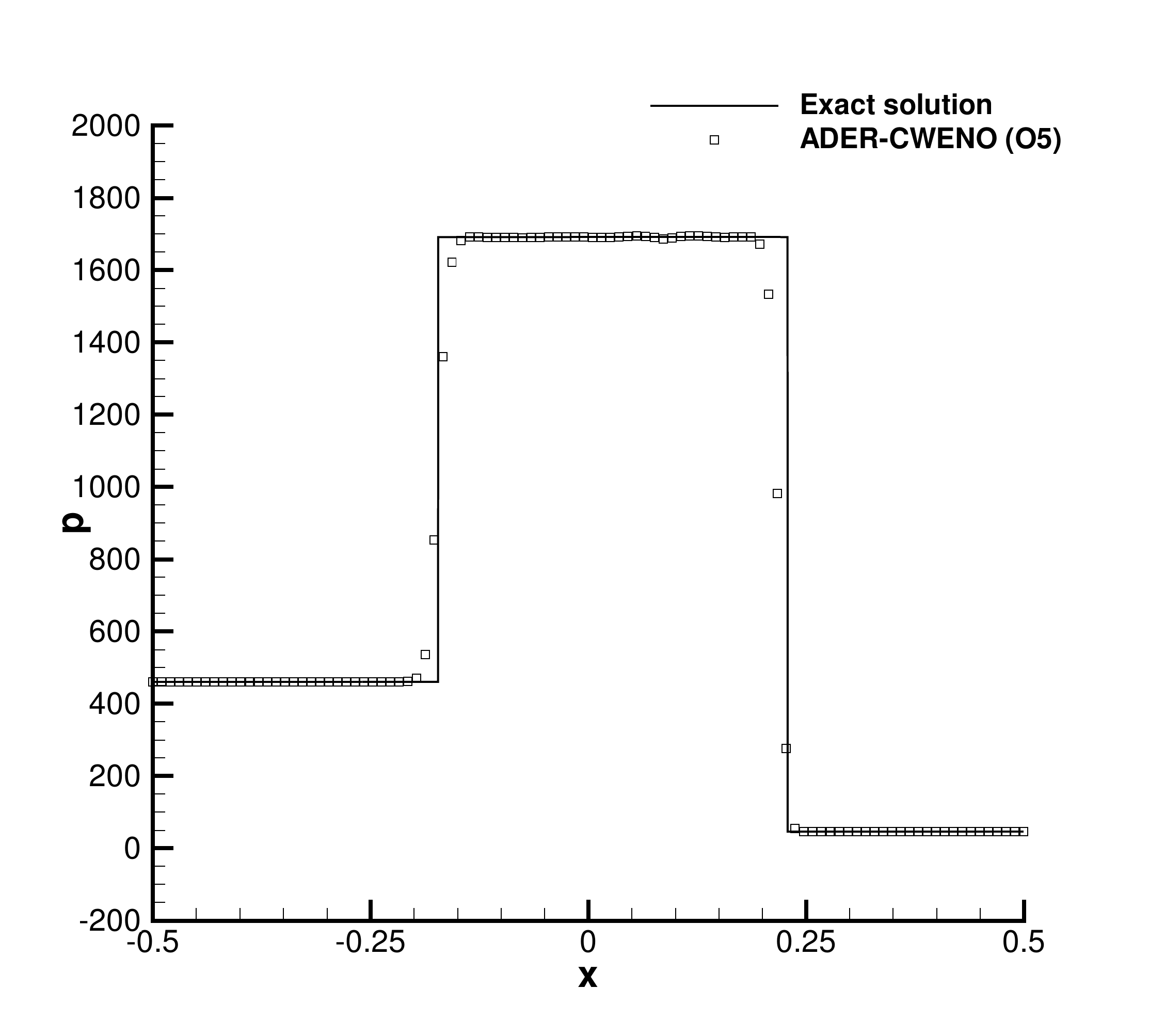} \\
\end{tabular} 
\caption{Numerical results for  Riemann problem RP4 ($M=4$) at time $t=0.035$: 3D view of the density distribution and comparison against the exact solution for density, velocity and pressure.}
\label{fig.RP5}
\end{center}
\end{figure}

Here we propose to apply the ADER-CWENO finite volume schemes for the solution of two well known one-dimensional Riemann problems, namely the Sod and the Lax shock tube test cases. The computational domain is the rectangular box $\Omega=[-0.5;0.5] \times [-0.05; 0.05]$ and the initial condition for the Euler equations \eqref{eulerTerms} is given in terms of a left state $\mathbf{U}_L$ and a right
state $\mathbf{U}_R$, separated by a discontinuity located in $x_d$. The ratio of specific heats is set to $\gamma=1.4$ and the detailed values for the initial states as well as the final simulation 
times and the location of the initial discontinuity are reported in Table \ref{tab.RP-IC}. The first two Riemann problems are the well known problems of Sod and Lax, respectively, while the difficult 
problems RP3 and RP4 are taken from the textbook of Toro \cite{ToroBook}. 
The numerical solution obtained on a rather coarse mesh composed of 2226 triangles with characteristic mesh size $h=1/100$ is depicted in Figures \ref{fig.Sod} - \ref{fig.RP5} 
and is compared against the exact solution of the Riemann problem of the Euler equations of compressible gas dynamics, see \cite{ToroBook} for details. 
We use $M=4$ and the Osher-type \eqref{eqn.osher} numerical flux for running the first two test cases, while the simple Rusanov flux is used for RP3 and RP4. Overall, a good agreement with the 
exact solution can be appreciated in all cases. The 1D plots in Figures \ref{fig.Sod} - \ref{fig.RP5} have been obtained from a one-dimensional cut through the reconstructed numerical 
solution $\w_i$ along the $x$-axis, evaluated at the final time on 100 equidistant sample points. 

\begin{table}[!t]
 \caption{Initial states left (L) and right (R) for the density $\rho$, velocity component $u$ and the pressure $p$ for the Riemann problems of the compressible Euler equations. 
          The final simulation times $t_{\textnormal{end}}$ and the initial position of the discontinuity $x_d$ are also given.} 
\begin{center} 
 \begin{tabular}{ccccccccc}
 \hline
 Case & $\rho_L$ & $u_L$ & $p_L$ & $\rho_R$ & $u_R$ & $p_R$ & $t_{\textnormal{end}}$ & $x_d$   \\ 
 \hline
 Sod    &  1.0      &  0.0       & 1.0     & 0.125      &  0.0        & 0.1      & 0.2    &  0.0    \\
 Lax    &  0.445    &  0.698     & 3.528   & 0.5        &  0.0        & 0.571    & 0.14   &  0.0    \\
 RP3    &  1.0      &  0.0       & 1000    & 1.0        &  0.0        & 0.01     & 0.012  &  0.1    \\
 RP4    &  5.99924  &  19.5975   & 460.894 & 5.99242    & -6.19633    & 46.095   & 0.035  & -0.2    \\ 
 \hline
 \end{tabular}
\end{center} 
 \label{tab.RP-IC}
\end{table}


\subsubsection{Explosion problems}
\label{ssec.EP}
The explosion problems can be seen as a multidimensional extension of the Sod test case presented in the previous section. The computational domain is given by the unit circle of radius $R=1$ in 
2D. In 3D we use the half-sphere of radius $R=1$ for $x \geq 0$, i.e. $\Omega_{3D}=\left\{ \x : x \geq 0 \, \wedge \, \left\| \x \right\| \leq R \right\}$. The initial condition is composed  of two different states reported in Table \ref{tab.RP-IC}, 
separated by a discontinuity at radius $R_d=0.5$. The ratio of specific heats is $\gamma = 1.4$ and the final time is $t_f=0.2$, so that the  shock  wave does not cross the external boundary 
of the domain, where a transmissive boundary condition is set. We run this problem in 2D and 3D in two different configurations: the first case uses a rather high order reconstruction, in order 
to show that the CWENO method can be at least in principle implemented for any order of accuracy, maintaining its ability to avoid spurious oscillations in the vicinity of shock waves also for 
high order reconstructions on unstructured meshes. The second case uses third order schemes ($M=2$), but on very fine meshes involving piecwise polynomial reconstructions with \textbf{billions} 
of degrees of freedom, in order to show that the ADER-CWENO method is also very well suited for the implementation on massively parallel distributed memory supercomputers. 

\begin{figure}[!htbp]
\begin{center}
\begin{tabular}{cc} 
\includegraphics[width=0.44\textwidth]{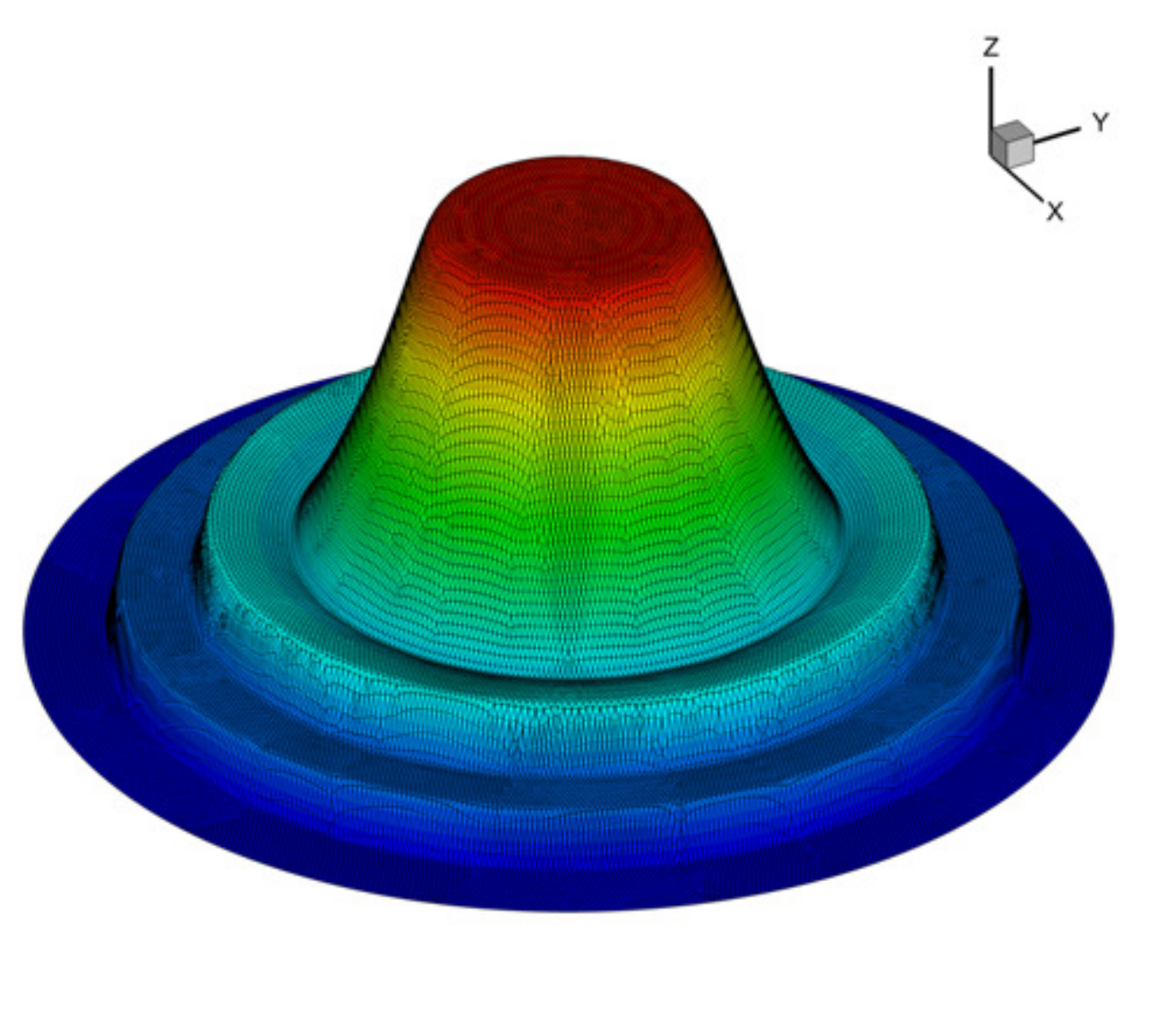}  &           
\includegraphics[width=0.44\textwidth]{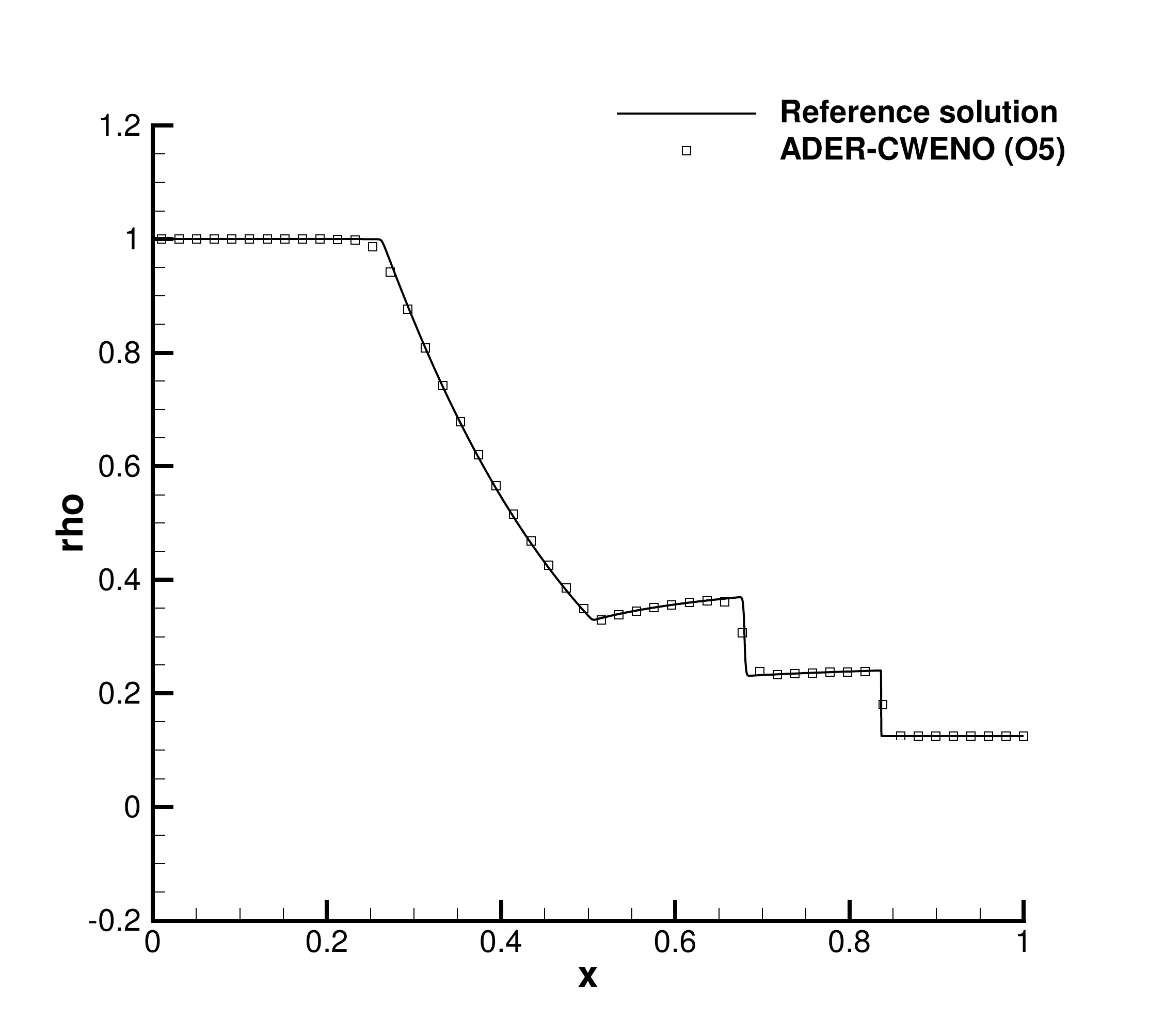} \\
\includegraphics[width=0.44\textwidth]{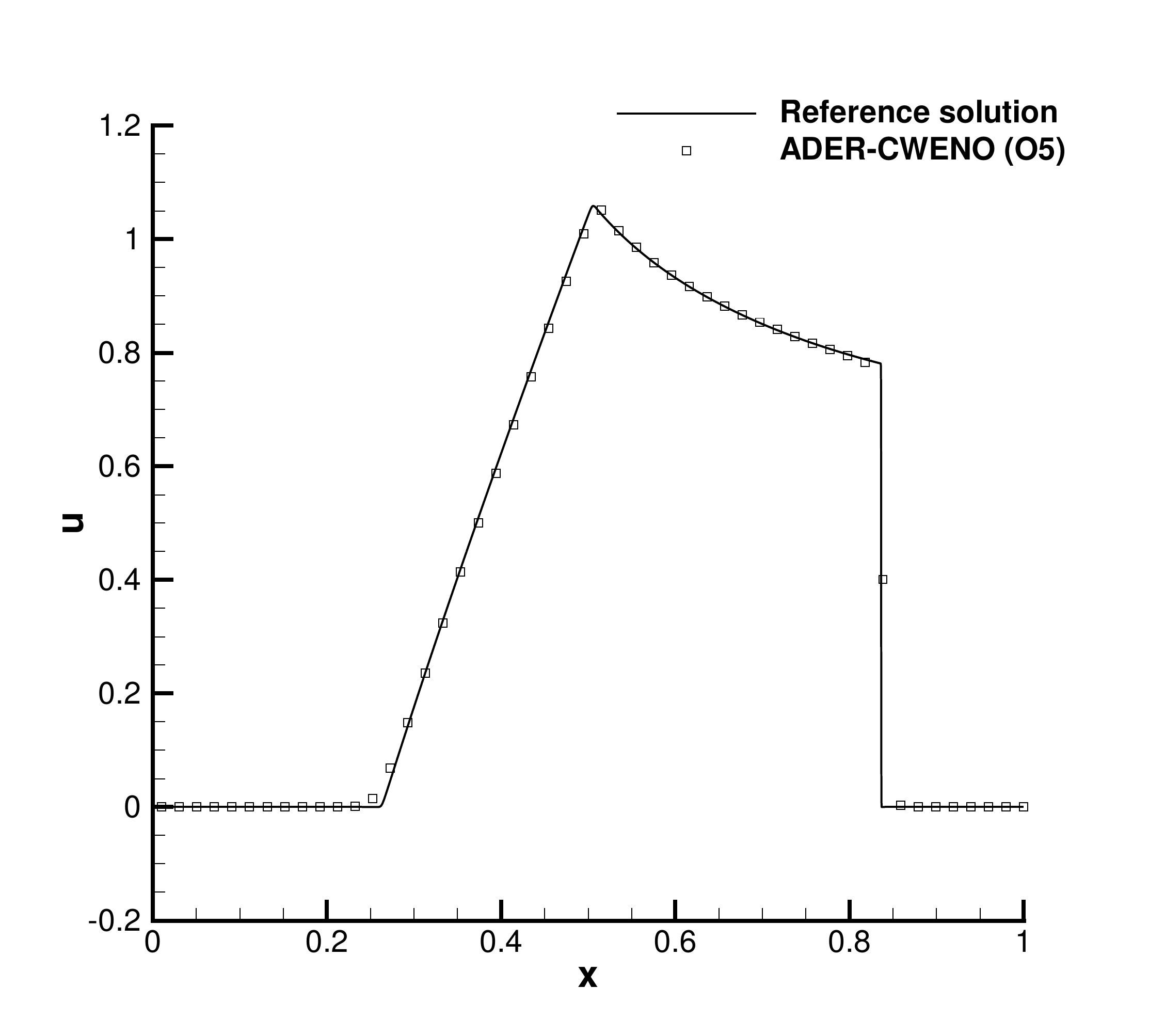}  &           
\includegraphics[width=0.44\textwidth]{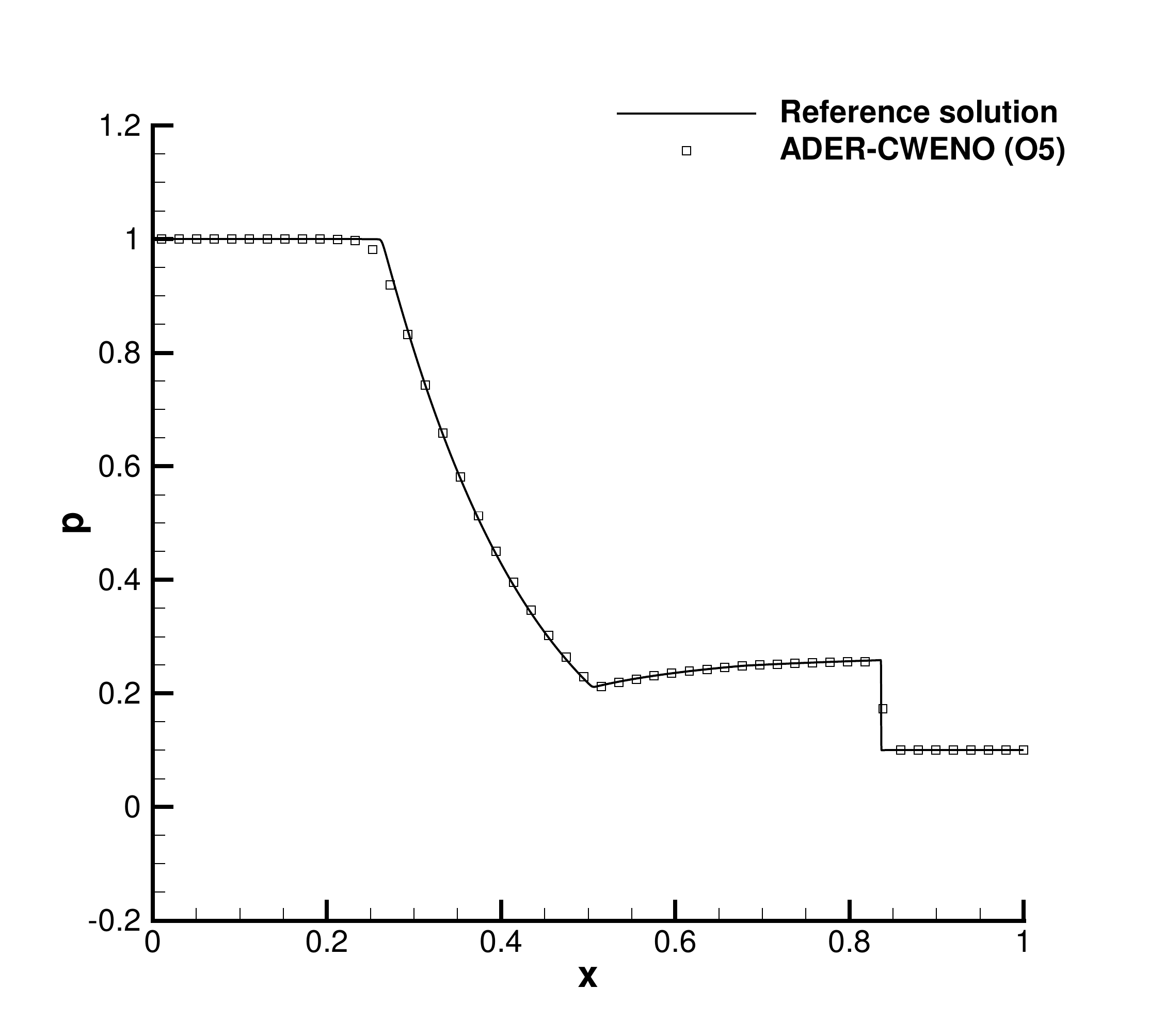} \\
\end{tabular} 
\caption{Numerical results for the coarse grid simulation of the two-dimensional explosion problem ($M=4$, 1,024,860 degrees of freedom) at time $t=0.2$. From top left to bottom right: 
3D view of the density distribution together with the unstructured triangular mesh and comparison against the reference solution for density, radial velocity and pressure.} 
\label{fig.EP2D}
\end{center}
\end{figure}

\begin{figure}[!htbp]
\begin{center}
\begin{tabular}{cc} 
\includegraphics[width=0.44\textwidth]{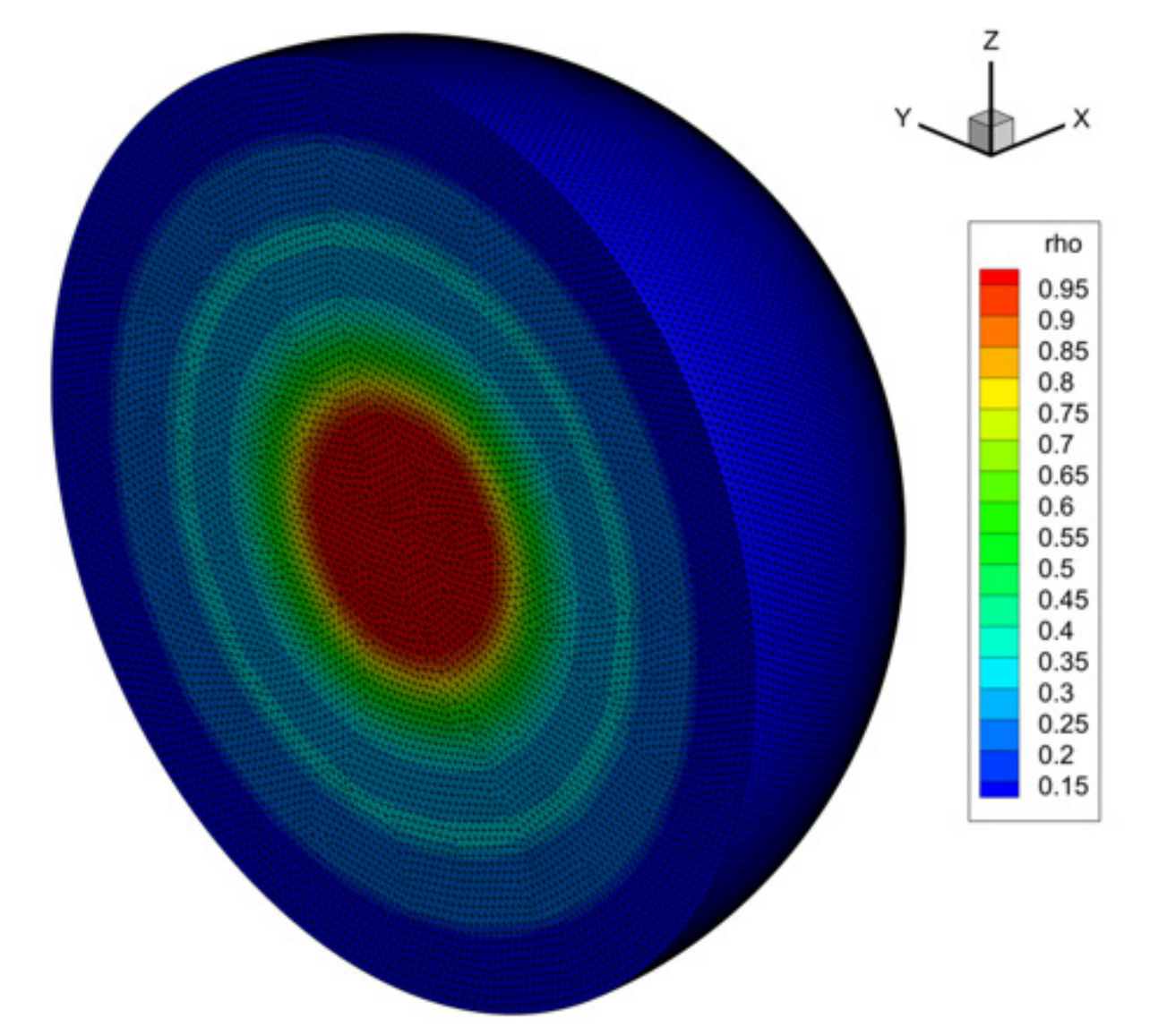}  &           
\includegraphics[width=0.44\textwidth]{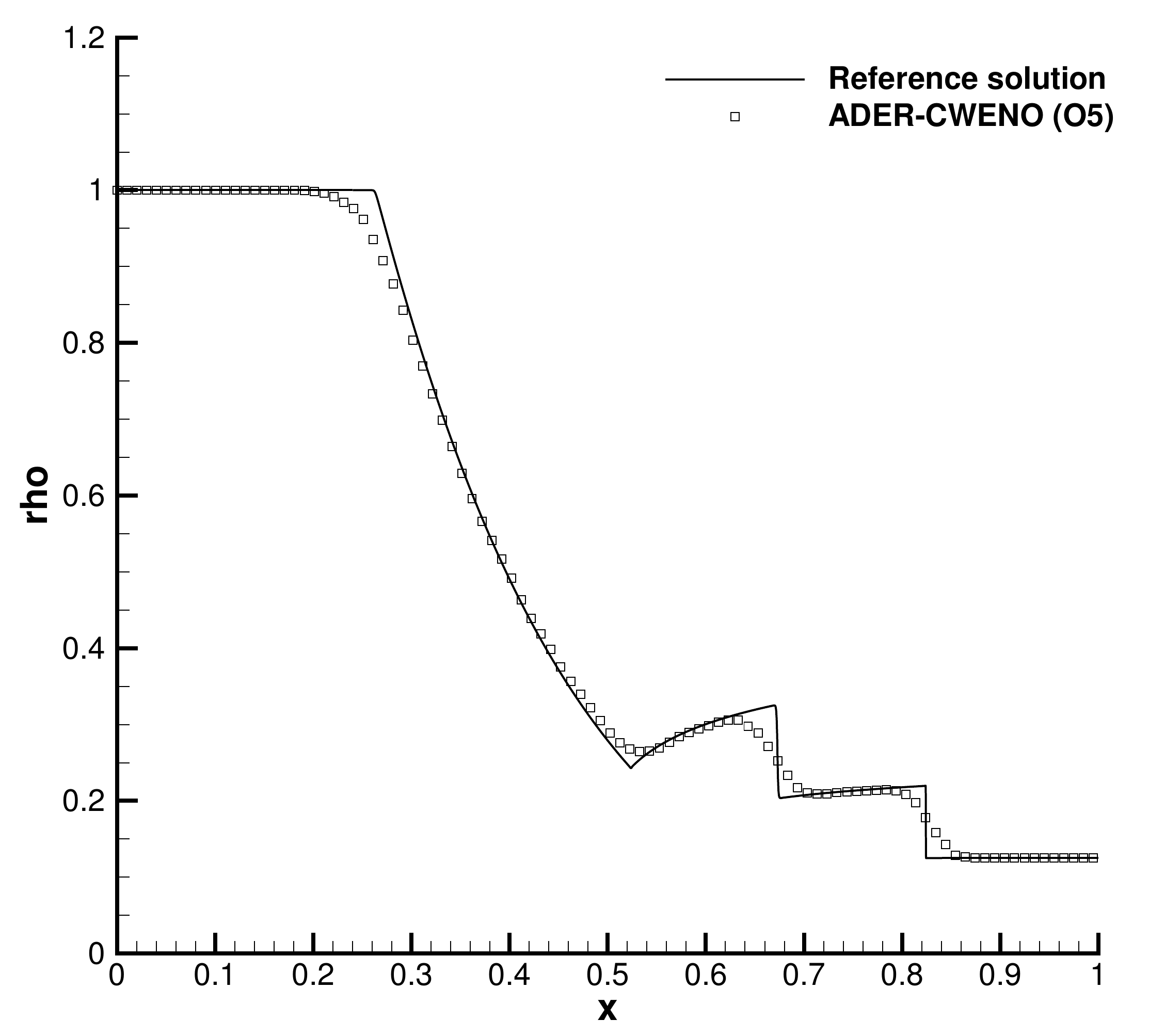} \\
\includegraphics[width=0.44\textwidth]{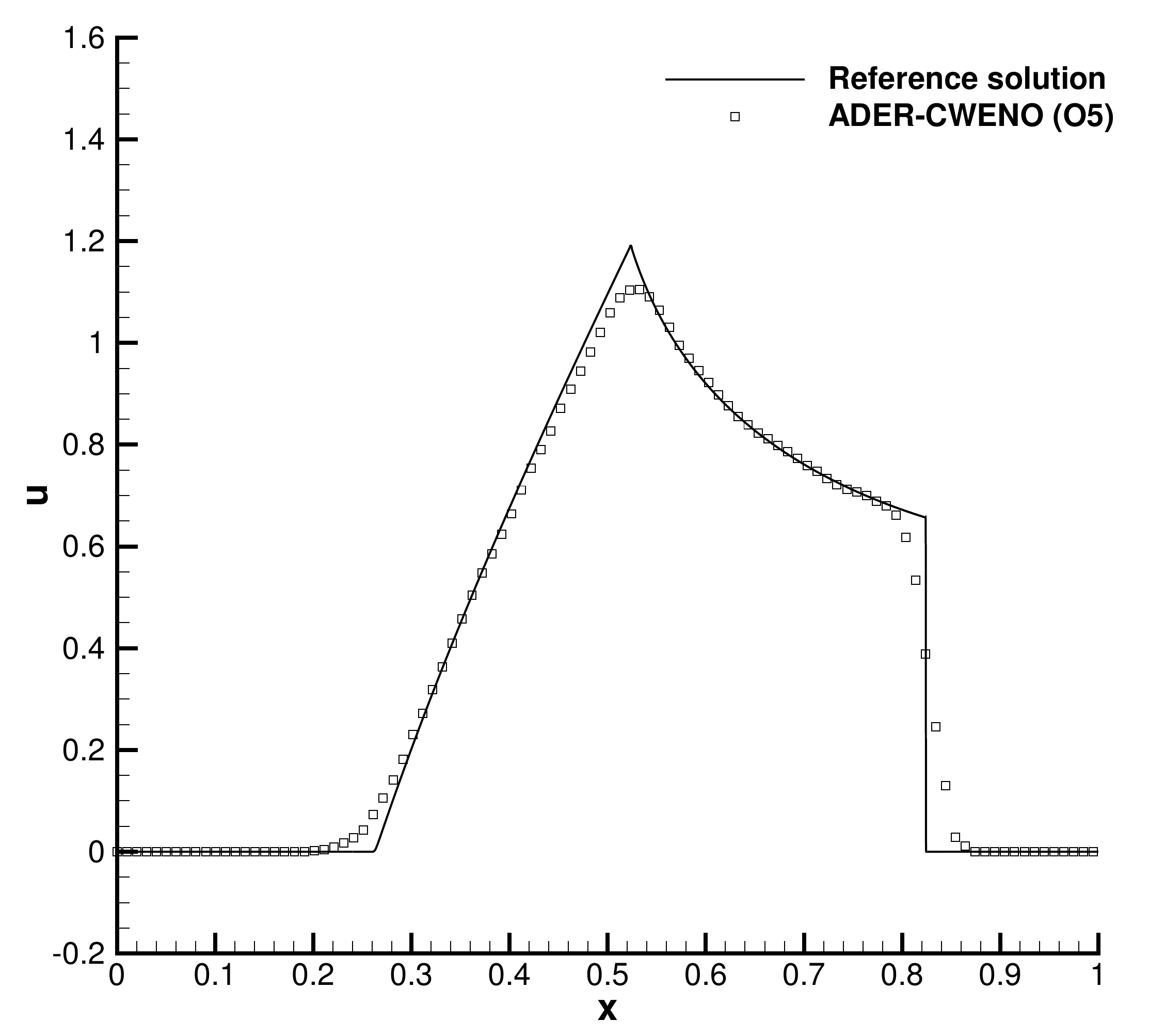}  &           
\includegraphics[width=0.44\textwidth]{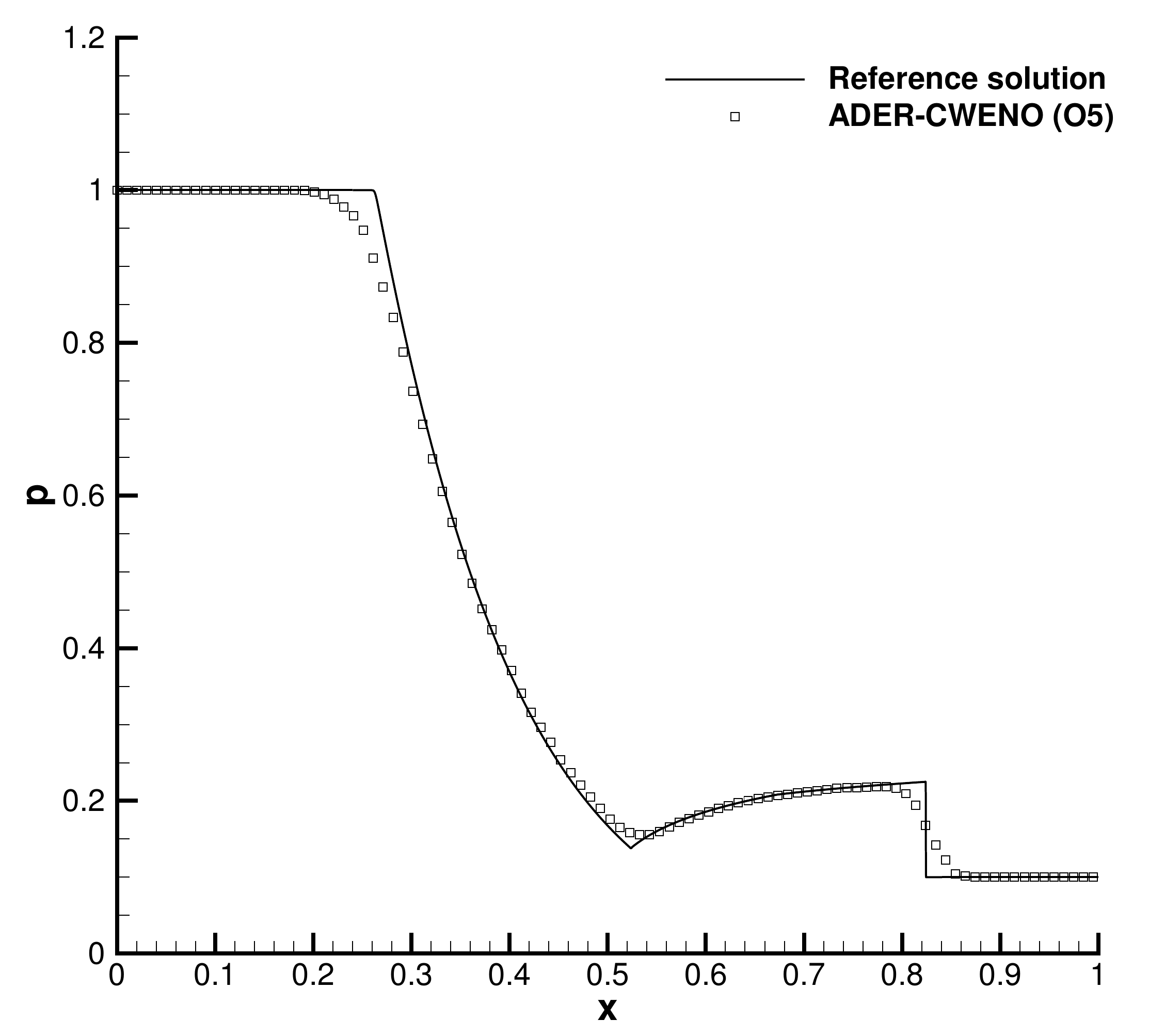} \\
\end{tabular} 
\caption{Numerical results for the coarse grid simulation of the three-dimensional explosion problem ($M=4$, 64,643,810 degrees of freedom) at time $t=0.2$. 
From top left to bottom right: density contours and computational mesh, comparison against the reference solution for density, horizontal velocity and pressure.}  
\label{fig.EP3D}
\end{center}
\end{figure}

\paragraph{Case I, high order reconstructions} In the first run, the numerical results have been obtained employing a fifth order CWENO reconstruction ($M=4$), together with the  
Osher-type flux \eqref{eqn.osher} using rather coarse meshes that discretize the computational domain with a total number of $N_E=68,324$ triangles in 2D ($h=1/100$) and $N_E=1,846,966$ 
tetrahedra in 3D ($h=1/50$), respectively. The results of this first run are depicted in Figures \ref{fig.EP2D} and \ref{fig.EP3D}, where also a comparison with the reference solution is 
given. 
We can observe a good agreement between the numerical results obtained with the high order ADER-CWENO schemes and the reference solution. 
As in \cite{Lagrange3D,ToroBook} the reference solution can be obtained by making use of the spherical symmetry of the problem and by solving a reduced one-dimensional system 
with geometric source terms using a classical second order TVD scheme on a very fine one-dimensional mesh. This test problem involves three different waves, namely one cylindrical or spherical 
shock wave that is running towards the external boundary of the domain, a rarefaction fan traveling in the opposite direction and an outward-moving contact wave in between. 

\paragraph{Case II, very fine mesh with more than one billion degrees of freedom} 

\begin{figure}[!htbp]
\begin{center}
\begin{tabular}{cc} 
\includegraphics[width=0.44\textwidth]{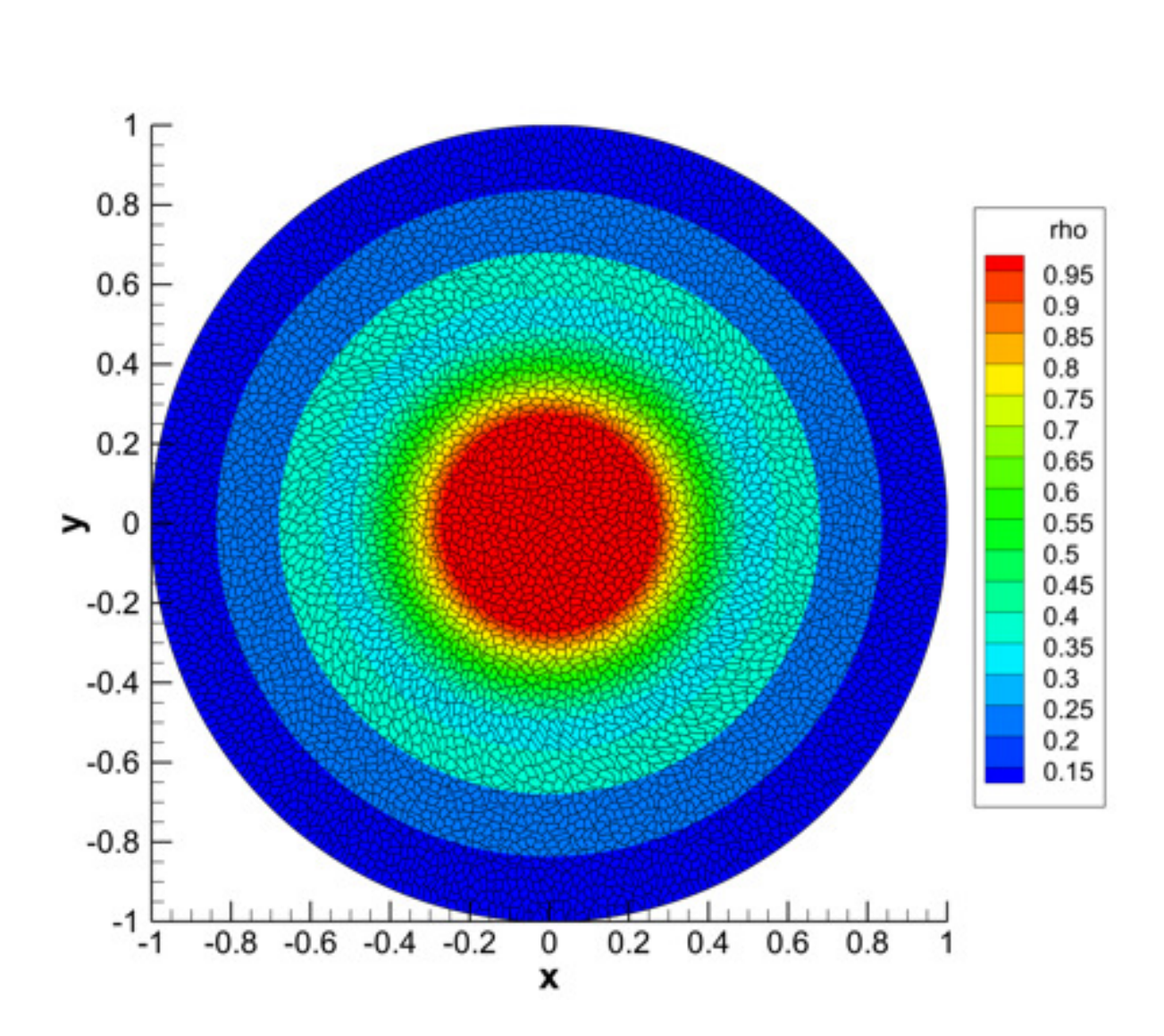}  &           
\includegraphics[width=0.44\textwidth]{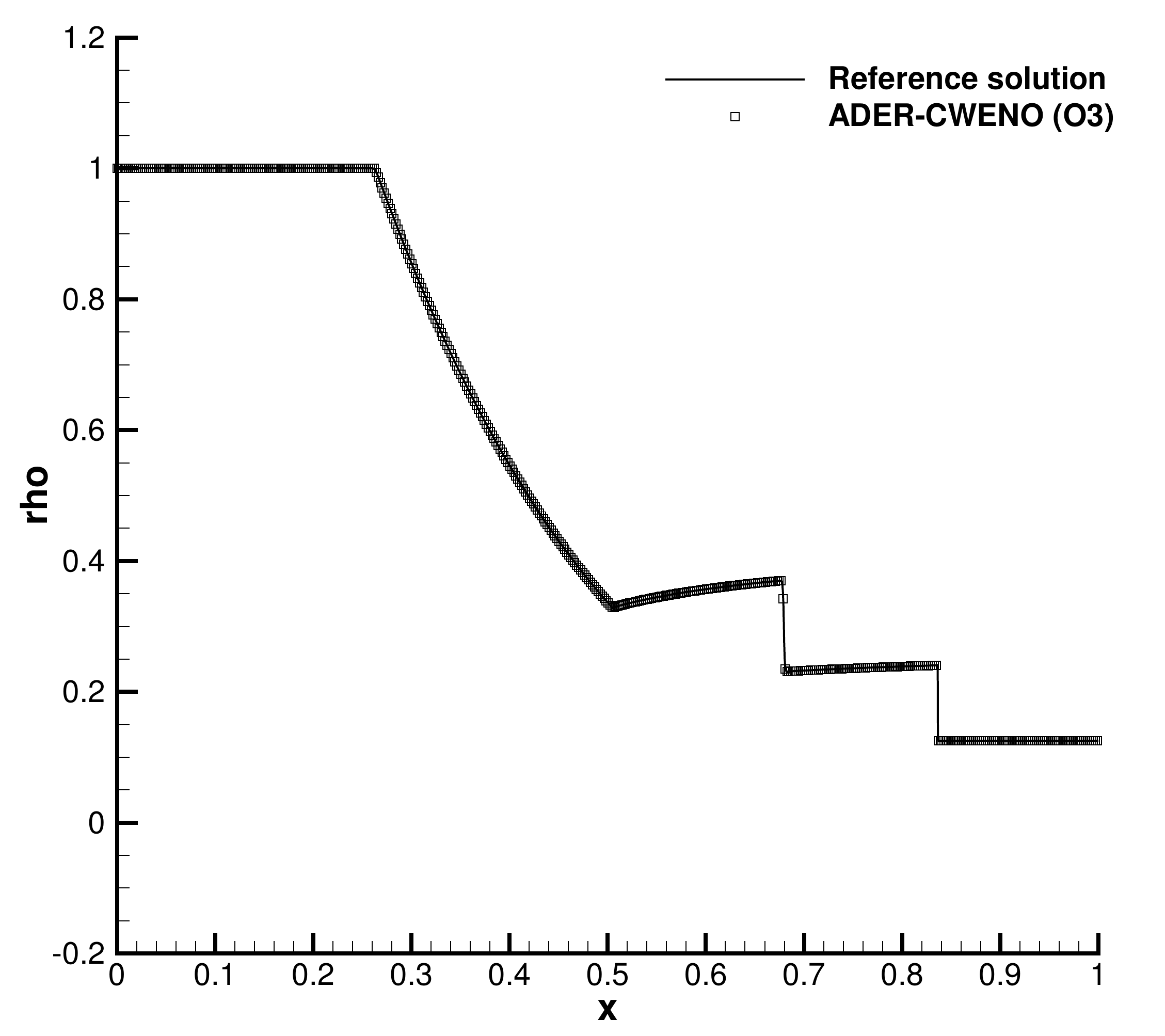} \\
\includegraphics[width=0.44\textwidth]{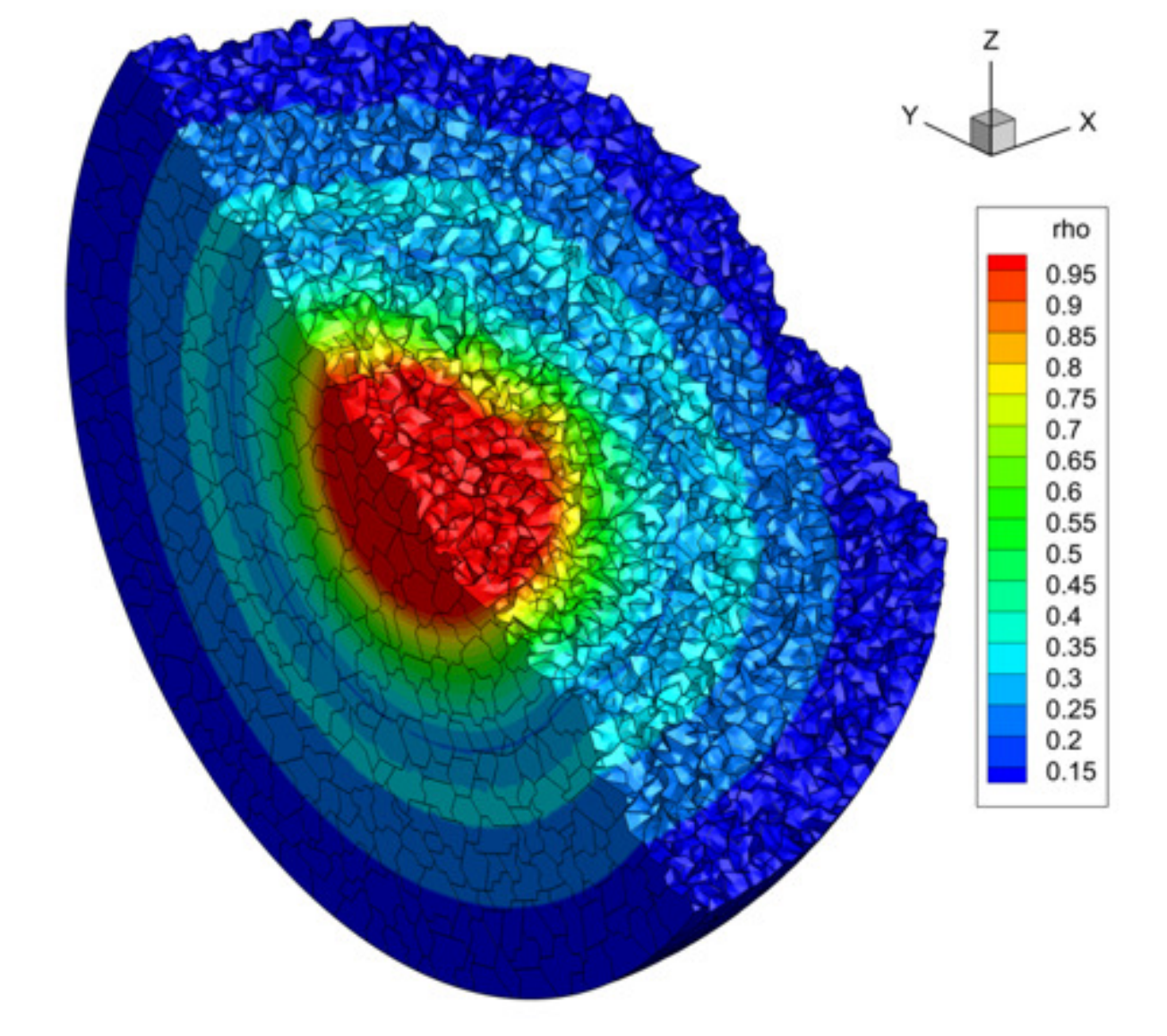}  &            
\includegraphics[width=0.44\textwidth]{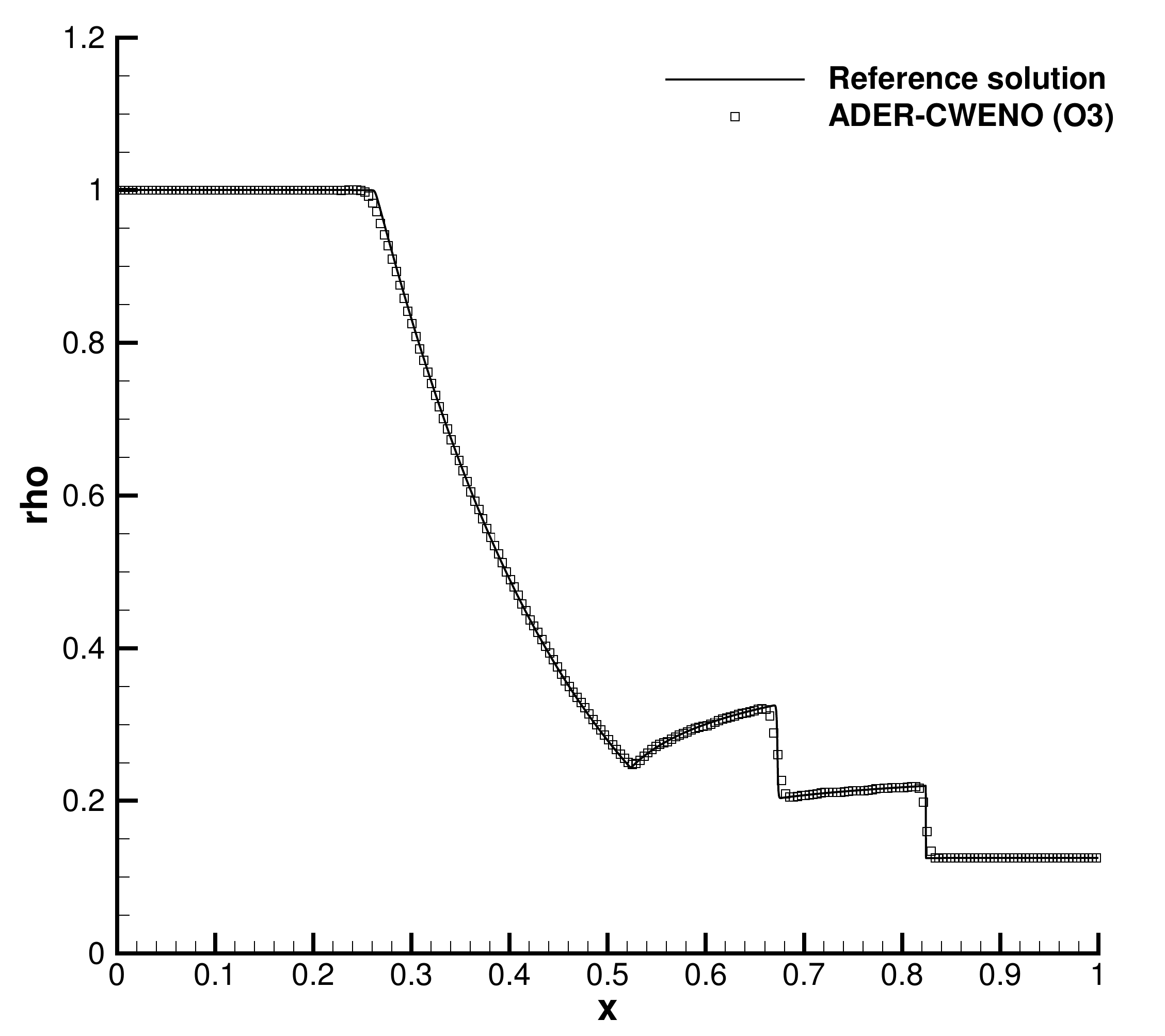} \\
\end{tabular} 
\caption{Numerical results obtained for the fine grid simulation of the explosion problem at time $t=0.2$. Top row: 2D simulation with 
$M=2$ and 1,038,622,800 degrees of freedom. 
Bottom row: 3D simulation with $M=2$ and 2,307,645,000 degrees of freedom. Density contours and domain decomposition onto the various MPI ranks (left). 
1D cut along the $x$-axis and comparison against the reference solution for the density (right).} 
\label{fig.EPfine}
\end{center}
\end{figure}

We now run the 2D and the 3D explosion problems again, but this time using two very fine unstructured grids. The fine 2D mesh is composed of a total number of 173,103,800 triangular elements of 
characteristic mesh spacing $h=1/5000$, while the 3D mesh consists of a total number of 230,764,500 tetrahedral elements with characteristic mesh spacing $h=1/250$. 
The two-dimensional simulation is run on 5600 CPUs and the three-dimensional one is run on 7168 CPUs at the SuperMUC supercomputer of the Leibniz Rechenzentrum (LRZ) in Munich, Germany. 
In both cases a third order CWENO reconstruction ($M=2$) is employed, hence leading to 6 and 10 degrees of freedom per element in 2D and 3D, respectively. Thus, the total number of spatial 
degrees of freedom for the representation of the piecewise polynomial CWENO reconstruction is 1,038,622,800 in 2D and 2,307,645,000 in 3D, respectively. 
A sketch of the domain decomposition onto the MPI ranks as well as a one-dimensional cut through the computational results along the $x$ axis are depicted for both cases in Figure \ref{fig.EPfine}. 
In both cases, an excellent agreement with the reference solution is obtained, as expected. 

\subsubsection{The two-dimensional double Mach reflection problem}
\label{ssec.DMR}
\begin{figure}[!htbp]
\begin{center}
\begin{tabular}{cc} 
\includegraphics[height=0.35\textwidth]{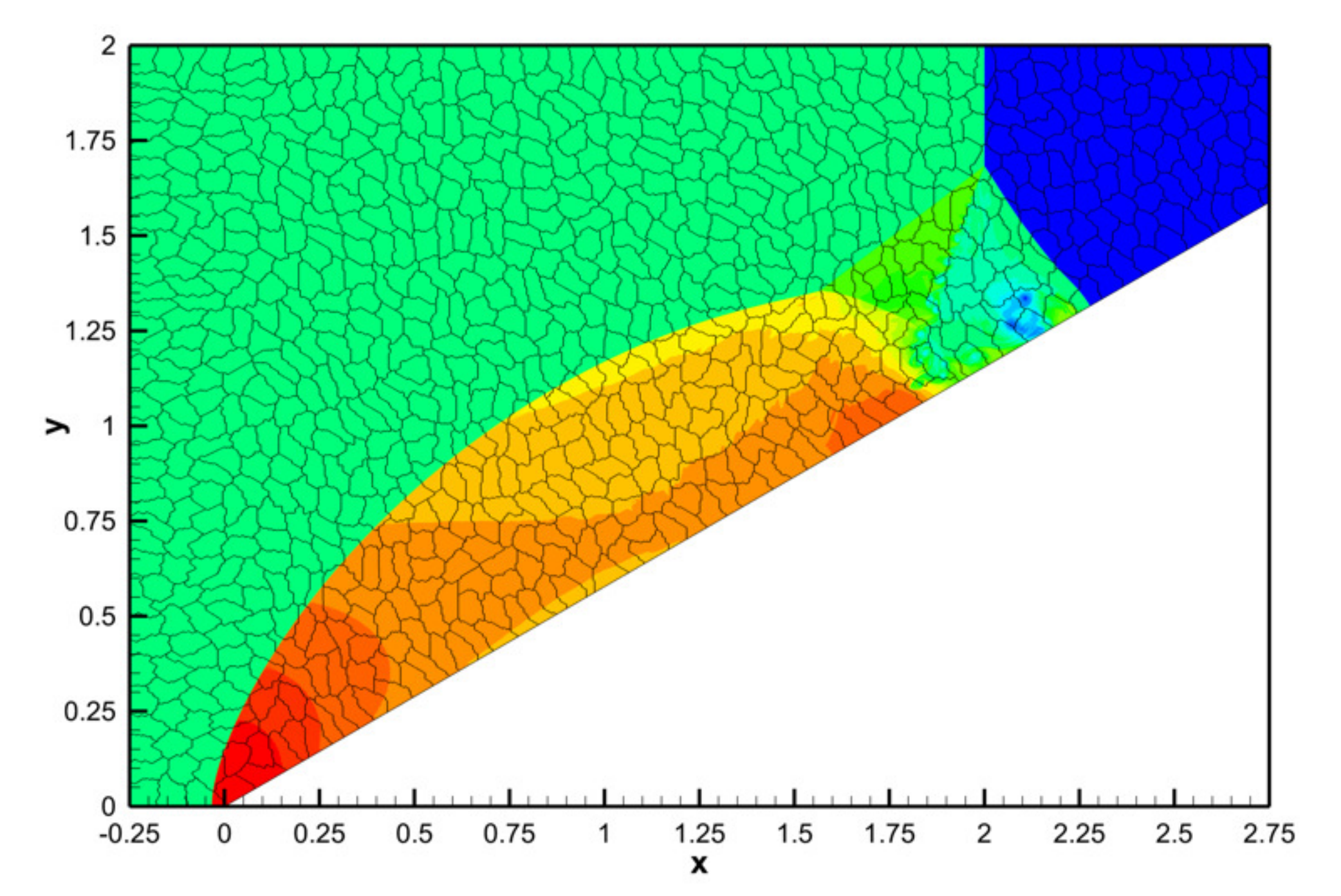}  &            
\includegraphics[height=0.35\textwidth]{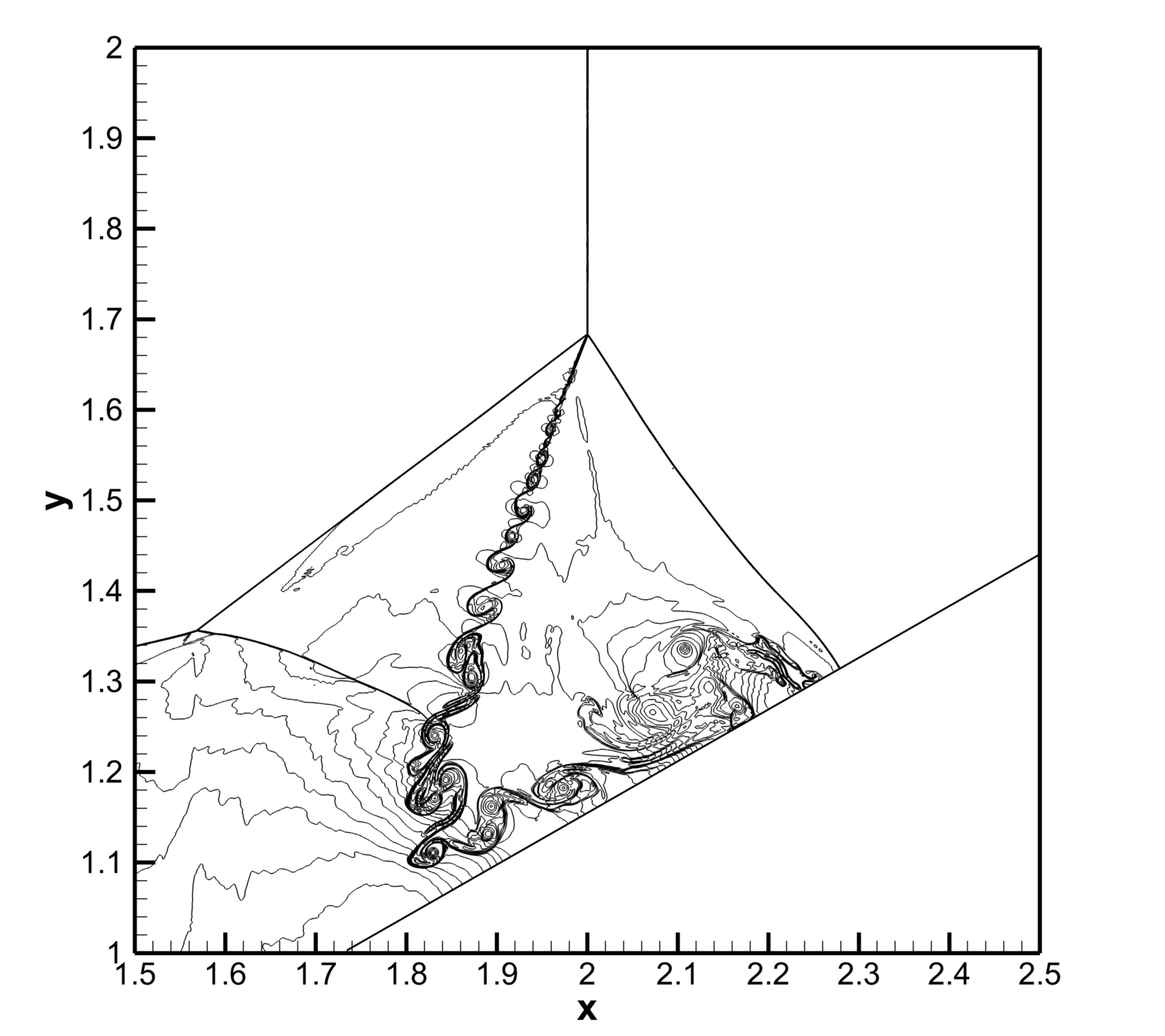} 
\end{tabular} 
\caption{Results for the double Mach reflection problem obtained with mesh spacing $h=1/2200$ ($M=2$, 260,645,616 degrees of freedom) 
at time $t=0.2$: density contours and domain decomposition onto 800 MPI ranks (left); zoom on the small-scale structures using 
41 equidistant density contour lines from 1.5 to 22.5 (right).}  
\label{fig.DMR}
\end{center}
\end{figure}

The double Mach reflection problem was originally proposed in \cite{woodwardcol84} and it considers a very strong shock wave that is moving along the $x-$direction of the computational domain, where a ramp with angle $\alpha=30^\circ$ is located. The shock Mach number is $M_s=10$ and small-scale structures are generated behind the shock wave that is impinging onto the ramp. The initial condition is  given by 
\begin{equation}
\U = \left\{ \begin{array}{lll} \left( 8.0, 8.25, 0, 116.5 \right), & \textnormal{ if } & x<x_0, \\ 
                                        \left( 1.4,0,0,1.0 \right), & \textnormal{ if } & x \geq x_0,  
                      \end{array}  \right. 
\label{eq:DMR_IC}
\end{equation}
where $x_0=0$ represents the initial location of the discontinuity. Since we use an unstructured mesh, the problem can be directly run in physical coordinates as in \cite{Dumbser2007204,SebastianShu},  without tilting the shock wave, as it is usually done for Cartesian codes. Slip wall boundary conditions are set on the upper and the lower side of the domain, while inflow and outflow boundaries are  imposed on the remaining sides. The final time is $t_f=0.2$ and the ratio of specific heats is $\gamma=1.4$. A fine unstructured mesh with characteristic mesh size of $h=1/2200$ composed of 
$N_E=43,440,936$ triangles is used together with a polynomial degree $M=2$ of the CWENO reconstruction and the Rusanov flux \eqref{eqn.rusanov}, leading to a total number of 260,645,616 degrees 
of freedom for the representation of the piecewise polynomial CWENO reconstruction. 
The computational results are depicted in Figure  
\ref{fig.DMR}. The small-scale structures produced by the roll-up of the shear layers behind the shock wave are clearly visible from the density contour lines in the right panel of Figure \ref{fig.DMR}, 
while the distribution of the computational domain onto 800 MPI ranks is highlighted in the left panel of Figure \ref{fig.DMR}. 
 
\subsubsection{The MHD rotor problem}
\label{ssec.MHDrotor}
\begin{figure}[!htbp]
\begin{center}
\begin{tabular}{cc} 
\includegraphics[width=0.44\textwidth]{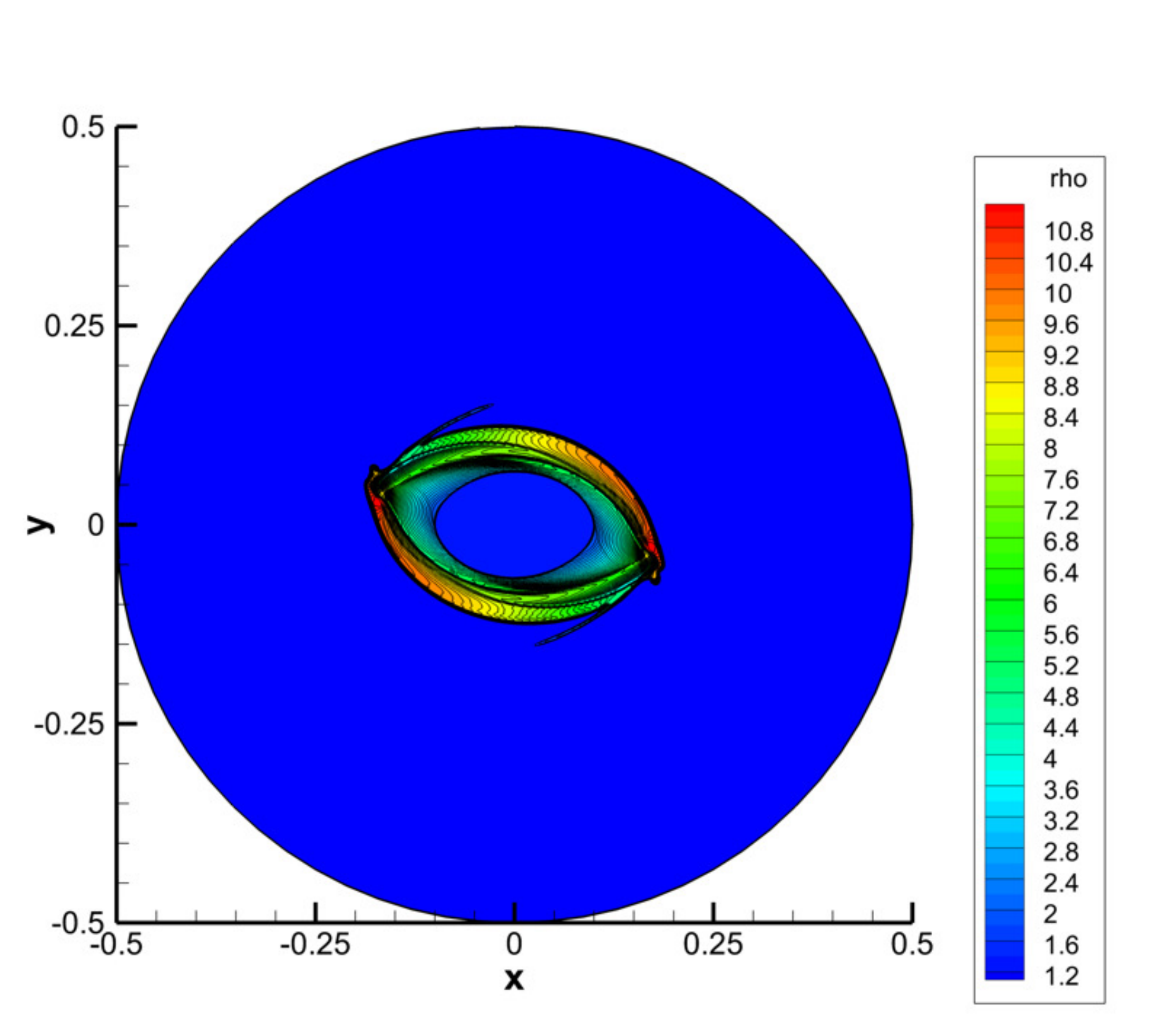}  &           
\includegraphics[width=0.44\textwidth]{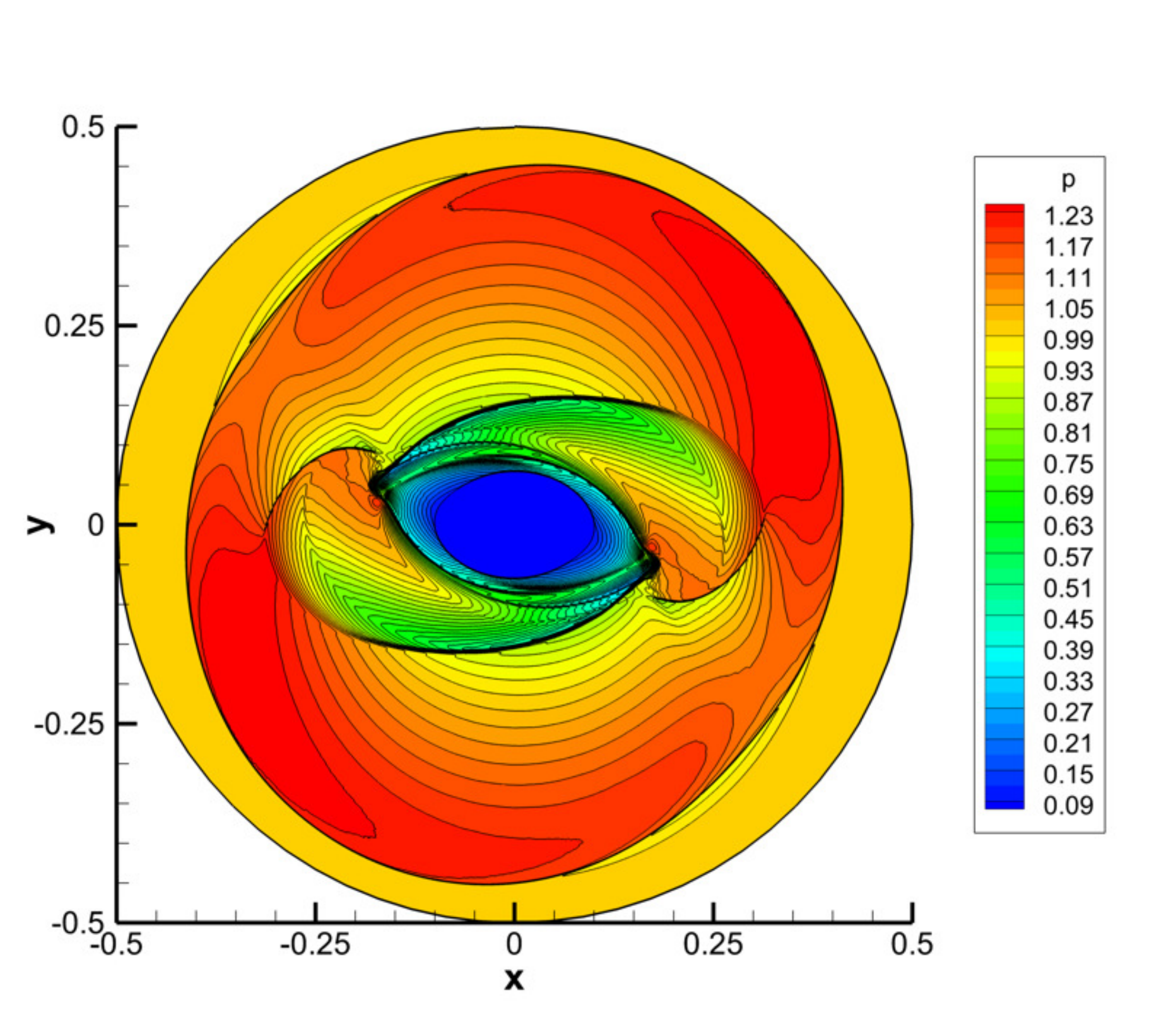} \\             
\includegraphics[width=0.44\textwidth]{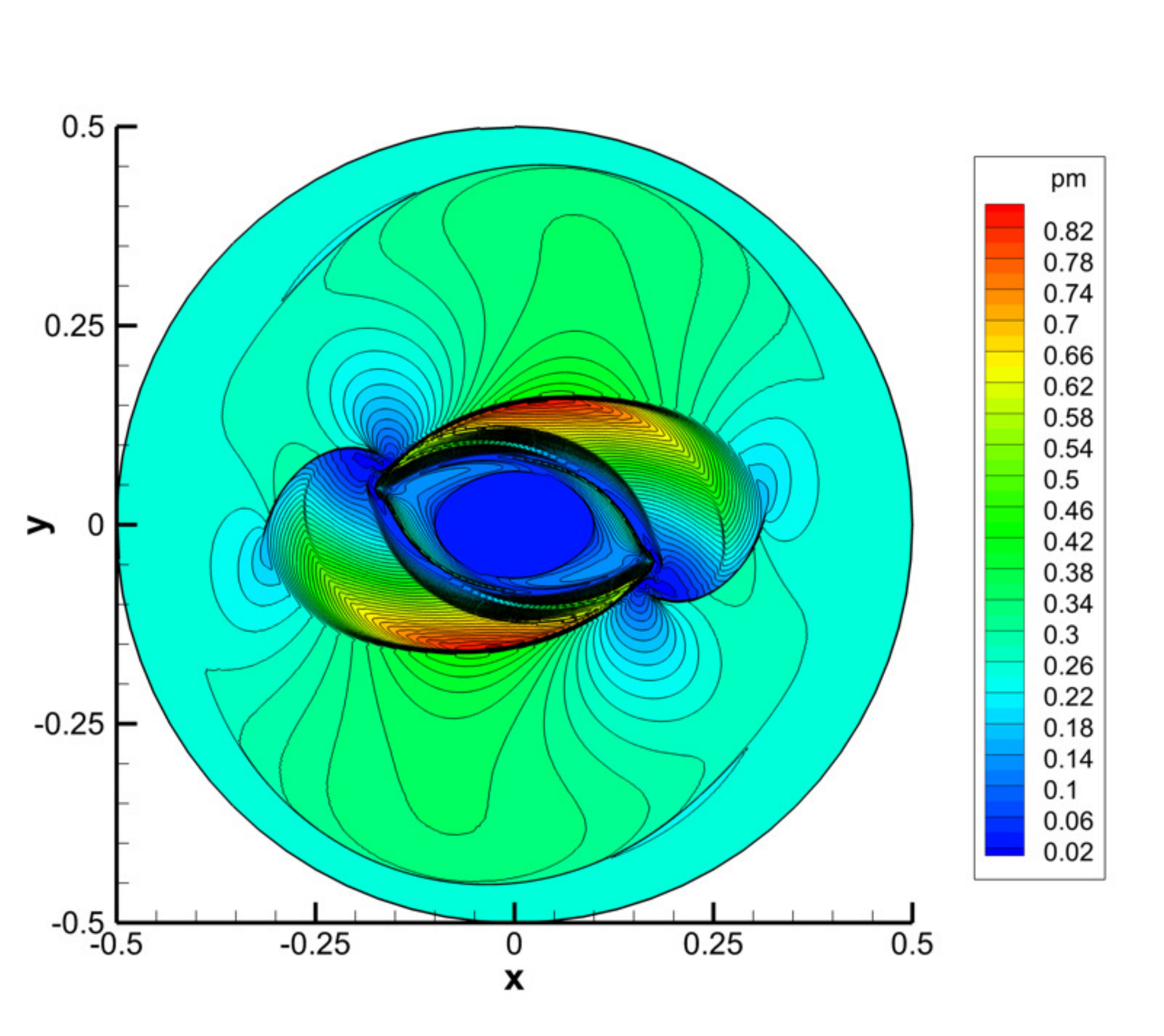}  &            
\includegraphics[width=0.44\textwidth]{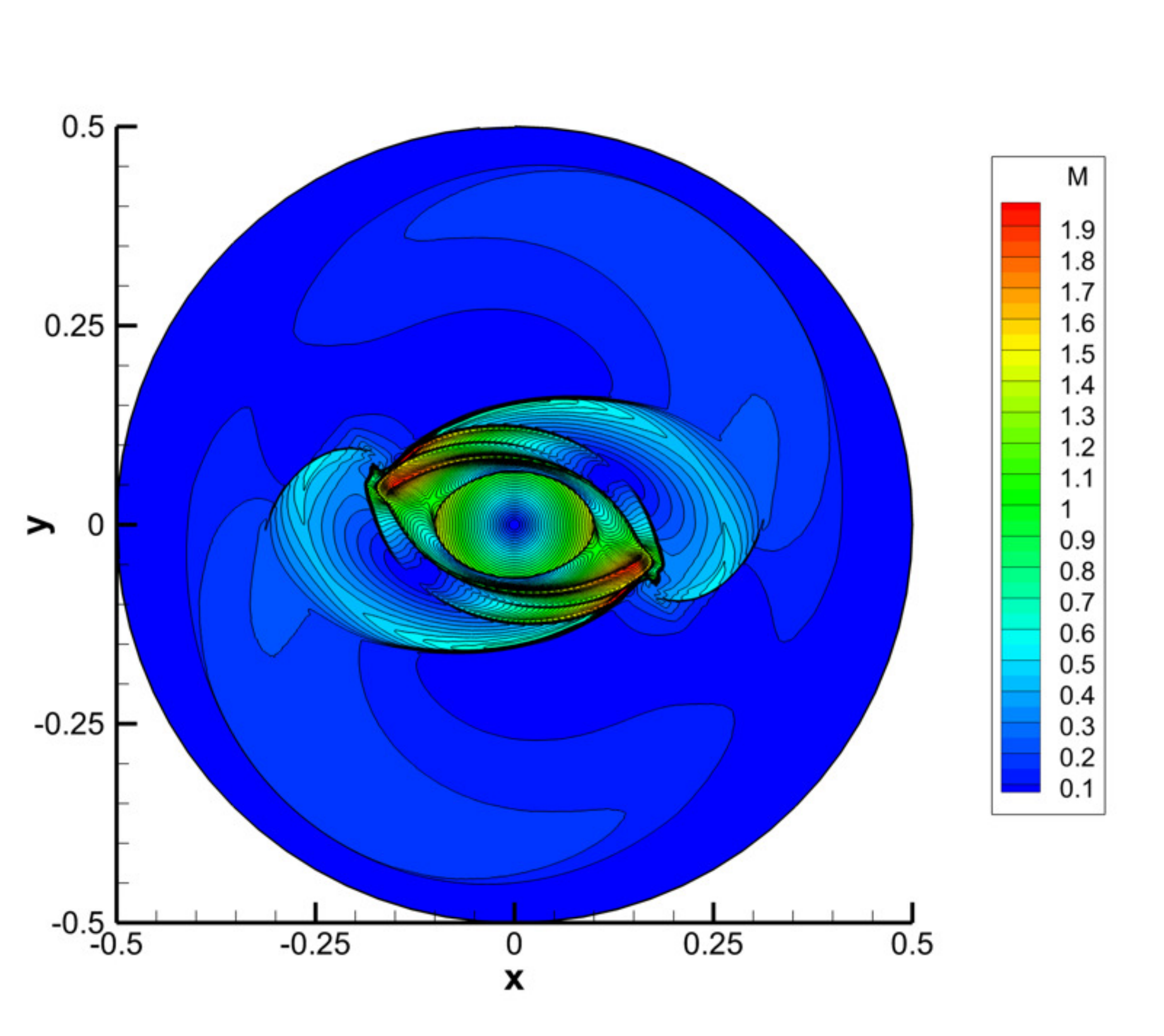}                
\end{tabular} 
\caption{Numerical results for the MHD rotor problem obtained with a third order unstructured ADER-CWENO scheme: density, pressure, magnetic pressure and Mach number at time $t=0.25$.} 
\label{fig.MHDRotor}
\end{center}
\end{figure}

The first test case for the ideal classical MHD equations \eqref{MHDTerms} is the MHD rotor problem \cite{BalsaraSpicer1999} that describes a high density fluid which is rotating inside 
a low density fluid at rest. The computational domain is the unit circle with radius $R=1$ and is split into an inner and an external region at radius $R_s=0.1$. The computational mesh  
consists of $N_E=8,596,566$ triangles of characteristic element size $h=1/2200$. Transmissive boundaries are imposed on the external side. The angular velocity $\omega=10$ of the rotor 
is constant for $r < R_s$, where $r=\sqrt{x^2+y^2}$ denotes the radial coordinate. The initial magnetic field $\mathbf{B}=(2.5,0,0)^T$ as well as the initial  pressure $p=1$ are 
constant in the entire domain and the divergence cleaning velocity of the GLM approach is set to $c_h=2$, while the ratio of specific heats is chosen to be $\gamma=1.4$. The final 
simulation time is $t_f=0.25$. In the inner region, i.e. for $r<R_s$, the initial fluid density is $\rho=10$, while we set $\rho=1$ for $r>R_s$. According to \cite{BalsaraSpicer1999}  
we adopt a linear taper bounded by $0.1 < r \leq 0.105$ for the velocity and the density field in order to smear the initial discontinuity located at radius $r=R_s$. We use the third 
order version ($M=2$) of our ADER-CWENO scheme with the Rusanov flux \eqref{eqn.rusanov} to run the MHD rotor problem. The numerical results are depicted in Figure \ref{fig.MHDRotor}, 
where we plot the distribution of density, fluid pressure, magnetic pressure as well as the Mach number. One can note an overall good qualitative agreement with the solutions presented 
for example in \cite{BalsaraSpicer1999,Dumbser2007204,MHDdivFree2015,DGLimiter2}. 

\subsubsection{The Orszag-Tang vortex system}
\label{ssec.orszag}
\begin{figure}[!htbp]
\begin{center}
\begin{tabular}{cc} 
\includegraphics[width=0.44\textwidth]{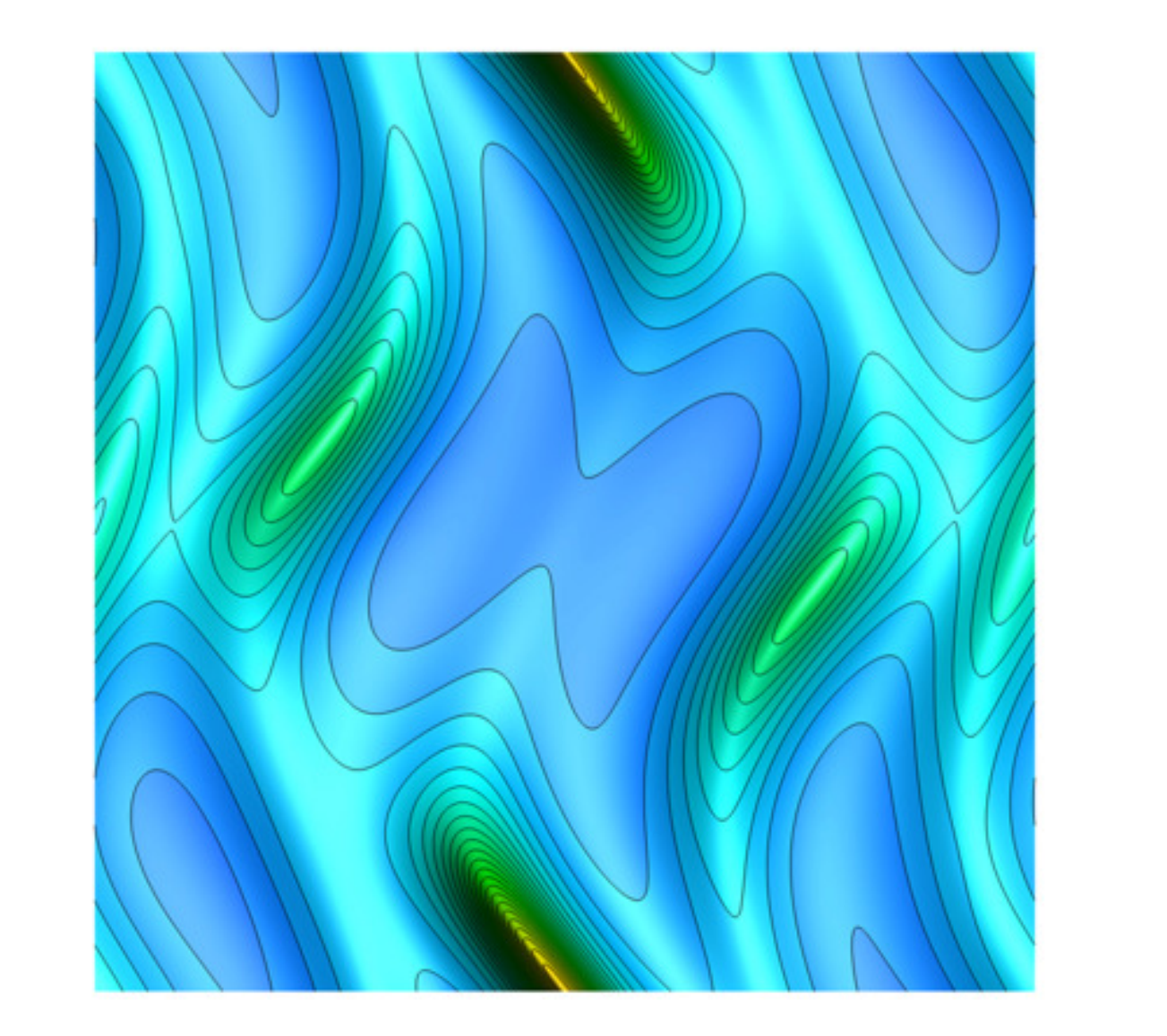}  &              
\includegraphics[width=0.44\textwidth]{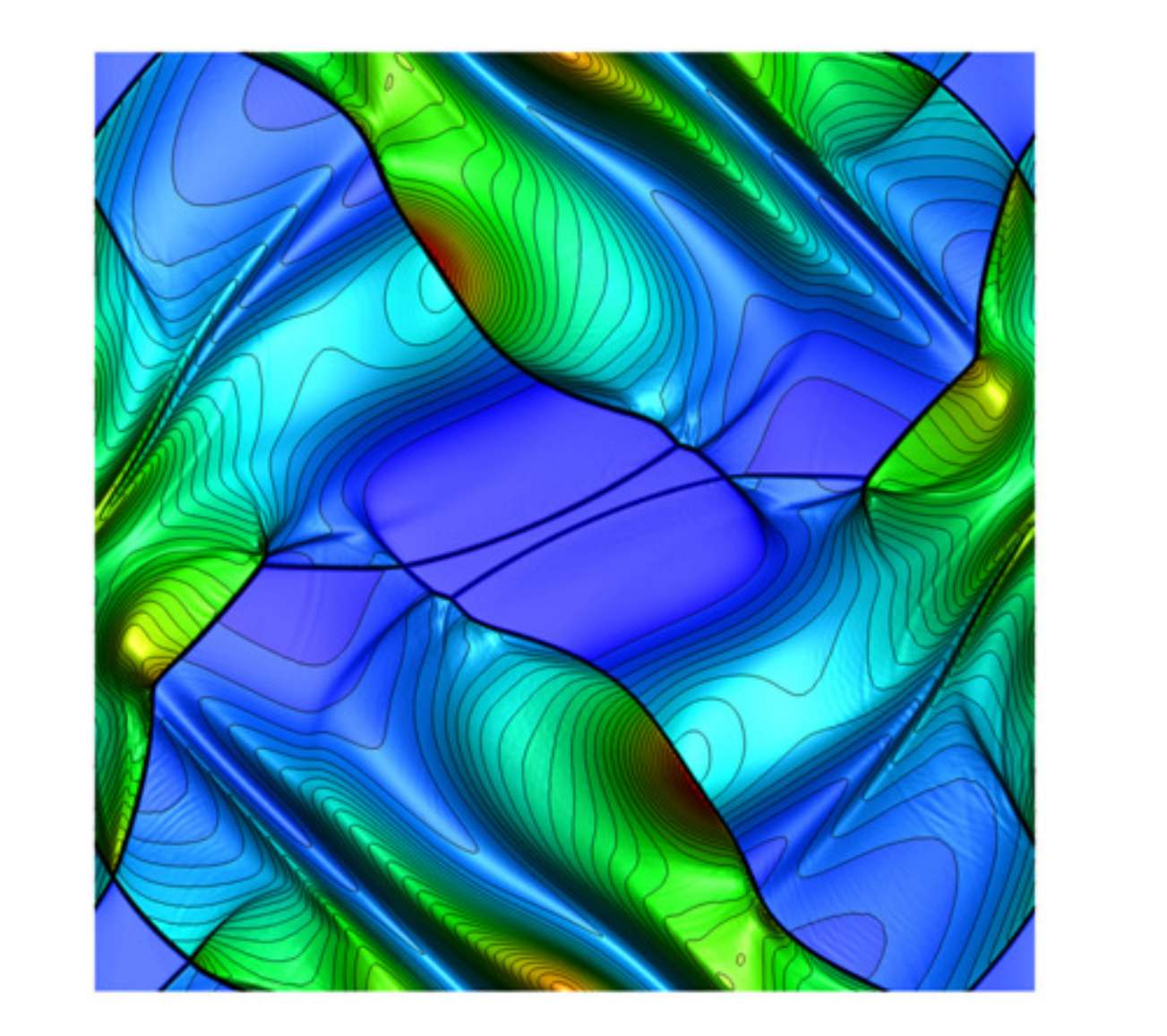} \\              
\includegraphics[width=0.44\textwidth]{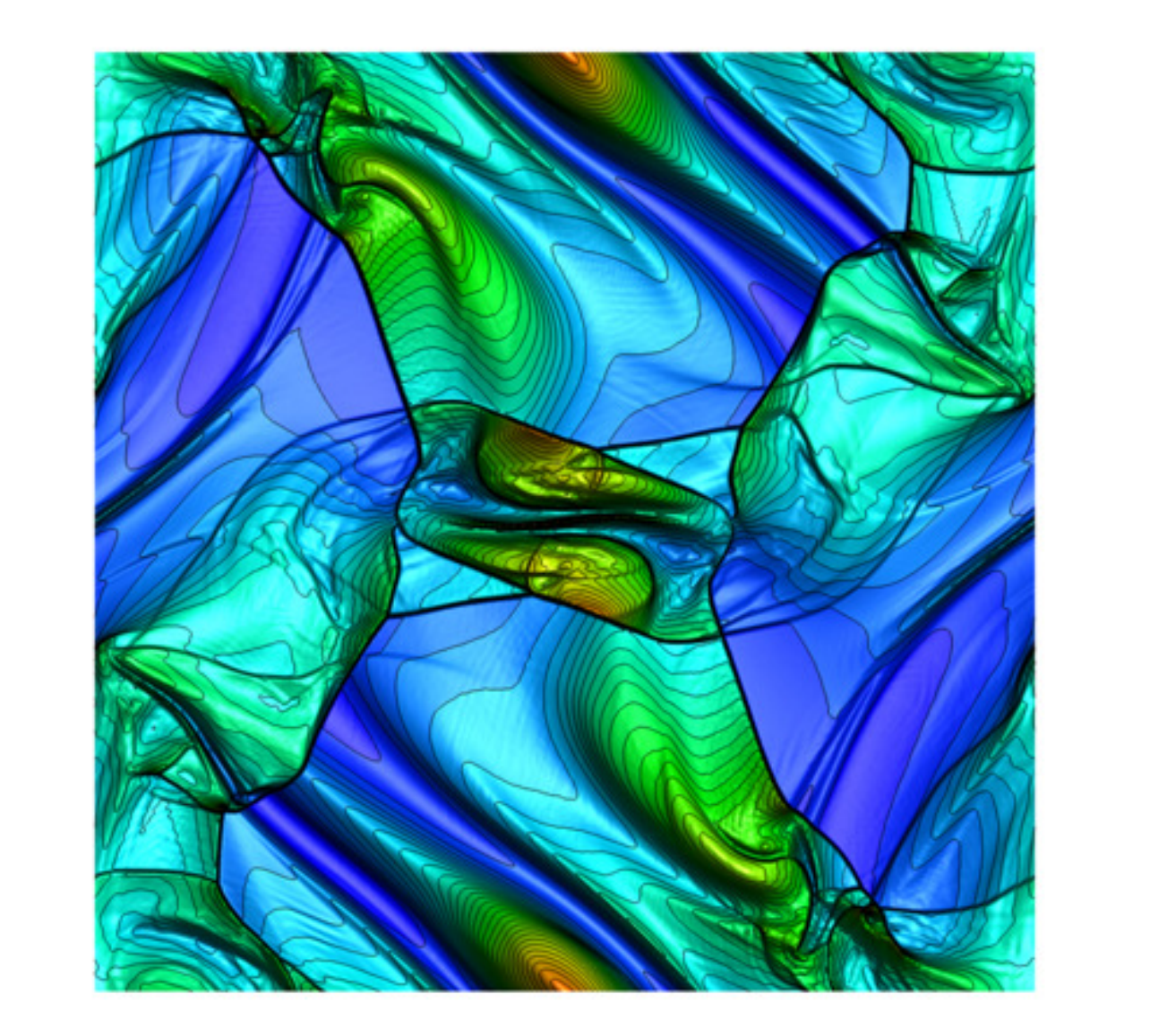}  &              
\includegraphics[width=0.44\textwidth]{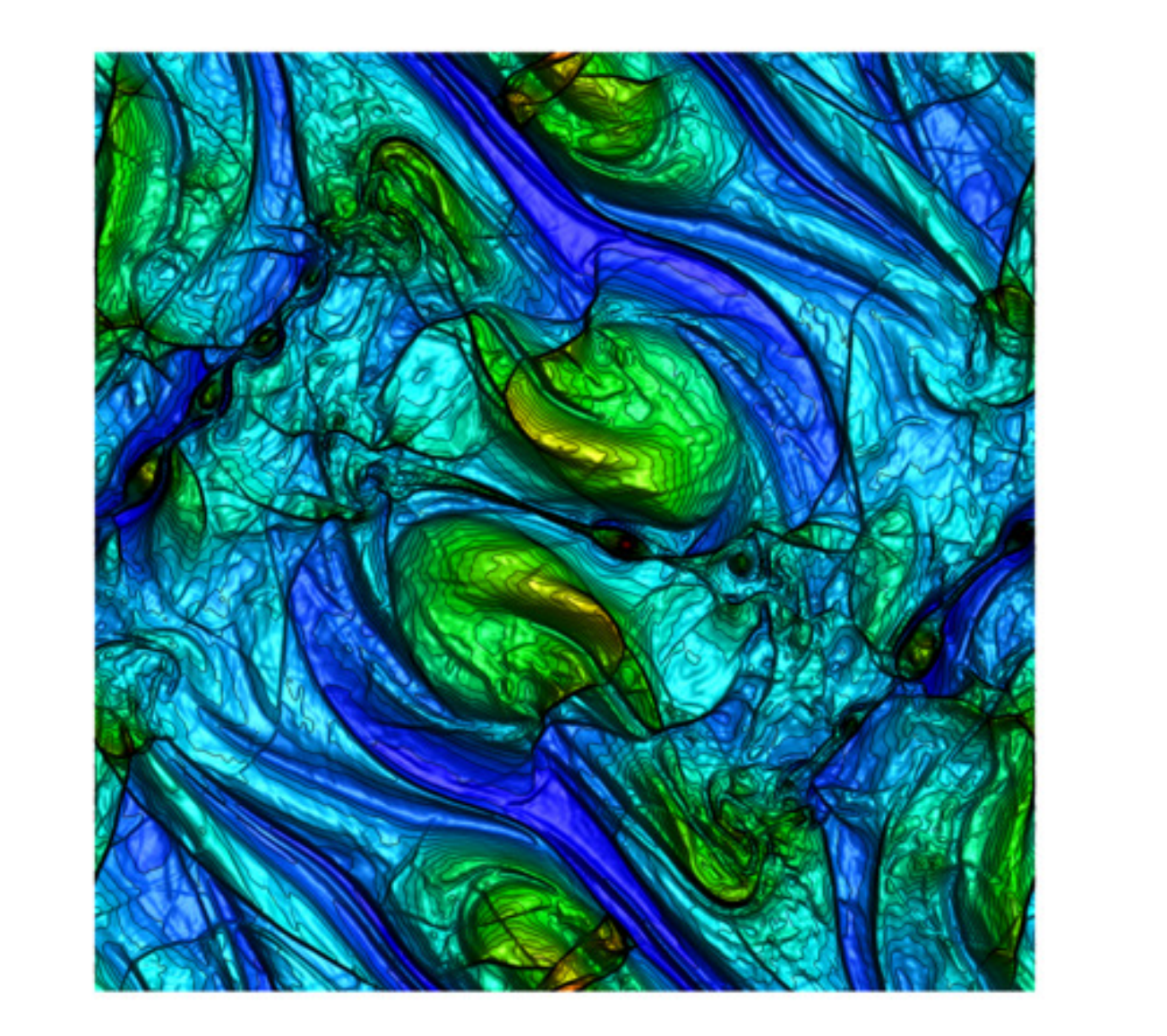}                 
\end{tabular} 
\caption{Pressure contours for the Orszag-Tang vortex problem obtained with a fourth order unstructured ADER-CWENO scheme at 
times $t=0.5$, $t=2.0$, $t=3.0$ and $t=5.0$.}
\label{fig.Orszag-p}
\end{center}
\end{figure}

We solve now the the vortex system of Orszag-Tang \cite{OrszagTang,DahlburgPicone,PiconeDahlburg} for the ideal MHD equations. The computational domain is the square $\Omega=[0;2\pi]^2$ with periodic boundaries everywhere and it is discretized with a characteristic mesh size of $h=2 \pi / 800$ with a total number of $N_E=1,453,240$ triangular elements. 
The initial condition in terms of primitive variables reads 
\begin{equation}
\U = \left( \gamma^2, -\sin(y), \sin(x), 0, \gamma, -\sqrt{4\pi} \sin(y), \sqrt{4\pi} \sin(2x), 0, 0 \right).
\label{eqn.Orszag-IC}
\end{equation}
with $\gamma=5/3$. The divergence cleaning speed is set to $c_h=2$ and the final time of the simulation is taken to be $t_f=5$. The numerical results have been obtained with $M=3$ and the Rusanov-type scheme \eqref{eqn.rusanov} and the evolution of the pressure distribution at output times $t=0.5$, $t=2$, $t=3$ and $t=5$ is shown in Figure \ref{fig.Orszag-p}. 33 equidistant pressure contours
are plotted in the interval $p \in [0.6,7]$. A good qualitative agreement with the solutions provided in \cite{Dumbser2007204,MHDdivFree2015,DGLimiter2} can be noted. The loss of symmetry of the solution
at the final time $t=5$ is due to the use of a non-symmetric unstructured mesh and the physically unstable nature of the flow, since for large times the spatial scales become very small and there is the onset of MHD turbulence. 

\subsection{Numerical results on moving meshes (Lagrangian schemes)}
\label{ssec.ALE}
In order to assess the robustness and the accuracy of the new ADER-CWENO method in the ALE framework on moving meshes, we have chosen a set of benchmark problems that involve very strong shock 
waves and sharp discontinuities in the context of the Euler equations of compressible gas dynamics. All simulations have been run using a fourth order CWENO reconstruction with $M=3$ and the 
Rusanov-type flux \eqref{eqn.rusanov}. 

\subsubsection{The Sedov problem}
\label{ssec.Sedov}
\begin{figure}[!htbp]
\begin{center}
\begin{tabular}{cc} 
\includegraphics[width=0.44\textwidth]{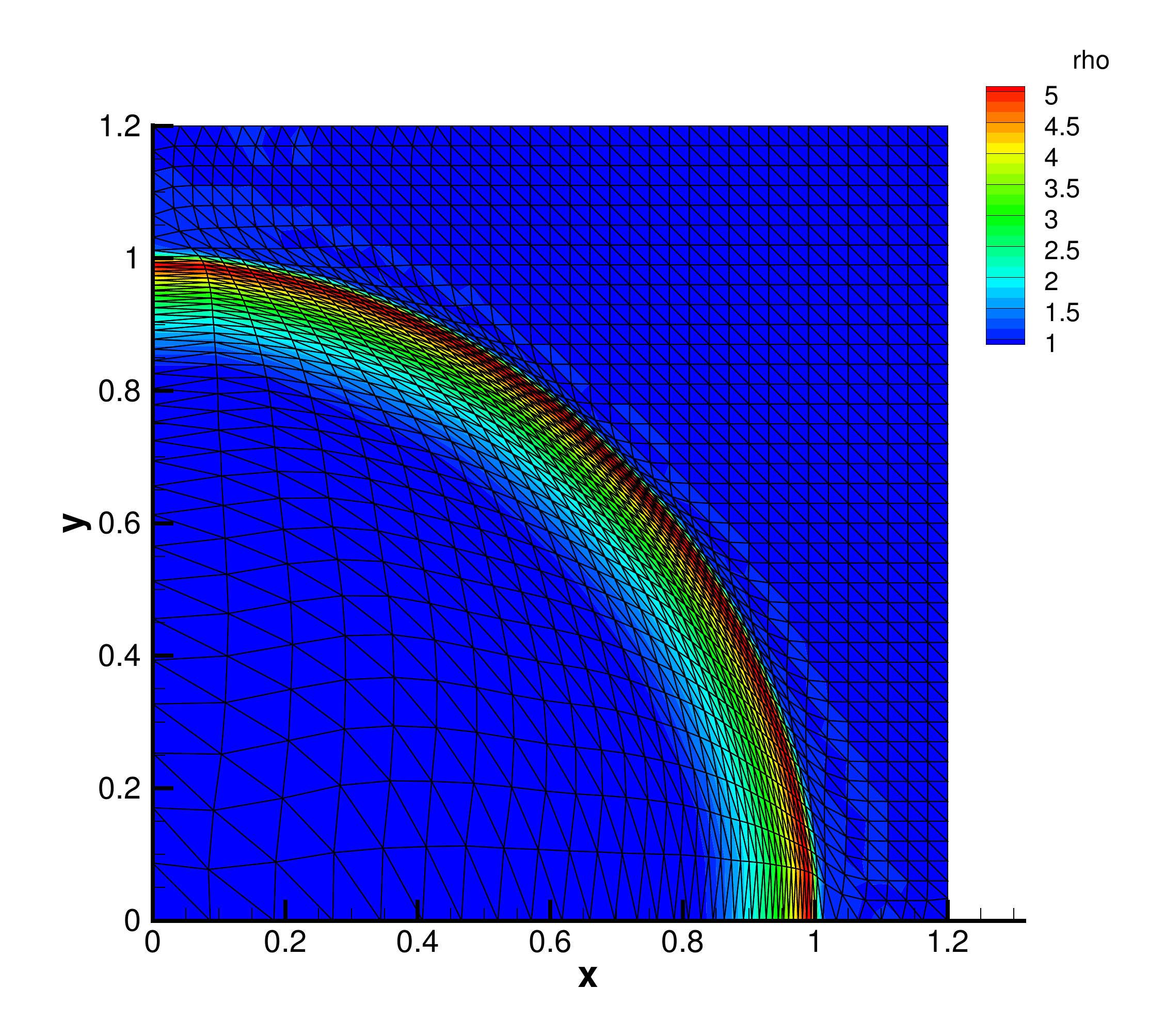}  &           
\includegraphics[width=0.44\textwidth]{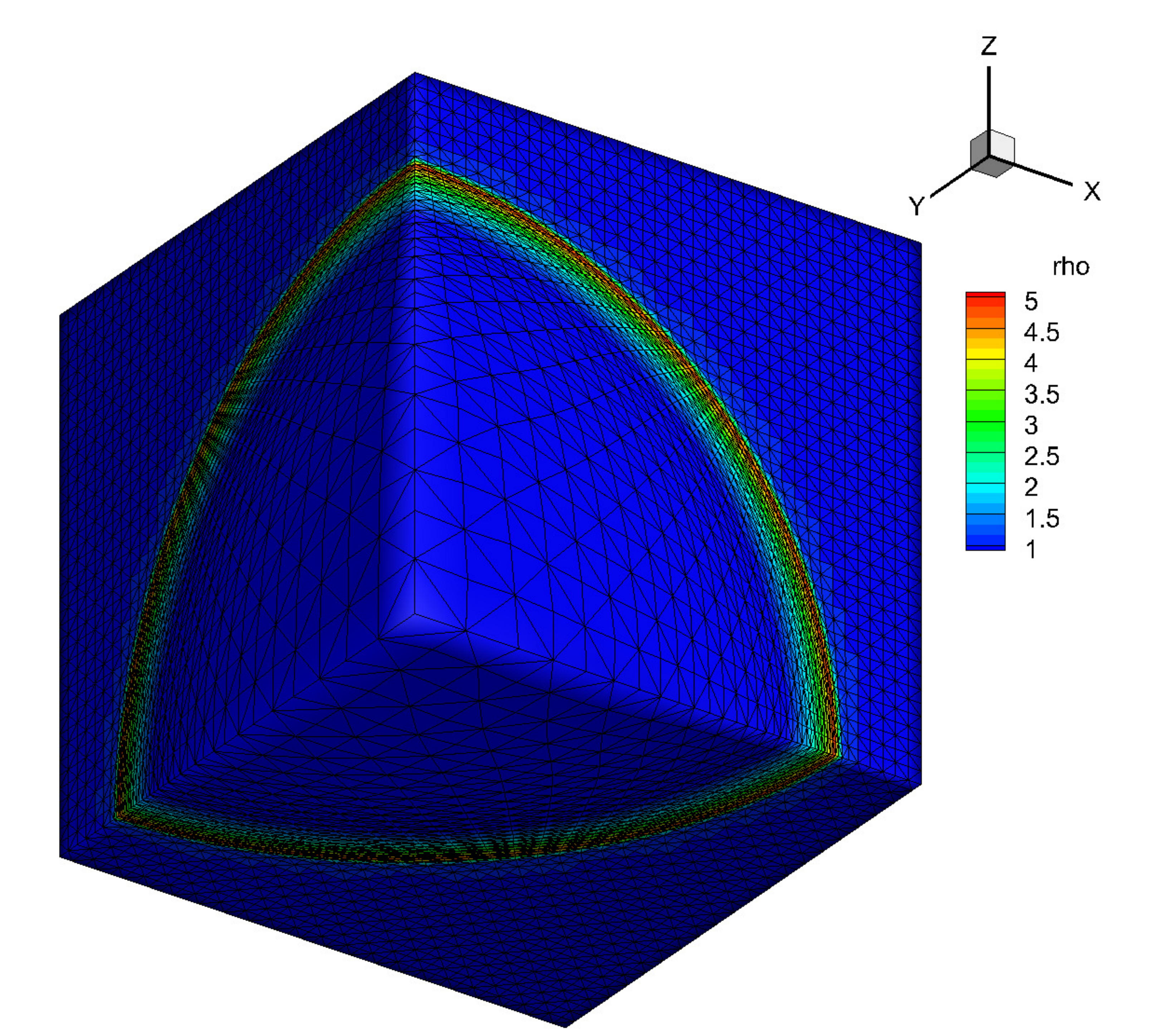}  
\end{tabular} 
\caption{Density distribution and mesh configuration at the final time $t_f=1$ of the Sedov problem in 2D (left) and in 3D (right).}
\label{fig.Sedov-rho}
\end{center}
\end{figure}

\begin{figure}[!htbp]
\begin{center}
\begin{tabular}{cc} 
\includegraphics[width=0.44\textwidth]{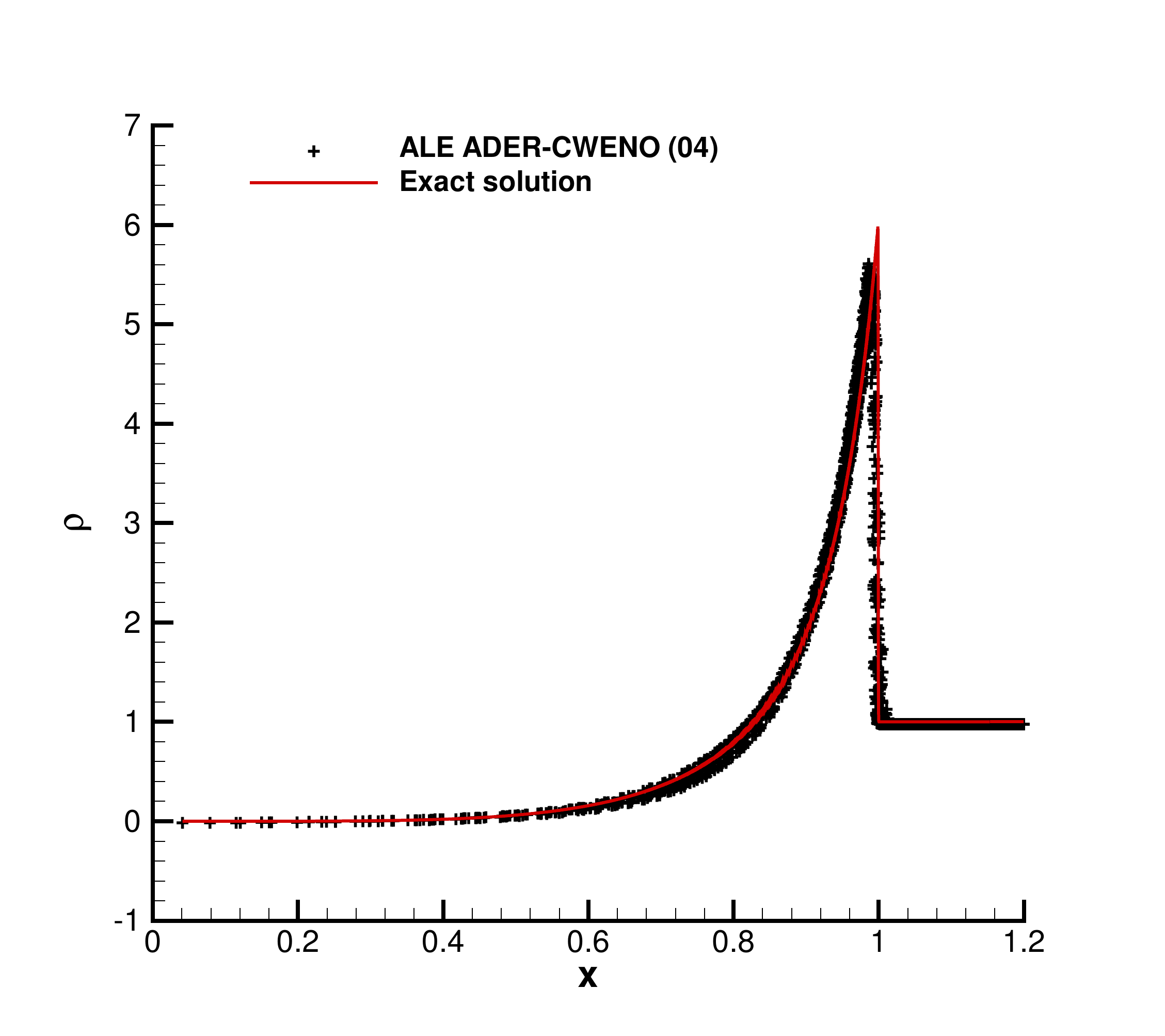}  &           
\includegraphics[width=0.44\textwidth]{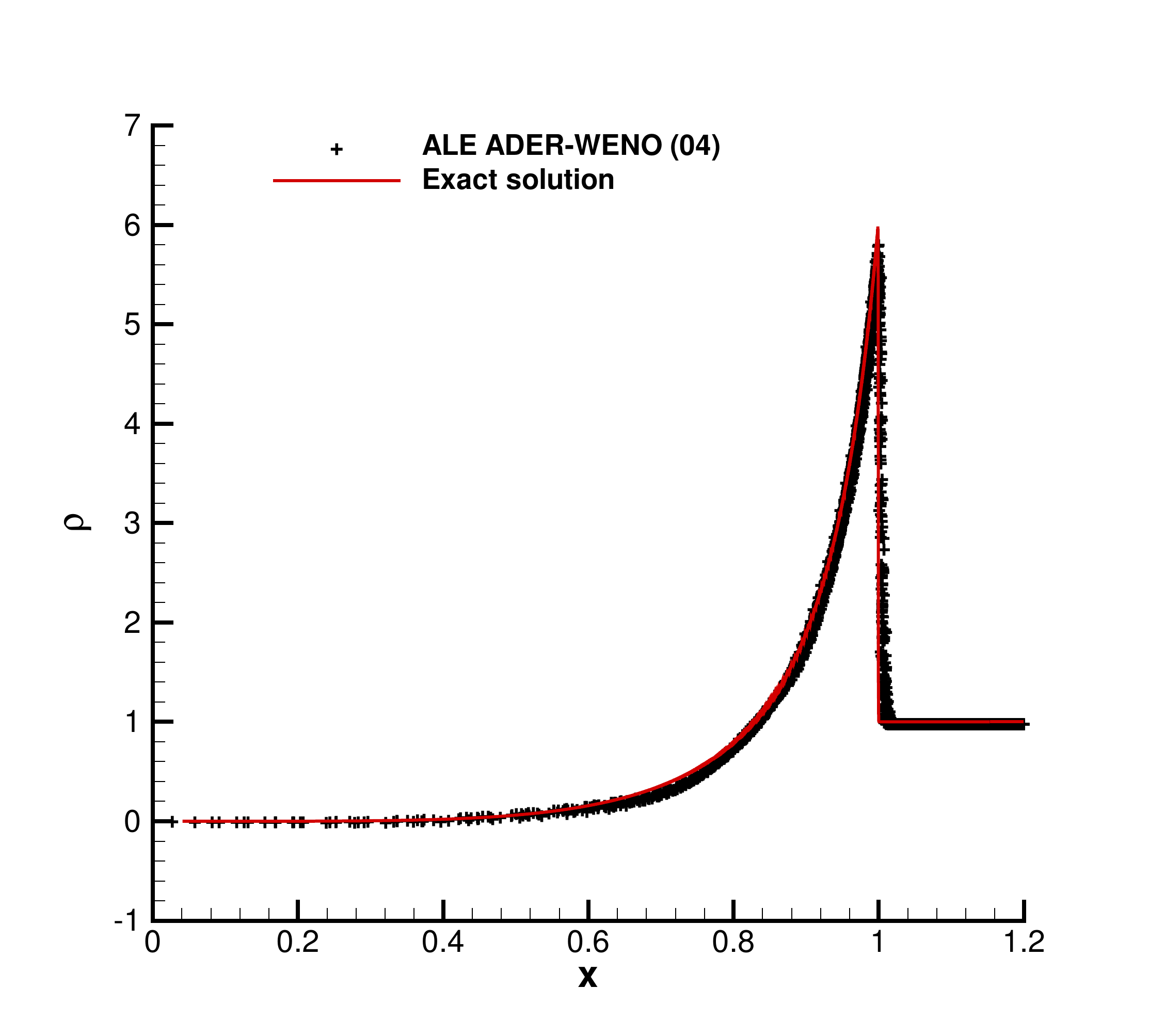} \\
\includegraphics[width=0.44\textwidth]{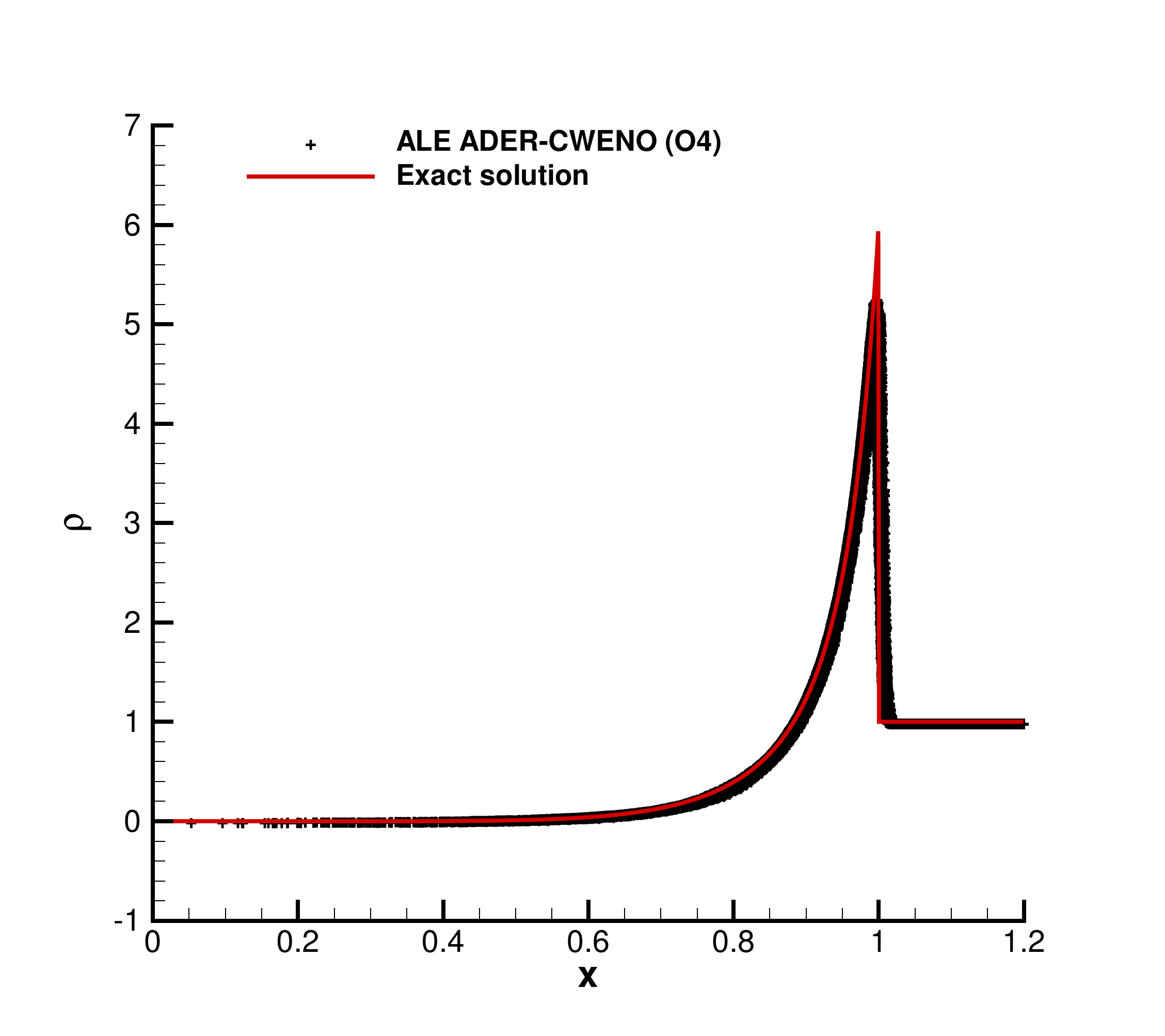}  &           
\includegraphics[width=0.44\textwidth]{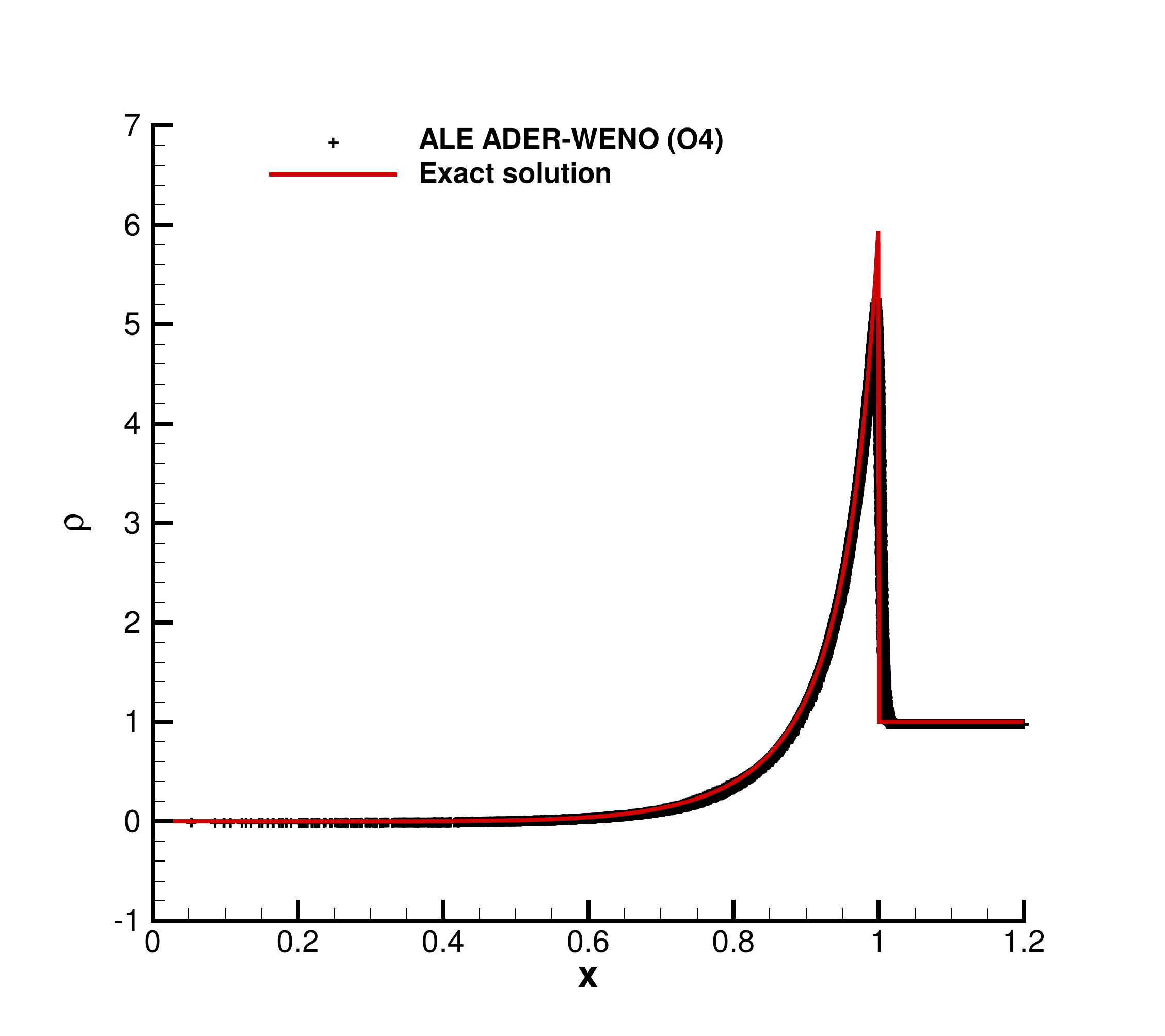}  
\end{tabular} 
\caption{Scatter plot of the cell density as a function of cell radius versus the exact solution for $d=2$ (top row) and $d=3$ (bottom row). The left column reports the results obtained with the ADER-CWENO schemes discussed in this paper, while the right column shows the same plots computed with the original ADER-WENO formulation \cite{Lagrange3D}.}
\label{fig.Sedov-scatter}
\end{center}
\end{figure}

This test problem is widespread in the literature \cite{Maire2009,Maire2009b,LoubereSedov3D} and it describes the evolution of a blast wave that is generated at the origin $\mathbf{O}=(x,y,z)=(0,0,0)$ of the computational domain $\Omega(0)=[0;1.2]^d$. An exact solution based on self-similarity arguments is available from \cite{Sedov,SedovExact} and the fluid is assumed to be an ideal gas with $\gamma=1.4$, which is initially at rest and assigned with a uniform density $\rho_0=1$. The initial pressure is $p_0=10^{-6}$ everywhere except in the cell $c_{or}$ containing the origin $\mathbf{O}$ where it is given by
\begin{equation}
p_{or} = (\gamma-1)\rho_0 \frac{E_{tot}}{\alpha \cdot V_{or}} \quad \textnormal{ with } \quad 
E_{tot} = \left\{ \begin{array}{l} 0.979264 \textnormal{ in 2D} \\
                                   0.851072 \textnormal{ in 3D} \end{array}  \right. .
\label{eqn.p0.sedov}
\end{equation}
$E_{tot}$ is the total energy which is concentrated at $\x=\mathbf{O}$ and the factor $\alpha$ takes into account the cylindrical or spherical symmetry in two or three space dimensions, hence $\alpha=4$ in 2D and $\alpha=8$ in 3D. The computational grid counts a total number of $N_E=3200$ triangles in 2D and $N_E=320000$ tetrahedra in 3D. Figure \ref{fig.Sedov-rho} depicts the density distribution as well as the mesh configuration at the final time of the simulation $t_f=1$, while Figure \ref{fig.Sedov-scatter} contains a scatter plot of the cell density as a function of cell radius versus the exact solution where all cells are represented. A comparison between the results obtained with the ADER-CWENO scheme and the original ADER-WENO formulation proposed in \cite{Lagrange3D} is shown and an excellent agreement with the exact solution is achieved. 

\subsubsection{The Saltzman problem}
\label{ssec.Saltzman}
\begin{figure}[!htbp]
\begin{center}
\begin{tabular}{cc} 
\includegraphics[width=0.55\textwidth]{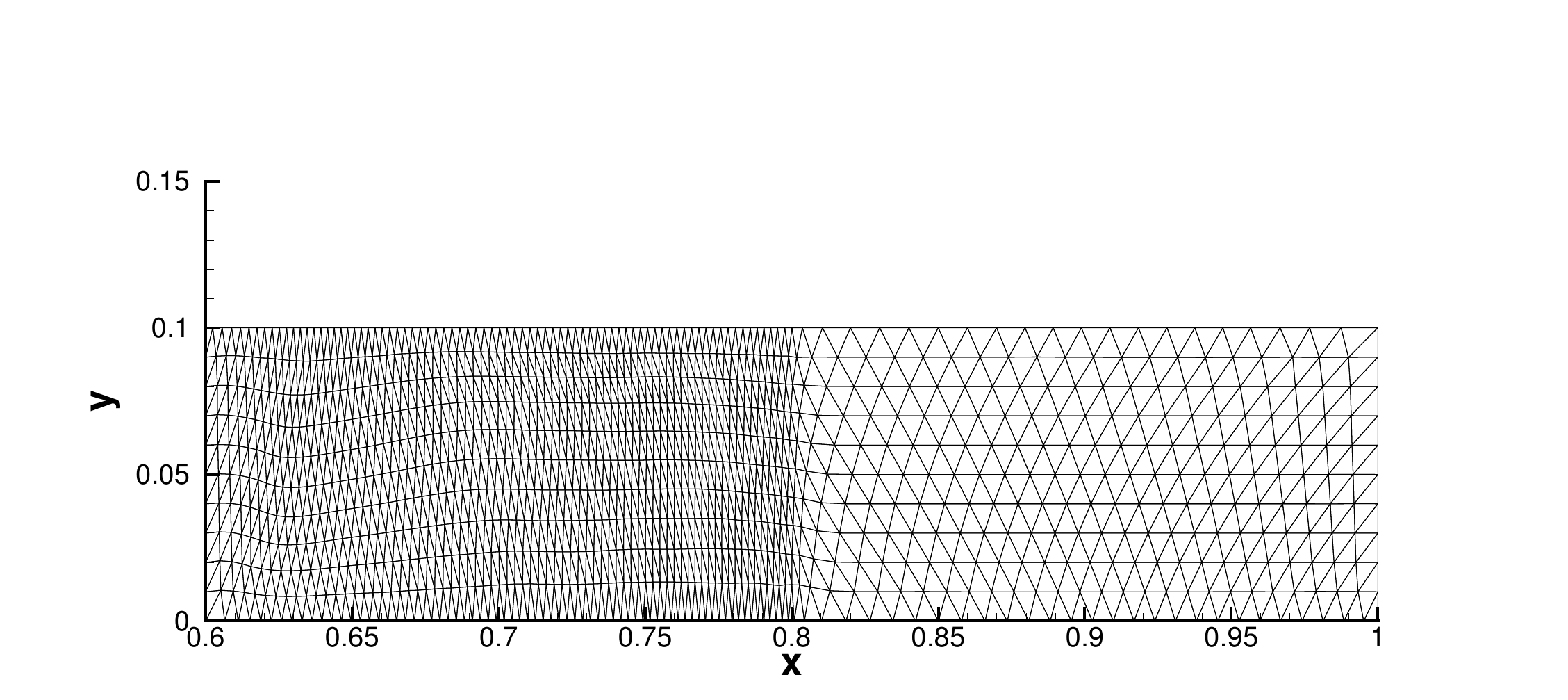}  &           
\includegraphics[width=0.35\textwidth]{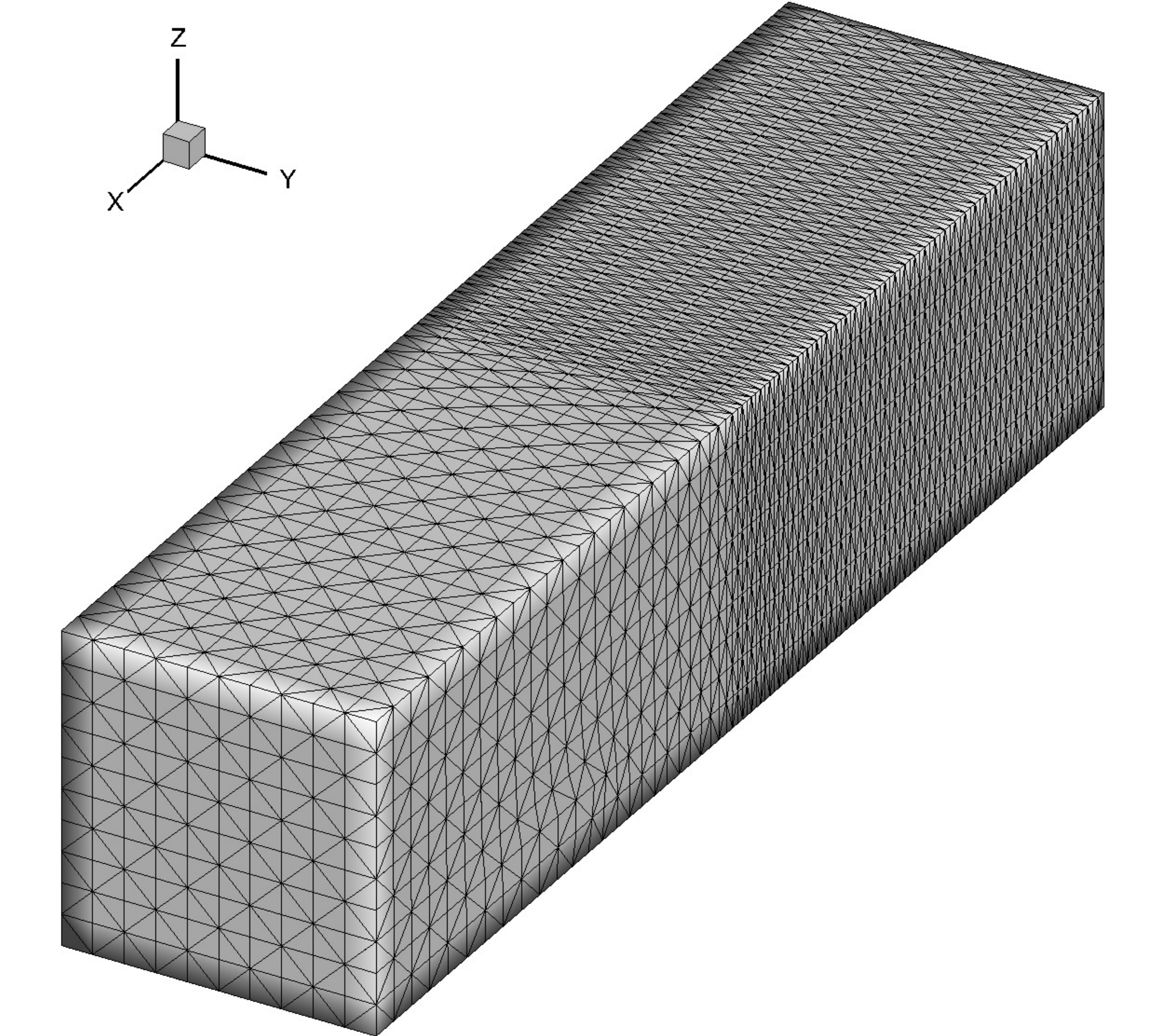}  
\end{tabular} 
\caption{Mesh configuration at the final time $t_f=0.6$ of the Saltzman problem in 2D (left) and in 3D (right).}
\label{fig.Saltzman-mesh}
\end{center}
\end{figure}

\begin{figure}[!htbp]
\begin{center}
\begin{tabular}{cc} 
\includegraphics[width=0.44\textwidth]{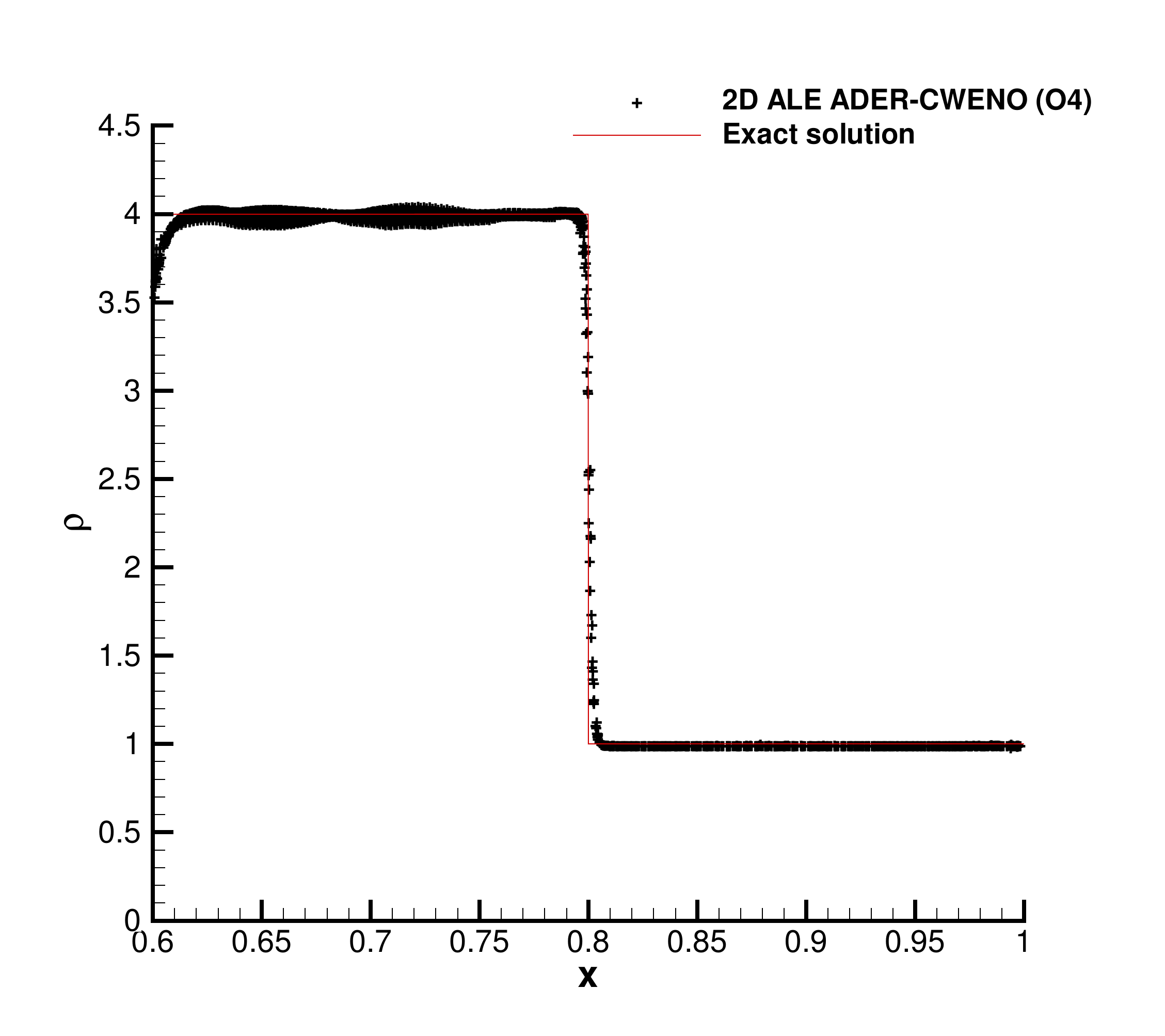}  &           
\includegraphics[width=0.44\textwidth]{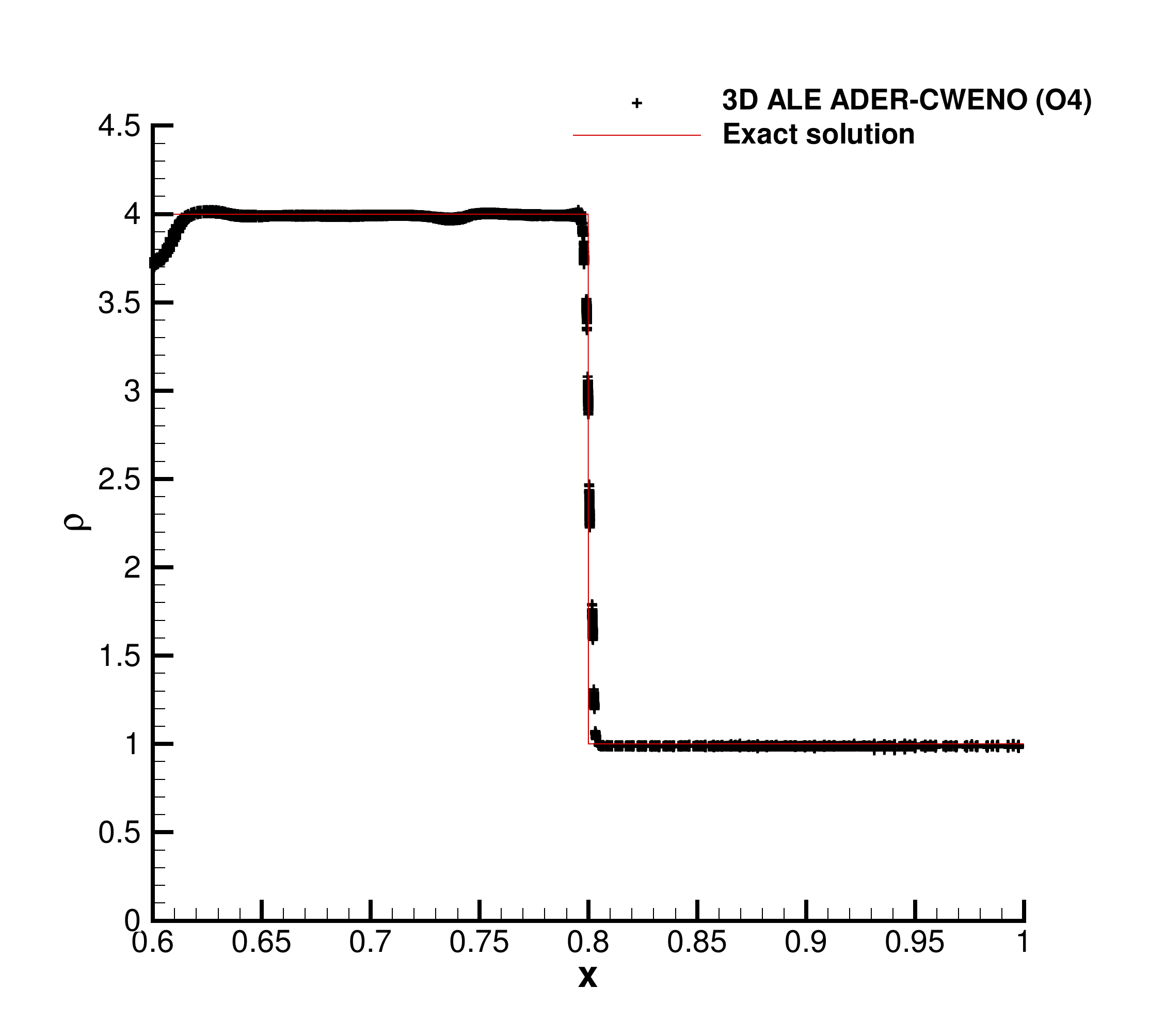} \\
\includegraphics[width=0.44\textwidth]{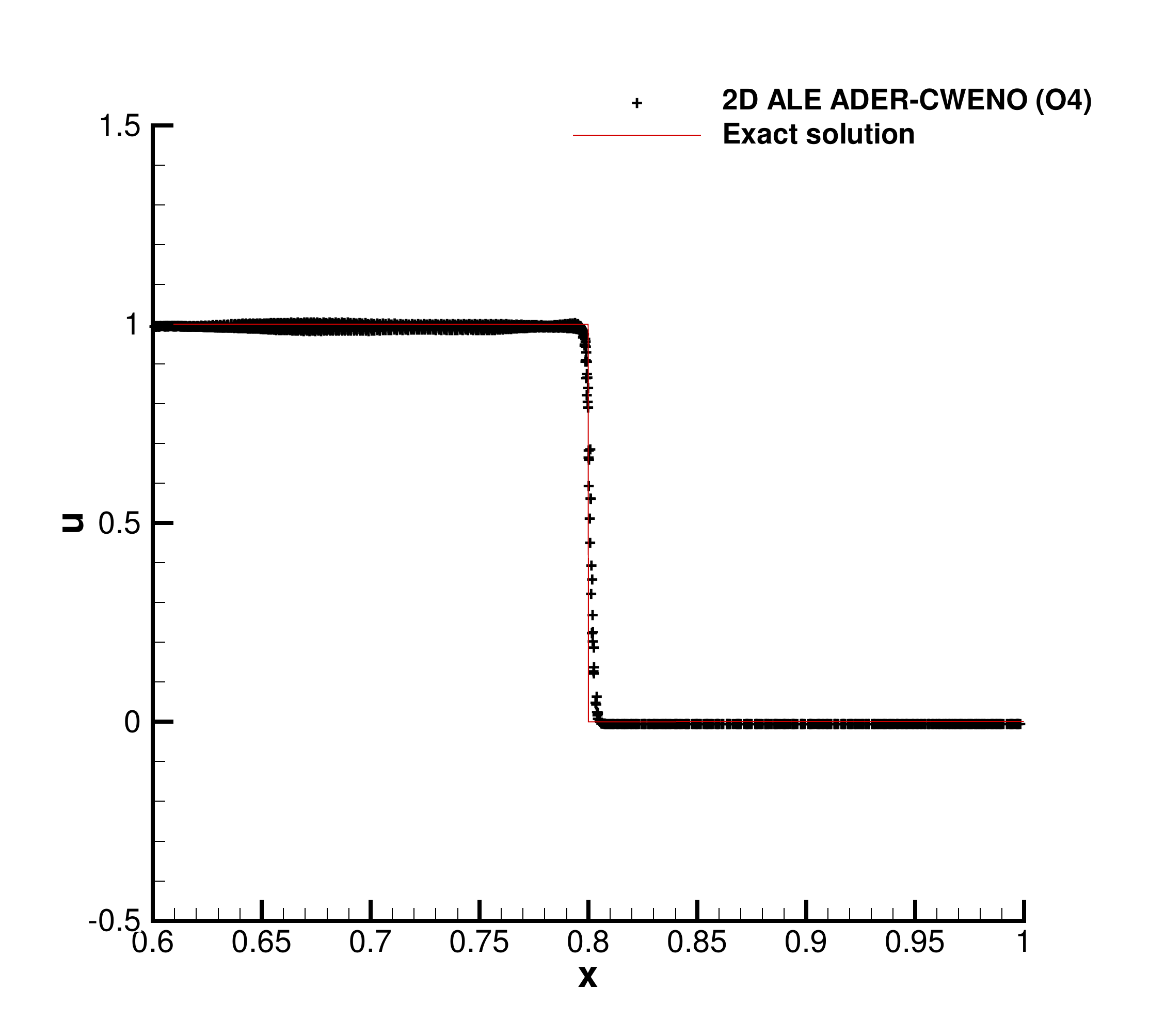}  &           
\includegraphics[width=0.44\textwidth]{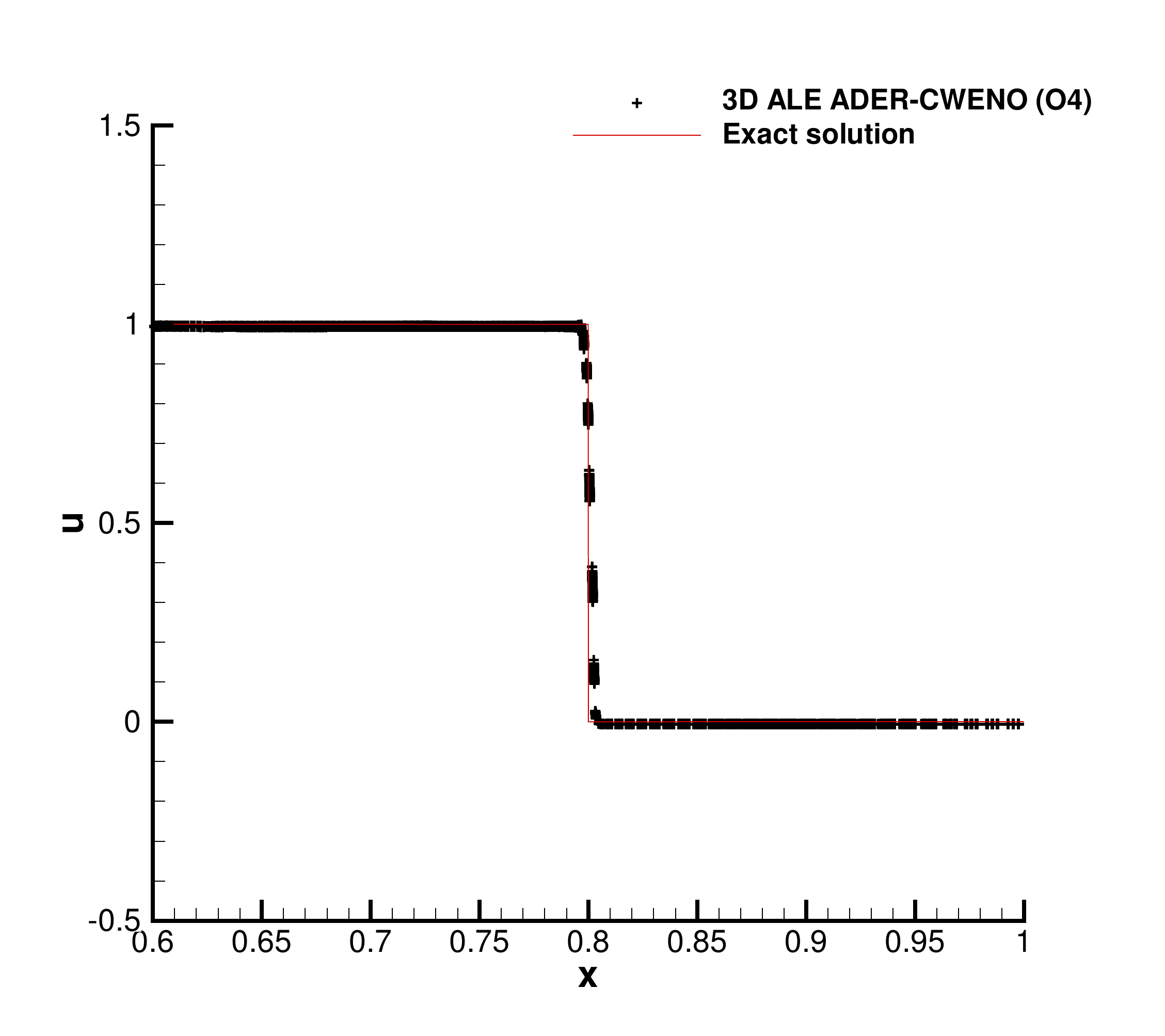}  
\end{tabular} 
\caption{Scatter plot of the cell density (top row) and cell velocity (bottom row) as a function of cell radius versus the exact solution for $d=2$ (left column) and $d=3$ (right column).}
\label{fig.Saltzman-scatter}
\end{center}
\end{figure}

In this problem \cite{SaltzmanOrg,SaltzmanOrg3D} a piston is traveling along the $x$ direction of the computational domain, given by $\Omega(0)=[0;1]\times[0;0.1]$ in 2D and $\Omega(0)=[0;1]\times[0;0.1]\times[0;0.1]$ in 3D. A shock wave is generated that propagates faster than the piston, thus highly compressing the control volumes of the grid. We set no-slip wall boundaries everywhere apart from the left boundary face which is a right-moving piston with velocity $\v_p=(1,0,0)$. The computational domain is first discretized using a characteristic mesh size of $h=1/100$ with quadrilateral and hexahedral elements that are subsequently split into two triangles and five tetrahedra, respectively, hence obtaining $N_E=2000$ in 2D and $N_E=50000$ in 3D. According to \cite{SaltzmanOrg3D,Caramana_CurlQ,Maire2009b} the mesh is initially skewed trough the transformation $\mathcal{R}$ in such a way that the faces and edges are in general not aligned with the main flow:
\begin{equation}
 \begin{array}{l}
 \mathcal{R}^{2D}=\left\{ \begin{array}{l}
  x' = x + \left( 0.1 - y \right) \sin(\pi x) \\
  y' = y
	\end{array}  \right. ,
	\\
	\\
	\mathcal{R}^{3D}=\left\{ \begin{array}{lcl}
  x' = x + \left( 0.1 - z \right) \left( 1 - 20y \right) \sin(\pi x) 	& \textnormal{for} & 0 \leq y\leq 0.05 \\
  x' = x + z\left( 20y - 1 \right) \sin(\pi x) 												& \textnormal{for} & 0.05 < y \leq 0.1 \\
  y' = y & & \\
	z' = z & &
	\end{array}  \right. .
	\end{array}
\end{equation}
The gas is initially at rest with density $\rho_0=1$ and pressure $p_0=10^{-4}$, the ratio of specific heats is $\gamma=5/3$ and the final time of the simulation is $t_f=0.6$, as done in \cite{chengshu2}. The exact solution can be computed by solving a one-dimensional Riemann problem \cite{ToroBook} and it is given by a post shock density of $\rho_e = 4.0$ with the shock front located at $x=0.8$. The final mesh configuration is depicted in Figure \ref{fig.Saltzman-mesh}, while in Figure \ref{fig.Saltzman-scatter} we show a scatter plot of cell density and velocity compared with the exact solution which is very well recovered by our ALE ADER-CWENO method. 

\subsubsection{The Noh problem}
\label{ssec.Noh}
\begin{figure}[!htbp]
\begin{center}
\begin{tabular}{cc} 
\includegraphics[width=0.44\textwidth]{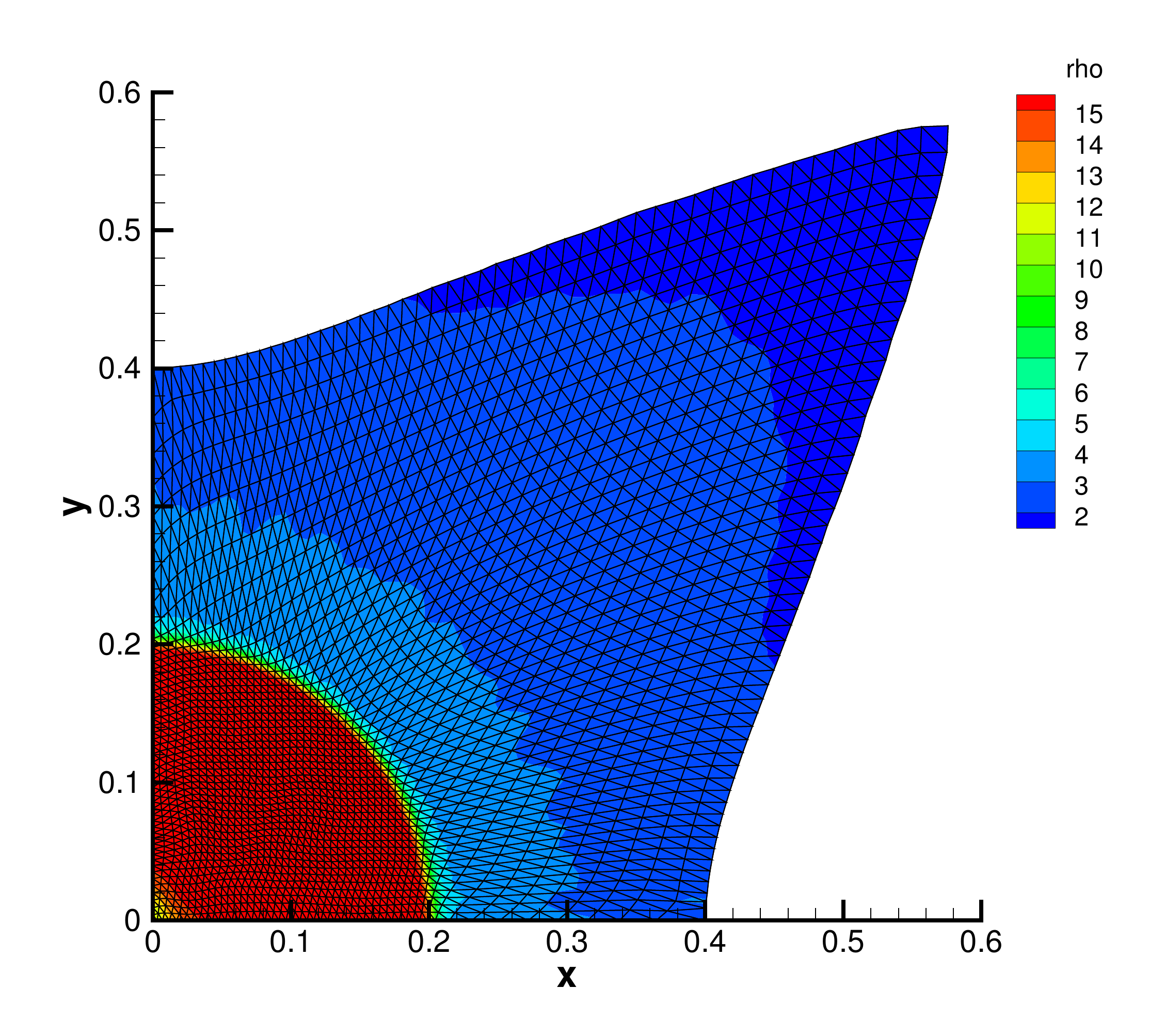}  &           
\includegraphics[width=0.44\textwidth]{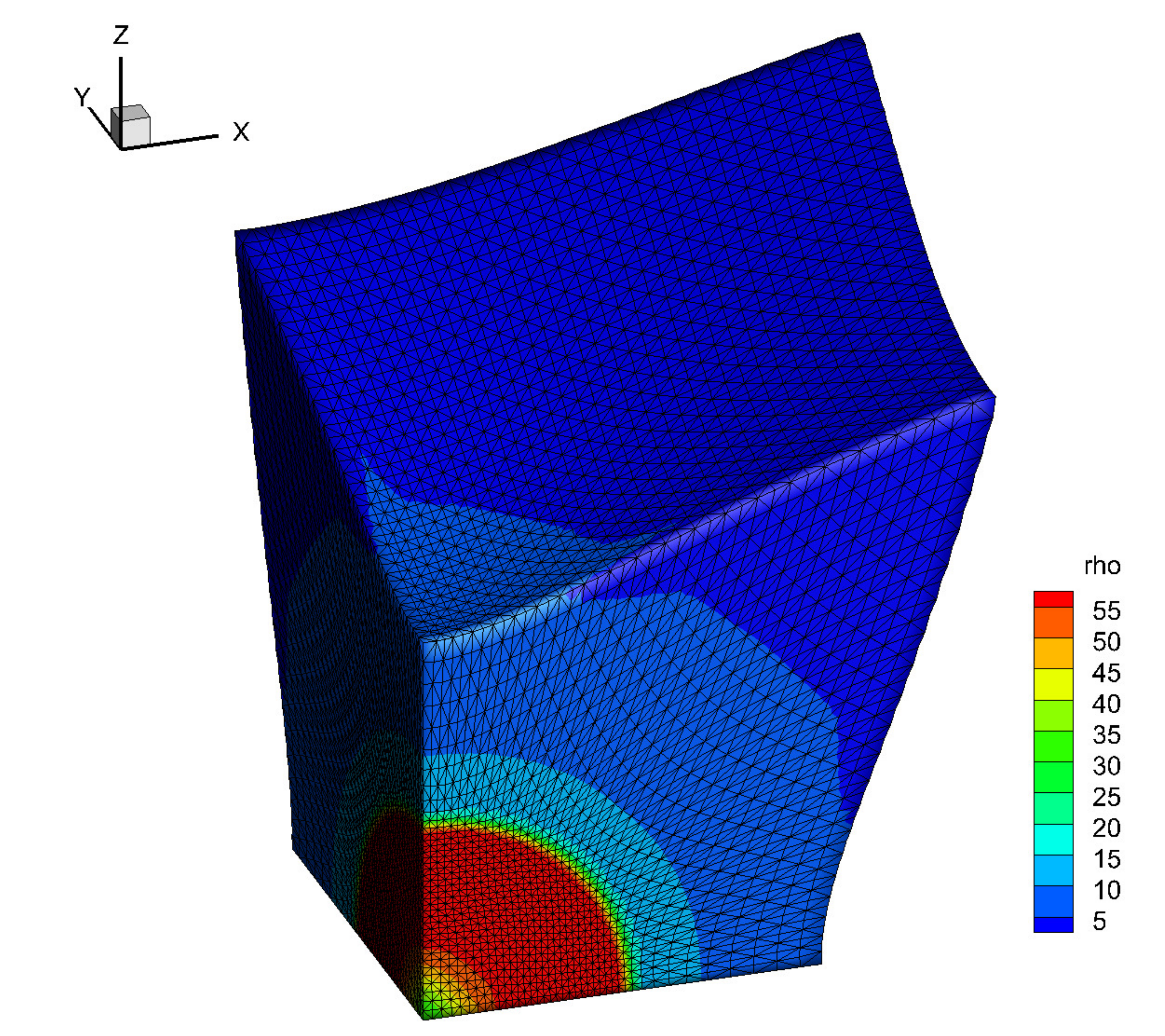} \\
\includegraphics[width=0.44\textwidth]{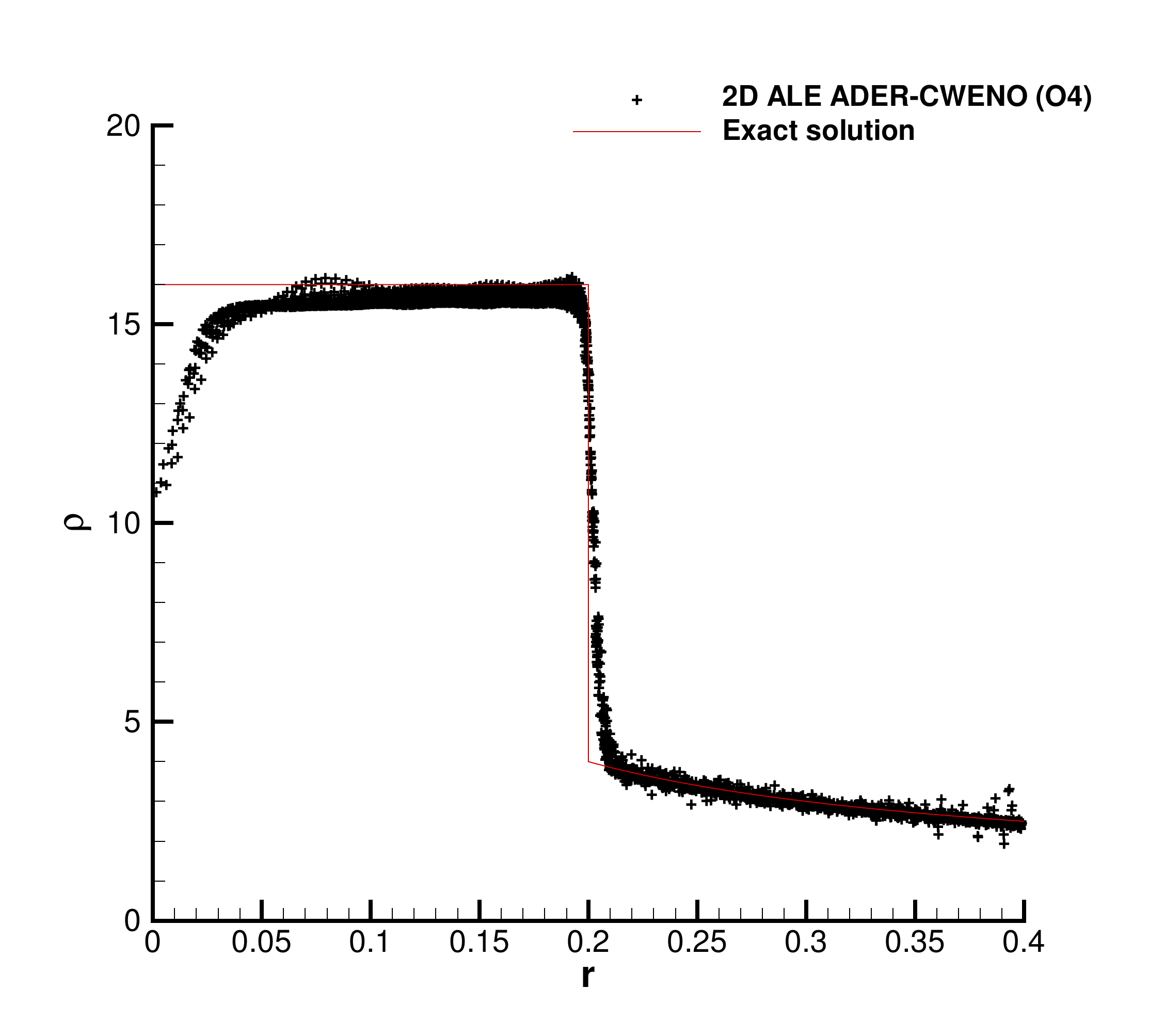}  &           
\includegraphics[width=0.44\textwidth]{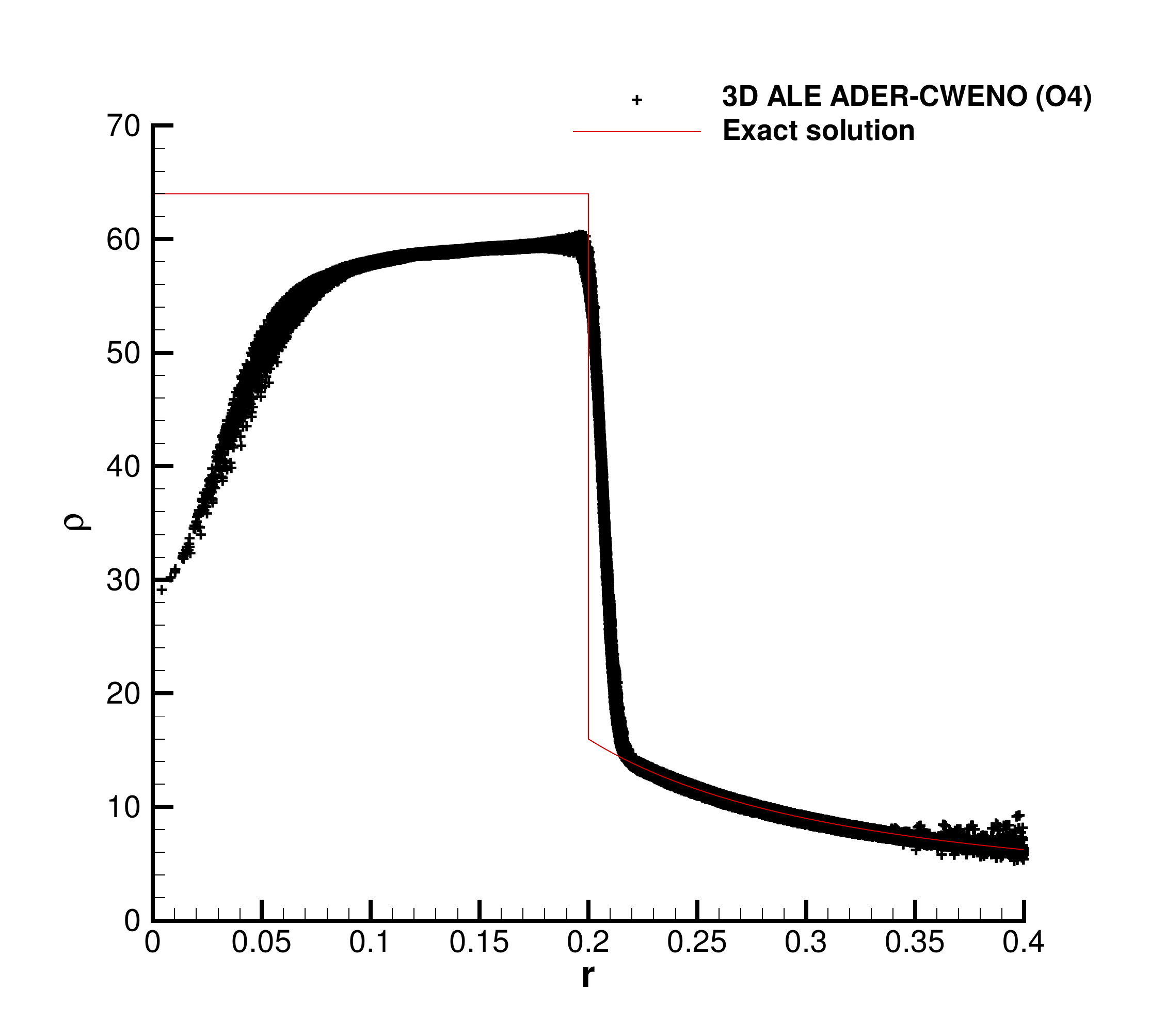}  
\end{tabular} 
\caption{Scatter plot of the cell density (top row) and mesh configuration with density distribution (bottom row) at output time $t_f=0.6$ in 2D (left column) and in 3D (right column).}
\label{fig.Noh}
\end{center}
\end{figure}

The Noh problem \cite{Noh} is a challenging test case used to validate moving mesh schemes in the regime of strong shock waves. The initial computational domain is the box $\Omega(0)=[0;1]^d$ which is filled with a perfect gas with $\gamma=5/3$. The initial values for density and pressure are $\rho_0=1$ and $p_0=10^{-6}$ everywhere, while the velocity field is given by $\v=-\frac{\x}{r}$ with 
$r=\sqrt{\mathbf{x}^2}$. A radial velocity field with $\left\| \mathbf{v} \right\|=1$ is moving the gas towards the origin of the domain $O=(0,0,0)$, where a circular (spherical) shock wave is generated and starts growing and expanding in the outward radial direction. As usual \cite{Maire2009,Maire2009b} the computational domain is assumed to be only a part of the entire cylindrical or spherical one,  therefore those sides which coincide with an axis of symmetry of the complete domain are addressed as internal faces and are given a no-slip wall condition, while the remaining sides are referred to as external faces where moving boundaries are imposed. The final time is chosen to be $t_f=0.6$, therefore according to the exact solution \cite{Noh} the shock wave is located at radius $R=0.2$ and the maximum density value is $\rho_f=16$ in 2D and $\rho_f=64$ in 3D, which occurs on the plateau behind the shock wave. Numerical results are shown in Figure \ref{fig.Noh} where a scatter plot of cell density compared with the analytical solution is plot as well as the initial and final mesh configuration with the corresponding density distribution. The cylindrical or spherical symmetry of the shock wave is very well preserved and located at the correct position. Close to the origin of the domain one notes the classical wall heating effect \cite{Rider:2000}. For the three-dimensional case the post-shock density is underestimated, but this result is in agreement with what is typically observed in the literature, see \cite{Dobrev2012,3Dgcl_Maire2016}. 

\subsubsection{Computational efficiency of the ALE ADER-CWENO algorithm}
\label{ssec.Efficiency}
In order to evaluate the benefits brought by the use of the new and more efficient central WENO (CWENO) reconstruction procedure introduced in Section \ref{ssec.cweno}, we have measured the total  computational time needed for running some test cases in multiple space dimensions. The results are compared with the original ADER-WENO schemes developed in \cite{Lagrange3D} and the 
results of the profiling are reported in Table \ref{tab.efficiency}. Let $N_E$ represent the total number of elements of the computational mesh and $N$ be the number of time steps needed to carry out
the simulation until the final time. Let furthermore $t_{CPU}$ denote the total computational time of the simulation measured in seconds and $\tau_E=\frac{t_{CPU}}{N_E \cdot N}$ be the time used per  element update. The final efficiency ratio is evaluated as 
\begin{equation}
\beta={\tau_E^{AW}}/{\tau_E^{ACW}},
\label{eqn.efficiency}
\end{equation}
where the ALE ADER-WENO scheme is referred to as \quotew{AW} and the algorithm presented in this article is addressed with \quotew{ACW}. One can note that the new ADER-CWENO scheme is more than 
two times faster than the original ADER-WENO algorithm. Concerning memory consumption, the unstructured WENO reconstruction schemes developed in \cite{Dumbser2007693,Dumbser2007204} need 
$n_s=1+2(d+1)$ stencils (one central stencil, $d+1$ forward sector stencils and $d+1$ reverse sector stencils), each of which produces a polynomial of degree $M$. In the new unstructured 
CWENO scheme introduced in this paper, we employ only $n_s=1+(d+1)$ stencils (one central stencil and $d+1$ forward sector stencils), but only the central stencil produces a polynomial of 
degree $M$, while the one-sided sectorial stencils build only polynomials of lower degree $M^s=1$. Therefore, the total number of matrix elements $\mathcal{E}$ that needs to be stored for 
the reconstruction matrices is $\mathcal{E}_{\textnormal{WENO}}=\left( 1 + 2(d+1) \right) \cdot d \mathcal{M}(M,d)^2$ for the unstructured WENO schemes used in 
\cite{Dumbser2007693,Dumbser2007204}, while it is only 
$\mathcal{E}_{\textnormal{CWENO}}=1 \cdot d \mathcal{M}(M,d)^2 + (d+1) \cdot (d+1)^2$ for the new unstructured CWENO schemes presented in this paper. For sufficiently high polynomial degrees 
$M$ one can essentially neglect the term $(d+1) \cdot (d+1)^2$ in the expression for $\mathcal{E}_{\textnormal{CWENO}}$. This means that for storing the reconstruction matrices on fixed
meshes, which is the main memory cost for unstructured Eulerian WENO schemes, the CWENO scheme in 2D needs approximately only $1/7$ of the memory of our original 
unstructured WENO schemes presented in \cite{Dumbser2007693,Dumbser2007204}, and it needs approximately only $1/9$ of the memory in 3D, which is a gain of almost one order of magnitude, 
without compromising accuracy, robustness or computational efficiency. 

\begin{table}[!htbp]
  \caption{Computational efficiency of the new ADER-CWENO schemes on moving meshes compared with the original direct ALE ADER-WENO algorithm \cite{Lagrange3D}. $N_E$ represents the total number of elements of the computational mesh, $N$ is the number of time steps needed to carry on the simulation until the final time, $t_{CPU}$ is the total computational time and $\tau_E=t_{CPU} / ( N_E N )$ is the time needed to update one element. Finally $\beta={\tau_E^{W}}/{\tau_E^{CW}}$ indicates the efficiency ratio. All the simulations use piecewise cubic reconstruction in space and time.} 
	\begin{center}
	\begin{small}
		\begin{tabular}{llccccccc}
		\hline
		\multicolumn{9}{c}{\textbf{2D test problems}} \\
		\hline
		Test case   &         & \multicolumn{3}{c}{ALE ADER-CWENO}                         & \multicolumn{3}{c}{ALE ADER-WENO}                         & $\beta$ 		\\
    \hline
								& $N_E$ &  $N$ & $t_{CPU}$            & $\tau_E^{ACW}$        & $N$            & $t_{CPU}$  & $\tau_E^{AW}$        &    \\
		\hline
		Sedov       & 3200  & $\!1292$ & $ 3.0\cdot 10^{3} $	& $ 7.2\cdot 10^{-4} $ & 1350 & $ 4.9\cdot 10^{3} $ & $ 1.1\cdot 10^{-3} $ & 1.6  \\
		Saltzman    & 2000  & $\!1715$ & $ 2.5\cdot 10^{3} $ & $ 7.2\cdot 10^{-4} $ & 1727 & $ 6.2\cdot 10^{3} $ & $ 1.8\cdot 10^{-3} $ & 2.5  \\
		Noh         & 5000  & $\!596$  & $ 2.1\cdot 10^{3} $ & $ 7.1\cdot 10^{-4} $ & 662  & $ 6.3\cdot 10^{3} $ & $ 1.9\cdot 10^{-3} $ & 2.7 \\
		Explosion		& 17340	& $\!298$	 & $ 3.9\cdot 10^{3} $ & $ 6.6\cdot 10^{-4} $ & 302  & $ 5.4\cdot 10^{3} $ & $ 1.0\cdot 10^{-3} $ & 1.6  \\
		\hline
		            &       &      &                      &                       &      &                      &                       & \textbf{2.1} \\
		\multicolumn{9}{c}{} \\						
		\hline
		\multicolumn{9}{c}{\textbf{3D test problems}} \\
		\hline
		Test case & $N_E$   & \multicolumn{3}{c}{ALE ADER-CWENO}                         & \multicolumn{3}{c}{ALE ADER-WENO}                         & $\beta$ 		\\
    \hline
								&       &  $N$ & $t_{CPU}$            & $\tau_E^{ACW}$        & $N$            & $t_{CPU}$  & $\tau_E^{AW}$        &    \\
		\hline
		Sedov       & 320000  & $\!2309$ & $ 4.0\cdot 10^{6} $ & $ 5.3\cdot 10^{-3} $ & 2623 & $ 8.8\cdot 10^{6} $ & $ 1.1\cdot 10^{-2} $ & 2.0  \\
		Saltzman    & 50000   & $\!2120$ & $ 6.0\cdot 10^{5} $ & $ 5.6\cdot 10^{-3} $ & 1934 & $ 1.4\cdot 10^{6} $ & $ 1.4\cdot 10^{-2} $ & 2.5  \\
		Noh         & 320000  & $\!1885$ & $ 3.1\cdot 10^{6} $ & $ 5.1\cdot 10^{-3} $ & 1886 & $ 4.6\cdot 10^{6} $ & $ 7.7\cdot 10^{-3} $ & 1.5 \\
		Explosion		& 1469472	& $\!588$  & $ 1.6\cdot 10^{7} $ & $ 1.9\cdot 10^{-2} $ & 658  & $ 5.6\cdot 10^{7} $ & $ 5.8\cdot 10^{-2} $ & 3.1  \\
		\hline
		            &       &      &                      &                       &      &                      &                       & \textbf{2.3}						
		\end{tabular}
  \end{small}
	\end{center}
	\label{tab.efficiency}
\end{table}

\section{Conclusions}
\label{sec.concl}
  
In this paper we have presented a novel arbitrary high order accurate central WENO reconstruction procedure (CWENO) in order to produce piecewise polynomials in space from known cell averages 
on unstructured simplex meshes in two and three space dimensions. CWENO considers a candidate polynomial of the desired order of accuracy with a centered stencil and a number of sectorial 
polynomials with directionally-biased stencils. A main difference from the original WENO approach is that the sectorial polynomials may have a smaller degree and thus their stencil can be 
chosen inside the stencil of the central optimal polynomial, giving rise to a reconstruction procedure with a very small total stencil. The method presented in this paper is thus much more
\textit{compact} compared to the unstructured WENO reconstructions used in \cite{Dumbser2007693,Dumbser2007204,Lagrange3D}. All the polynomials are then combined in the usual way 
at the aid of non-linear WENO weights in order to guarantee the non-oscillatory properties of the reconstruction. 

In this paper we have employed the new CWENO reconstruction as initial data for a fully-discrete one-step ADER finite volume method of orders up to five in both space and time on fixed 
Eulerian meshes as well as on moving Arbitrary-Lagrangian-Eulerian grids. As a result of the smaller stencils employed, when compared to analogous ADER-ALE schemes initialized with the 
classical WENO reconstruction procedure, the simulations can be completed with a speed-up between 1.5 and 3, depending on the test case. For fixed meshes, the memory consumption is 
almost between 7 and 9 times less. 

Since the ADER approach leads to a fully-discrete one-step scheme, only very few MPI communications are necessary when implementing the method on a massively parallel distributed memory supercomputer.
Only two main communication steps are necessary: the first step for the exchange of the cell averages between the CPUs before the CWENO reconstruction and the second step for the exchange 
of the resulting reconstruction polynomials across MPI domain boundaries. The very compact stencil of the CWENO approach helps to improve MPI efficiency compared to the WENO schemes presented 
in \cite{Dumbser2007693,Dumbser2007204}, 
since the stencil overlap across MPI domains is smaller and thus less data need to be exchanged. Compared with classical Runge-Kutta time stepping, the ADER approach allows this entire procedure
to be performed only \textit{once} per time step, while in Runge-Kutta based schemes, the nonlinear reconstruction and the associated MPI communication must be done in each Runge-Kutta sub-stage
again. We have shown two examples on more than one hundred million elements where the CWENO reconstruction leads to a piecewise polynomial data representation involving \textbf{billions} of degrees  
of freedom. To the knowledge of the authors, these are the largest simulations ever carried out so far with WENO schemes on unstructured meshes. 

In the future we plan to use the present CWENO reconstruction on unstructured meshes also within the novel \textit{a posteriori} subcell finite volume limiters for high order DG schemes 
recently proposed in \cite{DGLimiter1,DGLimiter3,DGLimiter4} on fixed and moving unstructured meshes.  
Future work may also concern applications of the new method to more complex geometries and to large-scale real-world flow problems as they appear in science and engineering, including 
also molecular and turbulent viscosity. 

\section*{Acknowledgments}
The presented research has been financed by the European Research Council (ERC) under the European Union's Seventh Framework 
Programme (FP7/2007-2013) with the research project \textit{STiMulUs}, ERC Grant agreement no. 278267. The authors acknowledge 
the kind support of the Leibniz Rechenzentrum (LRZ) in Munich, Germany, for awarding access to the SuperMUC supercomputer, 
where all the large scale test problems have been run. 

M.S. acknowledges the support of the \quotew{National Group for Scientific Computation (GNCS-INDAM)}.

\bibliographystyle{siamplain}
\bibliography{biblio}
\end{document}